%% file: ergm.estimation.tex
\documentclass[aos,preprint]{imsart}

\RequirePackage[OT1]{fontenc}
\RequirePackage{amsthm,amsmath}
\RequirePackage[numbers]{natbib}
\RequirePackage[colorlinks,citecolor=blue,urlcolor=blue]{hyperref}

\newcommand{\norm}[1]{\lVert#1\rVert}

\startlocaldefs
\numberwithin{equation}{section}
\theoremstyle{plain}

\endlocaldefs

\input{aospreamble.tex}

\newcommand{\bepsilon}{\bm{\delta}}
\renewcommand{\interior}{\textnormal{int}}
\newcommand{\BBB}{\mathbb{B}}

\usepackage{amsfonts,amsmath,amssymb,latexsym,bm,bbm}
\newcommand{\mbG}{\mathscr{G}}
\RequirePackage{amsthm,amsmath}
\usepackage{times}
\usepackage{natbib}

\newcommand{\deltaepsilon}{\delta(\epsilon)}

\renewcommand{\bta}{\bm{\eta}}

\usepackage{multirow}

\usepackage{mathtools}

\newcommand{\comp}{\mbox{comp}\,}
\renewcommand{\one}{\mathbbm{1}}
\renewcommand{\mM}{\mathbb{M}}

\newcommand{\etaspace}{\bm{\Xi}}
\newcommand{\fullspace}{\mathfrak{N}}
\renewcommand{\mT}{\mathbb{T}}
\usepackage[mathscr]{eucal}
\renewcommand{\bmu}{\bm{\mu}}
\renewcommand{\bA}{\bm{A}}
\renewcommand{\bt}{\bm{t}}

\renewcommand{\norm}[1]{\lVert#1\rVert}

\renewcommand{\mbR}{\mathbb{R}}

\renewcommand{\lip}{{\tiny\mbox{Lip}}}

\renewcommand{\lte}{&\leq&}
\renewcommand{\gte}{&\geq&}

\renewcommand{\btheta}{\bm{\theta}}
\renewcommand{\bTheta}{\bm{\Theta}}
\newcommand{\bThetas}{\bm{\Theta}^\star}
\renewcommand{\tends}{\to}

\renewcommand{\rint}{\mathop{\mbox{rint}}}

\renewcommand{\mS}{\mathbb{S}}
\renewcommand{\mbX}{\mathbb{X}}
\renewcommand{\mbY}{\mathbb{Y}}
\renewcommand{\uu}{x_i}
\renewcommand{\vv}{x_i^\star}

\renewcommand{\tv}{{\tiny\mbox{TV}}}

\renewcommand{\mB}{\mathscr{B}}

\usepackage{bbm,amsfonts,graphicx}

\newcommand{\ppp}{m}
\newcommand{\qqq}{q}
\newcommand{\diminf}{\ppp}
\newcommand{\rrr}{\ppp} 

\usepackage{multicol}

\renewcommand{\bx}{\bm{x}}
\renewcommand{\=}{&=&}

\renewcommand{\hide}[1]{}
\renewcommand{\ghost}[1]{}

\newcommand{\Thetad}{\Theta}
\newcommand{\bthetas}{\btheta^\star}
\renewcommand{\thetas}{\theta^\star}

\renewcommand{\bthetad}{{\btheta}_0}
\newcommand{\bThetad}{{\bTheta}_0}

\newcommand{\bthetah}{\widehat{\btheta}}

\renewcommand{\mX}{\mathbb{X}}

\renewcommand{\mY}{\mathbb{Y}}

\renewcommand{\be}{\begin{equation}\begin{array}{lllllllllllllllll}}
\renewcommand{\beno}{\begin{equation}\begin{array}{lllllllllllll}\nonumber}
\renewcommand{\ee}{\end{array}\end{equation}}
\renewcommand{\dis}{\displaystyle}
\renewcommand{\dsum}{\displaystyle\sum\limits}

\renewcommand{\dprod}{\displaystyle\prod\limits}

\renewcommand{\mR}{\mathbb{R}}

\renewcommand{\mbE}{\mathbb{E}}
\renewcommand{\mbP}{\mathbb{P}}
\renewcommand{\mA}{\mathscr{A}}
\newcommand{\mAs}{\mA}

\newcommand{\mJ}{\nabla_{\btheta}\, \bta(\btheta)}
\renewcommand{\mE}{\mathscr{E}(\delta(\epsilon))}
\renewcommand{\mbG}{\mathscr{G}(\omega)}

\newcommand{\longtitle}{Concentration and Consistency Results for Canonical and Curved Exponential-Family Models of Random Graphs}
\newcommand{\shorttitle}{Concentration and Consistency Results}

\newcommand{\ghoster}[0]{}

\begin{document}

\begin{frontmatter}

\title{\longtitle}
\runtitle{\shorttitle}

\ghoster{

\thankstext{T1}{\alert{Supported by NSF awards DMS-1513644 and DMS-1812119.}}

\begin{aug}
\author{\fnms{Michael} \snm{Schweinberger}\ead[label=e1]{m.s@rice.edu}}
\and
\author{\fnms{Jonathan} \snm{Stewart}\ead[label=e2]{jonathan.stewart@rice.edu}}
\affiliation{Rice University}
\address{
Michael Schweinberger, Jonathan Stewart\\
Department of Statistics\\
Rice University\\
6100 Main St, MS-138\\
Houston, TX 77005-1827\\
E-mail:\ m.s@rice.edu\\
jonathan.stewart@rice.edu
}
\end{aug}

\runauthor{Michael Schweinberger and Jonathan Stewart}

}

\begin{abstract}
Statistical inference for exponential-family models of random graphs with dependent edges is challenging.
We stress the importance of additional structure and show that additional structure facilitates statistical inference.
A simple example of a random graph with additional structure is a random graph with neighborhoods and local dependence within neighborhoods.
We develop the first concentration and consistency results for maximum likelihood and $M$-estimators of a wide range of canonical and curved exponential-family models of random graphs with local dependence.
All results are non-asymptotic and applicable to random graphs with finite populations of nodes,
although asymptotic consistency results can be obtained as well.
In addition,
we show that additional structure can facilitate subgraph-to-graph estimation,
and present concentration results for subgraph-to-graph estimators.
As an application,
we consider popular curved exponential-family models of random graphs,
with local dependence induced by transitivity and parameter vectors whose dimensions depend on the number of nodes.
\end{abstract}


\end{frontmatter}

\section{Introduction}
\label{sec:introduction}

Models of network data have witnessed a surge of interest in statistics and related areas \citep[e.g.,][]{Ko09a}.
Such data arise in the study of,
e.g., 
social networks,
epidemics,
insurgencies, 
and terrorist networks. 

Since the work of Holland and Leinhardt in the 1970s \citep[e.g.,][]{HpLs76},
it is known that network data exhibit a wide range of dependencies induced by transitivity and other interesting network phenomena \citep[e.g.,][]{ergm.book}.
Transitivity is a form of triadic closure in the sense that,
when a node $k$ is connected to two distinct nodes $i$ and $j$,
then $i$ and $j$ are likely to be connected as well,
which suggests that edges are dependent \citep[e.g.,][]{ergm.book}.
A large statistical framework for modeling dependencies among edges is given by discrete exponential-family models of random graphs,
called exponential-family random graphs \citep[e.g.,][]{FoSd86,WsPp96,Ha03p,SnPaRoHa04,HuGoHa08,ergm.book,Ha13,LaRiSa17}. 
Such models are popular among network scientists for the same reason Ising models are popular among physicists: 
Both classes of models enable scientists to model a wide range of dependencies of scientific interest \citep[e.g.,][]{ergm.book}.

Despite the appeal of the discrete exponential-family framework and its relationship to other discrete exponential-family models for dependent random variables \citep[e.g., Ising models and discrete Markov random fields,][]{Ch07a,RaWaLa10,BhMu18},
statistical inference for exponential-family random graphs is challenging.
One reason is that some exponential-family random graphs are ill-behaved \citep[e.g., the so-called triangle model,][]{Jo99,Ha03p,Sc09b,BaBrSl11,ChDi11,Bu10},
though well-behaved alternatives have been developed,
among them curved exponential-family random graphs \citep{SnPaRoHa04,HuGoHa08}. 
\alert{A second reason is that in most applications of exponential-family random graphs statistical inference is based on a single observation of a random graph with dependent edges.
Establishing desirable properties of estimators,
such as consistency,
is non-trivial when no more than one observation of a random graph with dependent edges is available.
While some consistency results have been obtained under independence assumptions \citep[][]{DiChSl11,XiNe11,RiPeFi13,ShRi11,Mu13,KrKo14,YaLeZh11,Yaetal18} and restrictive dependence assumptions \citep{XiNe11,ShRi11,Mu13}---as discussed in Section \ref{literature}---the existing consistency results do not cover the models most widely used in practice \citep[e.g.,][]{ergm.book}:
canonical and curved exponential-family random graphs with dependence among edges induced by transitivity and other interesting network phenomena \citep[][]{ergm.book}.}

We stress the importance of additional structure and show that additional structure facilitates statistical inference.
We consider here a simple and common form of additional structure,
called multilevel structure.
Network data with multilevel structure are popular in network science,
as the recent monograph of \citet{multilevelnetwork} and a growing number of applications demonstrate \citep[e.g.,][]{Wa13,ZaLo15,Lo16,slaughter2016multilevel,hollway2017multilevel,StScBoMo18}.
A simple form of multilevel structure is given by a partition of a population of nodes into subsets of nodes, called neighborhoods. 
In applications,
neighborhoods may correspond to 
school classes within schools,
departments within companies,
and units of armed forces.
It is worth noting that in multilevel networks the partition of the population of nodes is observed and models of multilevel networks attempt to capture dependencies within and between neighborhoods \citep[e.g.,][]{Wa13,ZaLo15,Lo16,slaughter2016multilevel,hollway2017multilevel,StScBoMo18},
whereas the well-known class of stochastic block models \citep{NkSt01} assumes that the partition is unobserved and that edges are independent conditional on the partition.

Additional structure in the form of multilevel structure offers opportunities in terms of statistical inference.
We take advantage of these opportunities to develop the first statistical theory which shows that statistical inference for many canonical and curved exponential-family random graphs with dependent edges is possible.
The main idea is based on a simple and general exponential-family argument that may be of independent interest.
It helps establish non-asymptotic probability statements about estimators of canonical and curved exponential families for dependent random variables under weak conditions,
as long as additional structure helps control the amount of dependence induced by the model and the sufficient statistics are sufficiently smooth functions of the random variables.
We exploit the main idea to develop the first concentration and consistency results for maximum likelihood and $M$-estimators of canonical and curved exponential-family random graphs with dependent edges,
under correct and incorrect model specifications.
All results are non-asymptotic and applicable to random graphs with finite populations of nodes,
although asymptotic consistency results can be obtained as well.
In addition,
we show that multilevel structure facilitates subgraph-to-graph estimation,
and present concentration results for subgraph-to-graph estimators.
As an application,
we consider popular curved exponential-family random graphs \citep{SnPaRoHa04,HuGoHa08},
with local dependence induced by transitivity and parameter vectors whose dimensions depend on the number of nodes.

These concentration and consistency results have important implications,
both in statistical theory and practice:
\bi
\item The most important implication is that statistical inference for transitivity and other network phenomena of great interest to network scientists is possible.
To date,
it has been widely believed that statistical inference for transitivity based on exponential-family random graphs is challenging \citep[e.g.,][]{ShRi11,ChDi11},
but additional structure in the form of multilevel structure facilitates it.
\item Network scientists can benefit from collecting network data with multilevel structure,
because multilevel structure can facilitate statistical inference for exponential-family random graphs with dependent edges.
\ei

Last,
but not least,
it is worth noting that these concentration and consistency results cover two broad inference scenarios:
\bi
\item {\em Inference scenarios with finite populations of nodes.}
In many applications of exponential-family random graphs,
there is a finite population of nodes and a population graph is assumed to have been generated by an exponential-family random graph model.
A common goal of statistical inference,
then,
is to estimate the parameters of the data-generating exponential-family random graph model based on a complete or incomplete observation of the population graph.
Our concentration results cover inference scenarios with finite populations of nodes,
when the whole population graph is observed or when neighborhoods are sampled and the subgraphs induced by the sampled neighborhoods are observed.
\item {\em Inference scenarios with populations of nodes growing without bound.}
In addition,
our concentration results can be used to obtain asymptotic consistency results by allowing the number of neighborhoods to grow without bound.
The resulting asymptotic consistency results resemble asymptotic consistency results in other areas of statistics,
albeit with two notable differences:
first,
the units of statistical analysis are subsets of nodes (neighborhoods) rather than nodes or edges;
and,
second,
the sizes of the units need not be identical,
but are similar in a well-defined sense.
\ei
Since the first application is more interesting than the second one,
we state all results with the first application in mind,
i.e.,
all results focus on random graphs with finite populations of nodes,
although we do mention some asymptotic consistency results along the way.

The remainder of our paper is structured as follows.
Section \ref{sec:model} introduces models.
Section \ref{mle} describes concentration and consistency results for maximum likelihood and $M$-estimators,
under correct and incorrect model specifications.
Section \ref{extendability} shows that multilevel structure facilitates subgraph-to-graph estimation.
A comparison with existing consistency results can be found in Section \ref{literature}.
Section \ref{sec:simulations} presents simulation results.

\section{Exponential-family random graphs with multilevel structure}
\label{sec:model}

We introduce exponential-family random graphs with multilevel structure.

A simple and common form of multilevel structure is a partition of a population of nodes into $K \geq 2$ non-empty subsets of nodes $\mA_1, \dots, \mA_K$,
called neighborhoods.
We note that in multilevel networks the partition of the population of nodes is observed \citep[e.g.,][]{Wa13,ZaLo15,Lo16,slaughter2016multilevel,hollway2017multilevel,StScBoMo18} and that some neighborhoods may be larger than others.
We consider random graphs with undirected edges that may be either absent or present or may have weights,
where the weights are elements of a countable set.
Extensions to random graphs with directed edges are straightforward.
Let $\bX = (\bX_k)_{k=1}^K$ and $\bY = (\bY_{k,l})_{k<l}^K$ be sequences of within- and between-neighborhood edge variables based on a sequence of neighborhoods $\mA_1, \dots, \mA_K$,
where $\bX_k = (X_{i,j})_{i\in\mA_k\, <\, j\in\mA_k}$ and $\bY_{k,l} = (Y_{i,j})_{i\in\mA_k,\, j\in\mA_l}$ ($k < l$) correspond to within- and between-neighborhood edge variables $X_{i,j} \in \mbX_{i,j}$ and $Y_{i,j} \in \mbY_{i,j}$,
taking on values in countable sets $\mbX_{i,j}$ and $\mbY_{i,j}$,
respectively.
We exclude self-edges,
assuming that $X_{i,i} = 0$ holds with probability $1$ ($i \in \mA_k$,\; $k = 1, \dots, K$),
and write 
$\mX_k = \prod_{i\in\mA_k\, <\, j\in\mA_k} \mbX_{i,j}$,\,
$\mX = \prod_{k=1}^K \prod_{i\in\mA_k\, <\, j\in\mA_k} \mbX_{i,j}$, 
and $\mY = \prod_{k<l}^K \prod_{i\in\mA_k,\, j\in\mA_l} \mbY_{i,j}$.

We assume that within-neighborhood edges $\bX$ are independent of between-neighborhood edges $\bY$,
i.e.,
\beno
\mbP(\bX \in \mathscr{X},\, \bY \in \mathscr{Y})
\= \mbP(\bX \in \mathscr{X})\, \mbP(\bY \in \mathscr{Y})
& \mbox{for all} & \mathscr{X}\times\mathscr{Y}\, \subseteq\, \mbX\times\mbY,
\ee
where $\mbP$ denotes a probability distribution with support $\mbX\times\mbY$.
We do not assume that edges are independent,
neither within nor between neighborhoods.

While in principle both within-neighborhood edge variables $\bX$ and between-neighborhood edge variables $\bY$ may be of interest,
we focus on within-\linebreak
neighborhood edge variables, 
which are of primary interest in applications \citep[e.g.,][]{multilevelnetwork,Wa13,ZaLo15,Lo16,slaughter2016multilevel,hollway2017multilevel,StScBoMo18}.
We therefore restrict attention to the probability law of $\bX$ and do not make assumptions about the probability law of $\bY$.
We assume that the parameter vectors of the probability laws of $\bX$ and $\bY$ are variation-independent, 
i.e.,
the parameter space is a product space,
so that statistical inference concerning the parameter vector of the probability law of $\bX$ can be based on $\bX$ without requiring knowledge of $\bY$.

The distribution of within-neighborhood edge variables $\bX$ is presumed to belong to an exponential family with local dependence,
defined as follows.

\s

{\em Definition. Exponential family with local dependence.}
An exponential family with local dependence is an exponential family of distributions with countable support $\mX$,
having densities with respect to counting measure of the form
\be
\label{local.ergm}
p_{\bta}(\bx)
\= \exp\left(\langle\bta,\, s(\bx)\rangle - \psi(\bta)\right)\, \nu(\bx)\s
\\
\= \exp\left(\dsum_{k=1}^K \langle\bta_k,\, s_k(\bx_k)\rangle - \psi(\bta)\right)\, \nu(\bx),
\ee
where
\beno
\psi(\bta)
\= \log\dsum_{\bx_1 \in \mX_1} \cdots \dsum_{\bx_K \in \mX_K} \exp\left(\dsum_{k=1}^K \langle\bta_k,\, s_k(\bx_k)\rangle\right)\, \nu(\bx)
\ee
and $\nu(\bx) = \prod_{k=1}^K \nu_k(\bx_k)$.

\s

In other words,
edges may depend on other edges in the same neighborhood,
but do not depend on edges in other neighborhoods \citep[][]{ScHa13}.
Here,
$\langle\bta,\, s(\bx)\rangle$ $=$\linebreak 
$\sum_{k=1}^K \langle\bta_k,\, s_k(\bx_k)\rangle$ is the inner product of a natural parameter vector $\bta \in \mR^{\ppp}$ and a sufficient statistic vector $s: \mX \mapsto \mR^{\ppp}$ while $\bta_k \in \mR^{\ppp_k}$ and $s_k: \mX_k \mapsto \mR^{\ppp_k}$ denote the natural parameter vector and sufficient statistic vector of neighborhood $\mA_k$,
respectively ($k = 1, \dots, K$).
The functions $\nu: \mbX \mapsto \mR^+ \cup \{0\}$ and $\nu_k: \mbX_k \mapsto \mR^+ \cup \{0\}$ ($k = 1, \dots, K$) along with the sample space $\mbX$ determine the reference measure of the exponential family.
A careful discussion of the choice of reference measure can be found in \citet{Kr11}.

We consider here a wide range of exponential families with local dependence.
A specific example of an exponential family with local dependence can be found in Section \ref{sec:curved}.
In the following,
we introduce selected exponential-family terms in order to distinguish exponential families from subfamilies of exponential families.
Subfamilies of exponential families give rise to distinct theoretical challenges,
and thus require a separate treatment.
We therefore introduce the classic exponential-family notions of full and non-full exponential families,
canonical and curved exponential families,
and minimal exponential families.
These exponential-family terms are taken from the monographs on exponential families by \citet{BN78} and \citet{Br86} and have been used in other recent works as well: see, e.g., \citet{LaRiSa17}, \citet{fienberg-2008}, and \citet{Ge09}.
To help ensure that parameters are identifiable,
we assume that exponential families of the form \eqref{local.ergm} are minimal in the sense of \citet{BN78} and \citet{Br86},
i.e.,
the closure of the convex hull of the set $\{s(\bx): \nu(\bx) > 0\}$ is not contained in a proper affine subspace of $\mR^{\ppp}$ \citep[e.g.,][p.\ 2]{Br86}.
It is well-known that all non-minimal exponential families can be reduced to minimal exponential families \citep[e.g.,][Theorem 1.9, p.\ 13]{Br86}.
We consider both full and non-full exponential families of the form \eqref{local.ergm}.
An exponential family $\{\mbP_{\bta},\, \bta \in \etaspace\}$ is full if $\etaspace = \fullspace$ and non-full if $\etaspace \subset \fullspace$,
where $\fullspace = \{\bta \in \mR^{\ppp}: \psi(\bta) < \infty\}$ is the natural parameter space,
i.e.,
the largest set of possible values the natural parameter vector $\bta$ can take on.
While full exponential families may be more convenient on mathematical grounds,
non-full exponential families---the most important example being curved exponential families \citep[e.g.,][]{Ef78,KaVo97}---offer parsimonious parameterizations of exponential families where the dimension $\ppp_k$ of neighborhood-dependent natural parameter vectors $\bta_k$ is an increasing function of the number of nodes in neighborhoods $\mA_k$ ($k = 1, \dots, K$),
and have turned out to be useful in practice \citep{SnPaRoHa04,HuGoHa08}.
A simple approach to generating non-full exponential families is to assume that $\bta: \interior(\bTheta) \mapsto \interior(\fullspace)$ is a known function of a parameter vector $\btheta \in \bTheta$,
where $\bTheta \subseteq \{\btheta \in \mR^{\qqq}: \psi(\bta(\btheta)) < \infty\}$,
$\interior(\bTheta)$ and $\interior(\fullspace)$ denote the interiors of $\bTheta$ and $\fullspace$,
respectively,
and $\qqq \leq \ppp$.
It is convenient to distinguish exponential families that can be reduced to canonical exponential families with natural parameter vectors of the form $\bta(\btheta) = \btheta$ from those that cannot.
An exponential family can be reduced to a canonical exponential family with $\bta(\btheta) = \btheta$ when the map $\bta: \interior(\bTheta) \mapsto \interior(\fullspace)$ is affine.
In other words,
if $\bta(\btheta) = \bA\, \btheta + \bm{b}$ with $\bA \in \mR^{\ppp\times\qqq}$ and $\bm{b} \in \mR^{\ppp}$,
then the exponential family can be reduced to a canonical exponential family with $\bta(\btheta) = \btheta$ by absorbing $\bA$ into the sufficient statistic vector and $\bm{b}$ into the reference measure.
We therefore call all exponential families with affine maps $\bta: \interior(\bTheta) \mapsto \interior(\fullspace)$ canonical exponential families,
and call all exponential families with non-affine maps $\bta: \interior(\bTheta) \mapsto \interior(\fullspace)$ curved exponential families.
We note that our definition of a curved exponential family is broader than the one used in \citet{Ef75,Ef78}, \citet[][pp.\ 81--84]{Br86}, and \citet{KaVo97}.
The main reason is that we do not restrict the map $\bta: \interior(\bTheta) \mapsto \interior(\fullspace)$ to be differentiable,
because our concentration and consistency results in Sections \ref{mle} and \ref{extendability} do not require differentiability.

Throughout,
we assume that the neighborhoods are of the same order of magnitude and that the neighborhood-dependent natural parameters $\eta_{k,i}(\btheta)$ are of the form $\eta_{k,i}(\btheta) = \eta_i(\btheta)$ ($i = 1, \dots, \ppp_k$, $k = 1, \dots, K$).
We define neighborhoods of the same order of magnitude as follows.

\s

\begin{samepage}
{\em Definition. Neighborhoods of the same order of magnitude.}
A sequence of neighborhoods $\mA_1, \dots, \mA_K$ is of the same order of magnitude if there exists a universal constant $A > 1$ such that $\max_{1 \leq k \leq K} |\mA_k| \leq A\, \min_{1 \leq k \leq K} |\mA_k|$ ($K = 1, 2, \dots$).\s

In other words,
the largest neighborhood size is a constant multiple of the smallest neighborhood size,
so that the sizes of neighborhoods may not be identical,
but are similar in a well-defined sense.
The definition is satisfied when the sizes of neighborhoods are bounded above.
When the number of neighborhoods $K$ grows and the sizes of neighborhoods grow with $K$,
then the definition implies that the sizes of neighborhoods grow at the same rate.
\end{samepage}
We note that when the neighborhoods are not of the same order of magnitude,
the natural parameters of neighborhoods may have to depend on the order of magnitude of neighborhoods \citep[e.g.,][]{KrHaMo11,KrKo14,BuAl15},
because there are good reasons to believe that small and large within-neighborhood subgraphs are not governed by the same natural parameters \citep{ShRi11,CrDe15,LaRiSa17}.
Size-dependent parameterizations have an important place in the exponential-family random graph framework,
and some promising size-dependent parameterizations have been proposed.
Most of them assume that natural parameters consist of size-invariant parameters and size-dependent deviations.
The size-dependent deviations may be size-dependent offsets \citep[e.g.,][]{KrHaMo11,KrKo14,BuAl15} or functions of size-dependent covariates \citep{slaughter2016multilevel}.
Some of those size-dependent deviations can be absorbed into the sufficient statistic vector and reference measure,
and are hence covered by our main concentration and consistency results in Sections \ref{mle} and \ref{extendability}.
However,
the topic of size-dependent parameterizations is an important topic in its own right,
and deserves a separate treatment that is beyond the scope of our paper.

The assumption $\eta_{k,i}(\btheta) = \eta_i(\btheta)$ ($i = 1, \dots, \ppp_k$, $k = 1, \dots, K$) implies that the exponential families considered here can be reduced to exponential families with natural parameter vectors of the form
\beno
\bta(\btheta)
\= \left(\eta_1(\btheta), \dots, \eta_{\diminf}(\btheta)\right)
\ee
and sufficient statistic vectors of the form
\beno
s(\bx)
\= \left(s_1(\bx), \dots, s_{\diminf}(\bx)\right),
\ee
where $s_i(\bx) = \sum_{k=1}^K s_{k,i}(\bx_k)$ ($i = 1, \dots, \diminf$) and $\diminf = \max_{1 \leq k \leq K} \ppp_k$.
We assume that the dimensions $\ppp_k$ of the neighborhood-dependent natural parameter vectors $\bta_k(\btheta)$ are non-decreasing functions of the sizes $|\mA_k|$ of neighborhoods $\mA_k$ ($k = 1, \dots, K$),
which implies that $\diminf = \max_{1 \leq k \leq K} \ppp_k$ is a non-decreasing function of $\norm{\mAs}_\infty = \max_{1 \leq k \leq K} |\mAs_k|$.
The dimensions $\ppp_k$ ($k = 1, \dots, K$) and $\diminf$\, do not depend on other quantities.
The dimension $\qqq$ of parameter vector $\btheta$ satisfies $\qqq \leq \ppp$,
as mentioned above.

\s

{\em Notation.}
To prepare the ground for the concentration and consistency results in Sections \ref{mle} and \ref{extendability},
we introduce mean-value parameterizations of exponential families along with additional notation.
Mean-value parameterizations facilitate concentration and consistency results, 
because concentration inequalities \citep{BoLuMa13} bound probabilities of deviations from means and the mean-value parameters of an exponential family are the means of the sufficient statistics,
defined by $\bmu(\bta(\btheta)) = \mbE_{\bta(\btheta)}\, s(\bX) \in \rint(\mM)$ \citep[][p.\ 2 and p.\ 75]{Br86}.
\alert{Here, 
$\mbE_{\bta(\btheta)}\, s(\bX)$ is the expectation of $s(\bX)$ with respect to exponential-family distributions $\mbP_{\bta(\btheta)}$ having densities of the form \eqref{local.ergm},
$\mM$ is the convex hull of the set $\{s(\bx): \nu(\bx) > 0\}$,
and $\rint(\mM)$ is the relative interior of $\mM$.
We denote the data-generating parameter vector by $\bthetas \in \interior(\bTheta)$ and write $\mbP \equiv \mbP_{\bta(\bthetas)}$ and $\mbE \equiv \mbE_{\bta(\bthetas)}$.
An open ball in $\mR^{v}$ ($v \geq 1$) with center $\bm{c} \in \mR^{v}$ and radius $\rho > 0$ is denoted by $\mB(\bm{c},\, \rho)$.}
We write $\norm{.}_1$, $\norm{.}_2$, and $\norm{.}_\infty$ to refer to the $\ell_1$-, $\ell_2$-, and $\ell_\infty$-norm of vectors,
respectively.
Throughout,
uppercase letters $A, B, C, C_1, C_2, \dots > 0$ denote universal constants,
which are independent of all other quantities of interest and may be recycled from line to line.
The function $d: \mX \times \mX \mapsto \{0, 1, \dots\}$ denotes the Hamming metric,
defined by
\beno
d(\bx_1, \bx_2)
\= \dsum_{k=1}^K\; \dsum_{i\in\mA_k\, <\, j\in\mA_k} \one(x_{1,i,j} \neq x_{2,i,j}), 
& (\bx_1, \bx_2) \in \mX \times \mX,
\ee
where $\one(x_{1,i,j} \neq x_{2,i,j})$ is an indicator function,
which is $1$ if $x_{1,i,j} \neq x_{2,i,j}$ and is $0$ otherwise.

\section{Concentration and consistency results: maximum likelihood and $M$-estimators}
\label{mle}

In many applications of exponential-family random graphs,
the parameter vector of primary interest is $\btheta$.
To estimate the parameter vector $\btheta$ of a wide range of full and non-full, curved exponential families under weak assumptions on the map $\bta: \interior(\bTheta) \mapsto \interior(\fullspace)$,
we consider an estimating function \citep{Go91}---a function $g: \bTheta \times \mX \mapsto \mR$ of both $\btheta$ and $\bX$---of the form
\be
\label{m.estimating}
g(\btheta;\; \widehat{\bmu(\bta(\bthetas))})
\= \norm{\widehat{\bmu(\bta(\bthetas))} - \bmu(\bta(\btheta))}_2,
& \btheta \in \bTheta,
\ee
which is an approximation of
\beno
g(\btheta;\; \bmu(\bta(\bthetas)))
\= \norm{\bmu(\bta(\bthetas)) - \bmu(\bta(\btheta))}_2,
& \btheta \in \bTheta,
\ee
where $\widehat{\bmu(\bta(\bthetas))} = s(\bX)$ is an estimator of the data-generating mean-value parameter vector $\bmu(\bta(\bthetas)) = \mbE_{\bta(\bthetas)}\, s(\bX) \in \rint(\mM)$.
The fact that $g(\btheta;\, \widehat{\bmu(\bta(\bthetas))})$ is an approximation of $g(\btheta;\, \bmu(\bta(\bthetas)))$ follows from the triangle inequality,
which shows that
\beno
\label{approx}
\big|g(\btheta;\; \widehat{\bmu(\bta(\bthetas))}) - g(\btheta;\; \bmu(\bta(\bthetas)))\big|
\lte \norm{\widehat{\bmu(\bta(\bthetas))} - \bmu(\bta(\bthetas))}_2,
& \btheta \in \bTheta,
\ee
along with the fact that,
under suitable conditions, 
$\norm{\widehat{\bmu(\bta(\bthetas))} - \bmu(\bta(\bthetas))}_2$ is small with high probability,
as we will show in Proposition \ref{existence.eta.mle} in Section \ref{sec:probability}.

Estimating function \eqref{m.estimating} has at least three advantages.
First,
estimating function \eqref{m.estimating} admits concentration and consistency results under weak assumptions on the map $\bta: \interior(\bTheta) \mapsto \interior(\fullspace)$.
Indeed,
the map $\bta: \interior(\bTheta) \mapsto \interior(\fullspace)$ satisfies the main assumptions of Section \ref{mle} as long as the map is one-to-one and continuous,
but it need not be differentiable.
The weakness of these assumptions implies that the main results of Section \ref{mle} cover a vast range of full and non-full exponential families---including, but not limited to curved exponential families---and it is possible to verify these assumptions in some of the most challenging applications,
as demonstrated in Section \ref{sec:curved}.
Second,
estimating function \eqref{m.estimating} is natural,
because maximum likelihood estimation of the data-generating natural parameter vector $\btas$ of an exponential family with natural parameter vector $\bta$ and sufficient statistic vector $s(\bx)$ can be based on the gradient of the loglikelihood function $\nabla_{\bta}\, \log p_{\bta}(\bx) = \widehat{\bmu(\btas)} - \bmu(\bta)$ provided $\bta \in \interior(\fullspace)$,
where $\widehat{\bmu(\btas)} = s(\bx)$ \citep[][Lemma 5.3, p.\ 146]{Br86}.
Therefore,
maximum likelihood estimation of $\btas$ can be based on estimating functions of the form $\norm{\widehat{\bmu(\btas)} - \bmu(\bta)}_2$.
By the parameterization-invariance of maximum likelihood estimators,
maximum likelihood estimation of functions of $\btas$,
such as $\bthetas$,
can be based on $\norm{\widehat{\bmu(\bta(\bthetas))} - \bmu(\bta(\btheta))}_2$ provided the map $\bta: \interior(\bTheta) \mapsto \interior(\fullspace)$ is one-to-one.
We note that estimating function \eqref{m.estimating} is chosen for mathematical convenience,
facilitating concentration and consistency results for maximum likelihood estimators of many full and non-full, curved exponential families under weak assumptions on the map $\bta: \interior(\bTheta) \mapsto \interior(\fullspace)$,
and is not chosen for computational convenience.
Last,
but not least,
the simple form of estimating function \eqref{m.estimating} helps determine when minimizers of \eqref{m.estimating} exist and are unique, and how the minimizers are related to each other when there is more than one minimizer.
These advantages are most useful in non-full exponential families,
in particular curved exponential families.

In the following,
we assume that the map $\bta: \interior(\bTheta) \mapsto \interior(\fullspace)$ is one-to-one.
A natural class of estimators is hence given by
\beno
\label{maxlik}
\widehat\btheta
\= \left\{\btheta \in \bTheta:\;\; g(\btheta;\; \widehat{\bmu(\bta(\bthetas))})\;\; =\;\; \inf\limits_{\dot\btheta\in\bTheta}\; g(\dot\btheta;\; \widehat{\bmu(\bta(\bthetas))})\right\}.
\ee
If the set $\bthetah$ is non-empty,
it may contain one element (e.g., in full exponential families) or more than one element (e.g., in non-full exponential families).
If the set $\bthetah$ contains more than one element,
then all elements of the set $\bthetah$ map to mean-value parameter vectors $\bmu(\bta(\bthetah))$ that have the same $\ell_2$-distance from $\widehat{\bmu(\bta(\bthetas))}$ by construction of estimating function \eqref{m.estimating}.
In addition,
if the set $\bthetah$ is non-empty,
then all minimizers $\bthetah$ of $g(\btheta;\; \widehat{\bmu(\bta(\bthetas))})$ are approximations of the minimizer of $g(\btheta;\; \bmu(\bta(\bthetas)))$.
The minimizer of $g(\btheta;\; \bmu(\bta(\bthetas)))$ is unique and is given by the data-generating parameter vector $\bthetas$ provided $\bmu(\bta(\bthetas))) \in \rint(\mM)$:
\beno
\label{maxlik}
\bthetas
\= \left\{\btheta \in \bTheta:\;\; g(\btheta;\; \bmu(\bta(\bthetas)))\;\; =\;\; \inf\limits_{\dot\btheta\in\bTheta}\; g(\dot\btheta;\; \bmu(\bta(\bthetas)))\right\}.
\ee
The data-generating parameter vector $\bthetas$ is the unique minimizer of $g(\btheta;\; \bmu(\bta(\bthetas)))$ 
$= \norm{\bmu(\bta(\bthetas)) - \bmu(\bta(\btheta))}_2$,
because $\bthetas \in \interior(\bTheta)$ and the map $\bta: \interior(\bTheta) \mapsto \interior(\fullspace)$ is one-to-one by assumption,
while the map $\bmu: \interior(\fullspace) \mapsto \rint(\mM)$ is one-to-one by classic exponential-family theory \citep[][Theorem 3.6, p.\ 74]{Br86}.
Therefore,
$\norm{\bmu(\bta(\bthetas)) - \bmu(\bta(\btheta))}_2 = 0$ holds if and only if $\btheta=\bthetas$.

In the remainder of the section,
we show that the estimator $\bthetah$ is close to the data-generating parameter vector $\bthetas$ with high probability under weak conditions.
We first sketch the main idea in Section \ref{sec:idea} and then discuss concentration and consistency results for maximum likelihood and $M$-estimators in Sections \ref{sec:curved} and \ref{sec:model.estimation},
respectively.
An application to popular curved exponential-family random graphs is presented in Section \ref{sec:curved}.

\subsection{Main idea: A non-asymptotic approach to full and non-full, curved exponential families for dependent random variables}
\label{sec:idea}

Establishing concentration and consistency results for estimators of full and non-full, curved exponential-family random graphs with dependent edges is non-trivial,
for at least three reasons.
First,
the map $\bta: \interior(\bTheta) \mapsto \interior(\fullspace)$ may not be affine and may not be differentiable.
Second,
in many non-full exponential families,
the mean-value parameter vector $\bmu(\bta(\btheta))$ is not available in closed form and there is no simple and known relationship between the mean-value parameter vector $\bmu(\bta(\bthetah))$ evaluated at $\bthetah$ and the sufficient statistic vector $s(\bX)$.
Third,
concentration results for functions of edges,
such as $s(\bX)$,
are more challenging when edges are dependent rather than independent.
As a result,
studying the behavior of the estimating function $\norm{s(\bX) - \bmu(\bta(\bthetah))}_2$ and its minimizer $\bthetah$ is not straightforward.

\begin{figure}[t]
\vspace{-.2cm}
\begin{center}
\includegraphics[scale=.65]{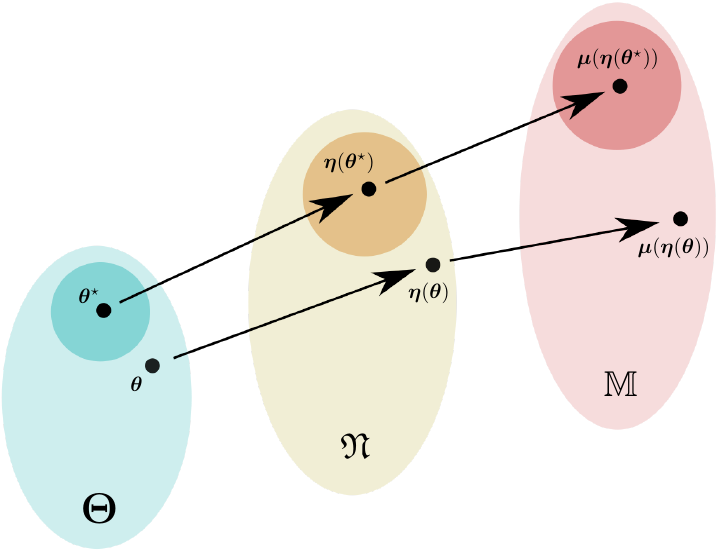}
\end{center}
\vspace{-.33cm}
\caption{\label{picture}The figure demonstrates the main assumption:
for all $\epsilon > 0$ small enough so that $\mathscr{B}(\bthetas, \epsilon) \subseteq \interior(\bTheta)$,
there exist $\gamma > 0$ and $\delta > 0$ such that $\btheta \not\in\, \mB(\bthetas, \epsilon)$ implies $\bta(\btheta) \not\in\, \mB(\bta(\bthetas),\, \gamma)$ and $\bmu(\bta(\btheta)) \not\in\, \mB(\bmu(\bta(\bthetas)), \delta)$.
}
\end{figure}

Our main idea is based on a simple and general exponential-family argument that may be of independent interest.
It helps establish non-asymptotic probability statements about estimators of full and non-full, curved exponential families for dependent random variables under weak conditions.

We make a single weak assumption:
for all $\epsilon > 0$ small enough so that $\mathscr{B}(\bthetas, \epsilon) \subseteq \interior(\bTheta)$,
there exist $\gamma > 0$ and $\delta > 0$ such that 
\beno
\btheta \not\in \mB(\bthetas, \epsilon)
\;\implies\; \bta(\btheta) \not\in \mB(\bta(\bthetas),\, \gamma)
\;\implies\; \bmu(\bta(\btheta)) \not\in \mB(\bmu(\bta(\bthetas)),\, \delta).
\ee
A graphical representation of the main assumption can be seen in Figure \ref{picture}.
A formal statement of the assumption can be found in Theorem \ref{nonfull} in Section \ref{sec:curved},
where $\gamma$ and $\delta$ depend on $\epsilon$ and $\delta$ depends on the sizes of neighborhoods $\mA_1, \dots, \mA_K$.
The main assumption is satisfied when the map $\bta: \interior(\bTheta) \mapsto \interior(\fullspace)$ is one-to-one and continuous,
but it need not be differentiable.
Note that the map $\bmu: \interior(\fullspace) \mapsto \rint(\mM)$ is one-to-one and continuous by classic exponential-family theory \citep[][Theorem 3.6, p.\ 74]{Br86}. 

The main assumption has a simple, but important implication.
As no element of $\fullspace$ outside of the ball $\mB(\bta(\bthetas),\, \gamma)$ maps to an element of the ball $\mB(\bmu(\bta(\bthetas)),\, \delta)$,
any element $\bmu(\bta(\btheta))$ of $\mB(\bmu(\bta(\bthetas)),\, \delta)$ must correspond to an element $\bta(\btheta)$ of $\mB(\bta(\bthetas),\, \gamma)$,
which in turn must correspond to an element $\btheta$ of $\mB(\bthetas, \epsilon)$,
so that
\beno
\bmu(\bta(\btheta)) \in \mB(\bmu(\bta(\bthetas)),\, \delta)
\;\implies\; \bta(\btheta) \in \mB(\bta(\bthetas),\, \gamma)
\;\implies\; \btheta \in \mB(\bthetas, \epsilon).
\ee
As a result, 
the probability of event $\bthetah \in \mB(\bthetas, \epsilon)$ can be bounded by bounding the probability of event $\bmu(\bta(\widehat\btheta))\; \in\; \mB(\bmu(\bta(\bthetas)),\, \delta)$.

A challenge,
which complicates probability statements about the event $\bmu(\bta(\widehat\btheta))$ $\in$ $\mB(\bmu(\bta(\bthetas)),\, \delta)$,
is that in many non-full exponential families $\bmu(\bta(\widehat\btheta))$ is not available in closed form and there is no simple and known relationship between $\bmu(\bta(\widehat\btheta))$ and the sufficient statistic vector $s(\bX)$.
To appreciate the difficulty of the problem,
suppose that $s(\bx) \in \rint(\mM)$ is observed,
so that $\widehat{\bmu(\bta(\bthetas))} = s(\bx) \in \rint(\mM)$.
The subset $\mM(\bTheta)$ of $\rint(\mM)$ induced by $\bTheta$ is defined by
\beno
\mM(\bTheta)
&=& \left\{\bmu' \in \rint(\mM): \mbox{ there exists } \btheta \in \bTheta \mbox{ such that } \bmu(\bta(\btheta)) = \bmu'\right\}.
\ee
In full exponential families,
$\mM(\bTheta) = \rint(\mM)$.
Thus,
there exists a minimizer $\bthetah \in \interior(\bTheta)$ of the estimating function $\norm{s(\bx) - \bmu(\bta(\btheta))}_2$ satisfying $\bmu(\bta(\bthetah)) = s(\bx)$.
In fact,
the minimizer is unique, 
because the maps $\bta: \interior(\bTheta) \mapsto \interior(\fullspace)$ and $\bmu: \interior(\fullspace) \mapsto \rint(\mM)$ are one-to-one \citep[][Theorem 3.6, p.\ 74]{Br86}.
As a result,
there is a simple and known relationship between $\bmu(\bta(\bthetah))$ and $s(\bx)$.
By contrast,
in non-full exponential families, 
$\mM(\bTheta)$ is a proper subset of $\rint(\mM)$,
because non-full exponential families are subfamilies of exponential families that exclude some natural parameter vectors along with the corresponding mean-value parameter vectors.
The problem, 
more often than not, 
is that it is unknown which mean-value parameter vectors are excluded,
because the mean-value parameter vectors are not available in closed form.
As a consequence,
there is no simple and known relationship between $\bmu(\bta(\widehat\btheta))$ and $s(\bx)$,
and it is not straightforward to determine where $\bmu(\bta(\widehat\btheta))$ is located in $\rint(\mM)$,
provided $\bmu(\bta(\widehat\btheta))$ is non-empty.
So bounding the probability of event $\bmu(\bta(\widehat\btheta))$ $\in$ $\mB(\bmu(\bta(\bthetas)),\, \delta)$ in non-full exponential families,
in particular curved exponential families,
is non-trivial.

But not all is lost.
Despite the challenge of characterizing $\mM(\bTheta) \subset \rint(\mM)$ and hence $\bmu(\bta(\bthetah)) \subset \mM(\bTheta)$,
it can be shown that,
under suitable conditions,
the set $\bmu(\bta(\bthetah)) \subset \mM(\bTheta)$ is non-empty and all elements of $\bmu(\bta(\bthetah))$ are close to $\bmu(\bta(\bthetas))$ with high probability.
Indeed,
when the set $\bthetah$ is non-empty,
each element of the set $\bthetah$ satisfies
\beno
\norm{\bmu(\bta(\bthetah)) - \bmu(\bta(\bthetas))}_2
\lte \norm{s(\bx) - \bmu(\bta(\bthetah))}_2 + \norm{s(\bx) - \bmu(\bta(\bthetas))}_2.
\ee
While characterizing $\mM(\bTheta)$ is difficult,
we do know one fact about $\mM(\bTheta)$:\linebreak
$\mM(\bTheta)$ contains the data-generating mean-value parameter vector $\bmu(\bta(\bthetas))$,
which implies that the $\ell_2$-distance of $\bmu(\bta(\bthetah))$ from $s(\bx)$ cannot exceed the $\ell_2$-distance of $s(\bx)$ from $\bmu(\bta(\bthetas)) \in \mM(\bTheta)$.
Thus,
we obtain the upper bound
\be
\label{triangle10}
\norm{\bmu(\bta(\bthetah)) - \bmu(\bta(\bthetas))}_2
\lte 2\, \norm{s(\bx) - \bmu(\bta(\bthetas))}_2.
\ee
If the right-hand side of \eqref{triangle10} can be shown to be small with high probability,
then the problem of bounding the probability of event $\bthetah \in \mB(\bthetas, \epsilon)$ can be converted into the problem of bounding the probability of the event that $s(\bX)$ is close to $\bmu(\bta(\bthetas)) = \mbE_{\bta(\bthetas)}\, s(\bX) \in \rint(\mM)$ in a well-defined sense.
All we need to bound the probabilities of those events are concentration results for the sufficient statistic vector $s(\bX)$,
which can be established as long as
(a) the dependence among edges is sufficiently weak;
and
(b) the sufficient statistics are sufficiently smooth functions of edges.
Concentration results are facilitated by additional structure that helps control the amount of dependence induced by the model and the smoothness of sufficient statistics.
We focus here on a simple form of additional structure in the form of multilevel structure,
which controls the dependence among edges by constraining it to neighborhoods.
But there are many other forms of additional structure that could help address (a) and (b),
e.g.,
other forms of multilevel structure or spatial structure.

The most important implication, then, is that statistical inference for many\linebreak 
exponential-family random graphs is possible and can be justified by statistical theory,
provided a suitable form of additional structure is available.
Indeed,
our main idea helps establish concentration and consistency results for estimators of
\bi
\item many full and non-full, curved exponential families;
\item many models with dependent edges;
\item finite populations of nodes;
\ei
as long as there is additional structure that helps address (a) and (b).

It is worth noting that verifying the main assumption is by no means trivial.
Its verification is easiest when the sufficient statistics are monotone functions of edges,
i.e.,
functions of a graph that do not decrease (increase) when edges are added to the graph.
It is less straightforward when the sufficient statistics are not monotone functions of edges,
as is the case with many popular curved exponential-family random graphs with geometrically weighted model terms \citep{SnPaRoHa04,HuGoHa08}.
But we demonstrate in Section \ref{sec:curved} that the main assumption can be verified even when the sufficient statistics are not monotone functions of a graph,
using curved exponential-family random graphs with geometrically weighted model terms as an example.

We make these ideas rigorous in Theorem \ref{nonfull} in Section \ref{sec:curved}.
An application to curved exponential-family random graphs is presented in Corollary \ref{corollary.curved} in Section \ref{sec:curved}.
More general results for $M$-estimators,
under correct and incorrect model specifications,
are mentioned in Section \ref{sec:model.estimation}.
To prepare the ground for these results,
we first introduce concentration results for sufficient statistics in Section \ref{sec:probability}.

\subsection{Concentration results for sufficient statistics}
\label{sec:probability}

To obtain concentration results for sufficient statistics,
we need concentration results for functions of random graphs with dependent edges.

Such concentration results are non-trivial for at least two reasons.
First,
exponential families of the form \eqref{local.ergm} may induce strong dependence within neighborhoods and the sizes of neighborhoods need not be identical.
Second,
exponential families of the form \eqref{local.ergm} can induce a wide range of dependencies within neighborhoods.
Therefore,
we need general-purpose concentration inequalities that cover a wide range of dependencies.

The following general-purpose concentration inequality addresses the challenges discussed above.
It shows that the dependence induced by exponential families of the form \eqref{local.ergm} may be strong within neighborhoods but is sufficiently weak overall to obtain concentration results as long as the neighborhoods are not too large.

\begin{samepage}

\begin{proposition}
\label{proposition.concentration}
Consider an exponential family with countable support $\mX$ and local dependence.
Let $f: \mX \mapsto \mbR$ be Lipschitz with respect to the Hamming metric $d: \mX \times \mX \mapsto \{0, 1, 2, \dots\}$ with Lipschitz coefficient $\norm{f}_{\lip} > 0$ and assume that $\mbE\, f(\bX)$ exists.
Then there exists $C > 0$ such that,
for all $t > 0$,
\beno
\mbP(|f(\bX) - \mbE\, f(\bX)|\, \geq\, t)
\lte 2\, \exp\left(- \dfrac{t^2}{C\, \sum_{k=1}^K {|\mA_k| \choose 2}\, \norm{\mA}_\infty^4\, \norm{f}_{\lip}^2}\right).
\ee
\end{proposition}

\end{samepage}

Proposition \ref{proposition.concentration} covers a wide range of exponential-family random graphs with local dependence.
The assumption that the function of interest is smooth,
in the sense that it is Lipschitz with respect to the Hamming metric,
is common in the concentration-of-measure literature:
see,
e.g.,
the concentration results for dependent random variables by \citet{Sa00},
\citet[][Theorem 4.3, p.\ 75]{Ch05},
and \citet{KoRa08}.
The smoothness assumption can be weakened by using divide-and-conquer strategies:
e.g.,
one may divide the domain of a function of interest into high- and low-probability regions and require the function to be smooth on high-probability regions,
but not on low-probability regions.
Such divide-and-conquer strategies were used by,
e.g.,
\citet{Vu02},
\citet{KiVu04},
and \citet[][Lemma 9]{Yaetal15}.
While exploring divide-and-conquer strategies for\linebreak
exponential-family random graphs with local dependence would be interesting,
Proposition \ref{proposition.concentration} suffices for the purpose of demonstrating that statistical inference for many exponential-family random graphs with local dependence is possible.

Proposition \ref{proposition.concentration} paves the way for concentration results for sufficient statistics.
Proposition \ref{existence.eta.mle} shows that the sufficient statistic vector $\widehat{\bmu(\bta(\bthetas))} = s(\bX)$ is close to the data-generating mean-value parameter vector $\bmu(\bta(\bthetas)) \in \rint(\mM)$ with high probability provided the number of neighborhoods $K$ is large relative to the size of the largest neighborhood $\norm{\mA}_\infty$ and the dimension $\diminf$ of $\bmu(\bta(\bthetas))$.

\begin{samepage}

\begin{proposition}
\label{existence.eta.mle}
Consider a full or non-full, curved exponential family with countable support $\mbX$ and local dependence.
Assume that there exists $A > 0$ such that
\be
\label{smoothness6}
\left\norm{s(\bx_1) - s(\bx_2)\right}_\infty
\lte A\; d(\bx_1, \bx_2)\, \norm{\mA}_\infty
\;\mbox{ for all }
(\bx_1, \bx_2) \in \mX \times \mX.
\ee
Then there exists $C > 0$ such that,
for all deviations of the form $t = \delta\, \sum_{k=1}^K {|\mA_k| \choose 2}^\alpha$ with $\delta > 0$ and $0 \leq \alpha \leq 1$,
\beno
\mbP\left(\norm{\widehat{\bmu(\bta(\bthetas))} - \bmu(\bta(\bthetas))}_2 \geq t\right)
\lte 2\, \exp\left(- \dfrac{\delta^2\; C\; K}{\diminf\, \norm{\mA}_\infty^{4\, (2 - \alpha)}} + \log \diminf\right).
\ee
\end{proposition}

\end{samepage}

\s

The smoothness condition \eqref{smoothness6} of Proposition \ref{existence.eta.mle} is satisfied as long as changing an edge cannot change the within-neighborhood sufficient statistics by more than a constant multiple of $\norm{\mA}_\infty$.
It is verified in Corollary \ref{corollary.curved} in Section \ref{sec:curved}.

\com {\em The relationship between $\alpha$ and sparsity.}
Proposition \ref{existence.eta.mle} shows how the concentration of sufficient statistics depends on the power $\alpha \in [0, 1]$ of deviations of size $t = \delta\, \sum_{k=1}^K {|\mA_k| \choose 2}^\alpha$.
The power $\alpha$ can be interpreted as the level of sparsity of a random graph,
with lower values of $\alpha$ corresponding to higher levels of sparsity.
The conventional definition of a sparse random graph is based on the scaling of the expected number of edges,
i.e.,
the sufficient statistic of Bernoulli random graphs with independent edges.
We use the term sparse random graph to refer to the scaling of the expectations of all sufficient statistics of exponential-family random graphs.
If $\alpha = 1$,
the within-neighborhood subgraphs may be called dense in the sense that the expectations of within-neighborhood sufficient statistics are non-negligible fractions of the number of edge variables ${|\mA_k| \choose 2}$ in neighborhood $\mA_k$ ($k = 1, \dots, K$).
Otherwise,
the within-neighborhood subgraphs may be called sparse.
Note that the interpretation of $\alpha$ in terms of sparsity makes more sense when the neighborhoods grow than when the neighborhoods are bounded above.
But,
regardless of whether the neighborhoods grow,
Proposition \ref{existence.eta.mle} shows how the concentration of sufficient statistics depends on the power $\alpha$.

\com {\em Sharpness.}
The concentration results discussed above are not, 
and cannot be sharp,
because these results cover many models and many dependencies.
It goes without saying that in special cases sharper results can be obtained.
We are here not interested in sharp bounds in special cases,
because the main appeal of the exponential-family framework is that it can capture many dependencies.

\subsection{Maximum likelihood estimators}
\label{sec:curved}

The main idea of Section \ref{sec:idea} is made rigorous in Theorem \ref{nonfull},
which establishes concentration results for maximum likelihood estimators of full and non-full, curved exponential-family random graphs with local dependence.

\begin{theorem}
\label{nonfull}
Consider a full or non-full, curved exponential-family random graph with countable support $\mbX$ and local dependence.
Let
\beno
\bTheta &\subseteq& \{\btheta\, \in\, \mR^{\qqq}:\, \psi(\bta(\btheta)) < \infty\}.
\ee
Assume that $\bthetas \in \interior(\bTheta)$.
Let $\bta: \interior(\bTheta) \mapsto \interior(\fullspace)$ be one-to-one and assume that,
for all $\epsilon > 0$ small enough so that $\mathscr{B}(\bthetas, \epsilon) \subseteq \interior(\bTheta)$,
there exists $\gamma(\epsilon) > 0$ such that,
for all $\btheta \in \bTheta \setminus \mB(\bthetas, \epsilon)$,
we have $\bta(\btheta) \in \fullspace \setminus \mB(\bta(\bthetas),\, \gamma(\epsilon))$.
In addition,
assume that there exist $\delta(\epsilon) > 0$ and $A > 0$ such that,
for all $\bta(\btheta) \in \fullspace \setminus \mB(\bta(\bthetas), \gamma(\epsilon))$,
\be
\label{nonfullid}
\norm{\bmu(\bta(\bthetas)) - \bmu(\bta(\btheta))}_2
&\geq& \delta(\epsilon)\, \dsum_{k=1}^K {|\mA_k| \choose 2}^\alpha
\;\mbox{ for some }\; 0 \leq \alpha \leq 1
\ee
and,
for all $(\bx_1, \bx_2) \in \mX \times \mX$,
\be
\label{smoothness4}
\left\norm{s(\bx_1) - s(\bx_2)\right}_\infty
\lte A\, d(\bx_1, \bx_2)\, \norm{\mA}_\infty.
\ee
Then,
for all $\epsilon > 0$ small enough so that $\mB(\bthetas, \epsilon) \subseteq \interior(\bTheta)$,
there exist $\kappa(\epsilon) > 0$ and $C > 0$ such that 
\be
\nonumber
\mbP\left(\widehat\btheta\; \in\; \mB(\bthetas, \epsilon)\right)
\;\geq\; 1 - 2\, \exp\left(- \dfrac{\kappa(\epsilon)^2\, C\, K}{\diminf\, \norm{\mA}_\infty^{4\, (2 - \alpha)}} + \log \diminf\right).
\ee
If the exponential family is full,
then $\bthetah$ is unique in the event $\widehat\btheta\, \in\, \mB(\bthetas, \epsilon)$.
\end{theorem}

\s

Theorem \ref{nonfull} shows that the estimator $\bthetah$ is close to the data-generating parameter vector $\bthetas$ with high probability provided the number of neighborhoods $K$ is large relative to $\norm{\mA}_\infty$ and $\diminf$.
An application of Theorem \ref{nonfull} to popular curved exponential-family random graphs can be found in Corollary \ref{corollary.curved}.
These concentration results cover inference scenarios with a finite population of nodes and a population graph generated by an exponential-family random graph model,
and assume that the population graph can be observed.
Inference scenarios where the population graph cannot be observed but subgraphs of the population graph can be observed are considered in Section \ref{extendability}.
Asymptotic consistency results can be obtained by allowing the number of neighborhoods $K$ to increase without bound.
We discuss asymptotic consistency results in Remark \ref{rates.of.convgence} following Corollary \ref{corollary.curved}.

Conditions \eqref{nonfullid} and \eqref{smoothness4} of Theorem \ref{nonfull} are verified in Corollary \ref{corollary.curved}.
As pointed out in Section \ref{sec:idea},
the assumptions are satisfied when the map $\bta: \interior(\bTheta) \mapsto \interior(\fullspace)$ is one-to-one and continuous,
but it need not be differentiable.
Condition \eqref{nonfullid} is an identifiability assumption and covers both sparse ($0 \leq \alpha < 1$) and dense ($\alpha = 1$) within-neighborhood subgraphs.
The power $\alpha$ can be interpreted as the level of sparsity of a random graph,
as explained in Section \ref{sec:probability}.
Theorem \ref{nonfull} shows that sparsity comes at a cost,
because the probability of event $\bthetah \not\in \mB(\bthetas, \epsilon)$ decays slower when the within-neighborhood subgraphs are sparse rather than dense.
The fact that sparsity weakens concentration results is well-known in the concentration-of-measure literature on random graphs with independent edges \citep[e.g.,][]{JaRu02,KiVu04}.
Condition \eqref{smoothness4} is a smoothness condition,
which is satisfied as long as changing an edge cannot change the within-neighborhood sufficient statistics by more than a constant multiple of $\norm{\mA}_\infty$.

It is worth noting that in full exponential families the set $\bthetah$ contains a single element when it is non-empty,
whereas in non-full exponential families it may contain more than one element.
A pleasant feature of estimating function \eqref{m.estimating} is that, with high probability, the minimizers $\bthetah$ of \eqref{m.estimating} do not give rise to global minima that are separated by large distances,
under the assumptions made.
The reason is that,
if the set $\bthetah$ contains more than one element,
then all elements of the set $\bthetah$ map to mean-value parameter vectors $\bmu(\bta(\bthetah))$ whose $\ell_2$-distance from $\widehat{\bmu(\bta(\bthetas))}$ is identical and whose $\ell_2$-distance from $\bmu(\bta(\bthetas))$ is bounded above by
\beno
\label{triangle0}
\norm{\bmu(\bta(\bthetah)) - \bmu(\bta(\bthetas))}_2
\lte 2\; \norm{\widehat{\bmu(\bta(\bthetas))} - \bmu(\bta(\bthetas))}_2,
\ee
as explained in Section \ref{sec:idea}.
By Proposition \ref{existence.eta.mle},
$\widehat{\bmu(\bta(\bthetas))}$ is close to $\bmu(\bta(\bthetas))$ with high probability provided the number of neighborhoods $K$ is sufficiently large.
Therefore,
all elements of the set $\bmu(\bta(\bthetah))$ are close to $\bmu(\bta(\bthetas))$ with high probability and hence,
by the identifiability conditions of Theorem \ref{nonfull},
all elements of the set $\bthetah$ are close to $\bthetas$ with high probability.

\s

{\em Applications.}
We present two applications of Theorem \ref{nonfull}.
An application to canonical exponential-family random graphs can be found in Appendix \ref{sec:canonical} \citep[see the supplement,][]{Sc15}.
Here,
we focus on curved exponential-family random graphs with geometrically weighted model terms,
which are popular in practice \citep[e.g.,][]{ergm.book} but are challenging on theoretical grounds.

As a specific example,
consider curved exponential-family random graphs with support $\mbX = \{0, 1\}^{\sum_{k=1}^K \binom{|\mA_k|}{2}}$ and within-neighbor\-hood edge and geometrically weighted edgewise shared partner terms \citep{SnPaRoHa04,HuGoHa08}.
Such models are based on sufficient statistics of the form 
\beno
\label{ss}
s_{k,1}(\bx_k) 
\= \dsum_{a \in \mA_k\, <\, b \in \mA_k}\, x_{a,b}\s
\\
s_{k,i+1}(\bx_k)
\= \dsum_{a \in \mA_k\, <\, b \in \mA_k} x_{a,b}\; f_{a,b,i}(\bx_k),\;\;\;
i = 1, \dots, |\mA_k| - 2,
\ee
where $f_{a,b,i}(\bx_k) = \one(\sum_{c \in \mA_k,\, c \neq a,b} x_{a,c} \, x_{b,c}$ $=$ $ i)$ is an indicator function,
which is $1$ if nodes $a$ and $b$ are both connected to $i$ other nodes in neighborhood $\mA_k$ and is $0$ otherwise ($k = 1, \dots, K$).
The natural parameters are of the form
\beno
\label{np}
\eta_{k,1}(\btheta) 
\= \theta_1\s
\\
\eta_{k,i+1}(\btheta)
\= \exp(\vartheta)\, \left[1 - \left(1 - \exp(-\vartheta)\right)^i\right],\;\;\;
i = 1, \dots, |\mA_k| - 2,
\ee
where $\vartheta > 0$ controls the rate of decay of the geometric sequence $(1 - \exp(-\vartheta))^i$, $i = 1, 2, \dots$ ($k = 1, \dots, K$).
For convenience, 
we consider here the parameterization $\theta_2 = \exp(-\vartheta) \in (0, 1)$,
so that $\btheta = (\theta_1, \theta_2) \in \mR \times (0, 1)$.
Such model terms are called geometrically weighted terms,
because the natural parameters $\eta_{k,i+1}(\btheta)$ are based on the geometric sequence $(1 - \exp(-\vartheta))^i$, $i = 1, 2, \dots$

While complicated,
such models are able to capture transitivity in neighborhoods.
As explained in Section \ref{sec:introduction},
transitivity is one of the more interesting network phenomena,
and induces dependence among edges.
There are many models of transitivity,
some of which are well-posed while others are ill-posed.
An example of a model that is ill-posed in the large-graph limit is the so-called triangle model \citep[e.g.,][]{Jo99,Ha03p,Sc09b,BaBrSl11,ChDi11}.
The triangle model is a canonical exponential-family random graph model with the number of edges and triangles as sufficient statistics;
note that a triangle in a random graph corresponds to three distinct nodes such that all three pairs of nodes are connected by edges.
Compared with the triangle model,
the curved exponential-family random graph model described above makes more reasonable assumptions:
\bi
\item The curved exponential-family random graph model exploits multilevel structure to constrain the dependence among edges induced by transitivity to neighborhoods,
i.e., subsets of nodes.
By contrast,
the triangle model does not restrict transitivity to subsets of nodes,
and allows each edge to depend on many other edges in the random graph.
\item The curved exponential-family random graph model implies that within neighborhoods,
for each pair of nodes,
the value added by additional triangles to the log odds of the conditional probability of an edge decays at a geometric rate \citep[e.g.,][]{Hu08,StScBoMo18}.
As a result,
the model encourages triangles within neighborhoods,
but discourages too many of them.
By contrast,
the added value of additional triangles under the triangle model is constant,
so that the triangle model with a positive triangle parameter places more probability mass on graphs with more triangles (among graphs with the same number of edges).
\ei
While the problematic assumptions underlying the triangle model lead to undesirable behavior in large random graphs \citep[e.g.,][]{Jo99,Ha03p,Sc09b,BaBrSl11,ChDi11},
curved exponential-family random graphs with geometrically weighted edgewise shared partner terms have turned out to be well-behaved \citep[e.g.,][]{HuGoHa08,Sc09b} and have been widely used \citep[see, e.g.,][]{SnPaRoHa04,HuHa04,Hu08,HuGoHa08,ergm.book,StScBoMo18}.
A full-fledged discussion of these complex models is beyond the scope of our paper.
We therefore refer the interested reader to the above-cited literature and focus here on concentration and consistency results.

Curved exponential-family random graphs with within-neighborhood edge and geometrically weighted edgewise shared partner terms are popular in practice but are challenging on theoretical grounds,
for several reasons.
First,
the dimension of the natural parameter vector $\bta(\btheta) \in \mR^{\norm{\mA}_\infty - 1}$ is an increasing function of the number of nodes in the largest neighborhood(s), 
$\norm{\mA}_\infty$.
Second,
the natural parameter vector $\bta(\btheta) \in \mR^{\norm{\mA}_\infty - 1}$ is a non-affine function of a lower-dimensional parameter vector $\btheta \in \mR \times (0, 1)$.
Third,
the mean-value parameter vector $\bmu(\bta(\btheta))$ is not available in closed form.
Finally,
the sufficient statistics $s_{2}(\bx), \dots, s_{\norm{\mA}_\infty - 1}(\bx)$ are not monotone functions of graphs,
which complicates the verification of the main assumption of Theorem \ref{nonfull},
as mentioned in Section \ref{sec:idea}.
Despite these challenges,
it is possible to verify all conditions of Theorem \ref{nonfull} and obtain the following concentration result.
It shows that the estimator $\bthetah$ is close to the data-generating parameter vector $\bthetas$ with high probability provided $K$ is large relative to $\norm{\mA}_\infty^6\, \log \norm{\mA}_\infty$.

\s

\begin{samepage}

\begin{corollary}
\label{corollary.curved}
Consider a curved exponential-family random graph with within-neighborhood edge and geometrically weighted edgewise shared partner terms.
Let $\bTheta = \mR \times (0, 1)$ and assume that $\bthetas \in \interior(\bTheta)$.
Then all conditions of Theorem \ref{nonfull} are satisfied and hence,
for all $\epsilon > 0$ small enough so that $\mB(\bthetas, \epsilon) \subseteq \interior(\bTheta)$,
there exist $\kappa(\epsilon) > 0$ and $C > 0$ such that 
\be
\nonumber
\mbP\left(\widehat\btheta\; \in\; \mB(\bthetas, \epsilon)\right)
&\geq& 1 - 2\, \exp\left(- \dfrac{\kappa(\epsilon)^2\, C\, K}{\norm{\mA}_\infty^6} + \log\, \norm{\mA}_\infty\right),
\ee
provided $|\mA_k| \geq 4$ ($k = 1, \dots, K$) and $K \geq 2$.
\end{corollary}

\end{samepage}

\s

Corollary \ref{corollary.curved} is the first concentration result for estimators of exponential-family random graphs with dependence among edges induced by transitivity,
one of the more interesting network phenomena.
In addition,
it is the first concentration result for curved exponential-family random graphs with geometrically weighted model terms,
which are popular in practice \citep[e.g.,][]{SnPaRoHa04,HuGoHa08}.
The concentration result assumes that each neighborhood $\mA_k$ consists of $|\mA_k| \geq 4$ nodes,
because $\bthetas$ is not identifiable when $|\mA_k| \leq 3$ ($k = 1, \dots, K$).
As mentioned above,
the set $\bthetah$ may contain more than one element,
but all elements of the set $\bthetah$ are close to $\bthetas$ with high probability provided $K$ is large relative to $\norm{\mA}_\infty^6\, \log \norm{\mA}_\infty$.
Concentration results for other curved exponential-family random graphs with geometrically weighted model terms \citep[e.g.,][]{SnPaRoHa04,HuGoHa08} can be established along the same lines.

\com {\em Asymptotic consistency results.}
\label{rates.of.convgence}
As pointed out in Section \ref{sec:introduction},
we state all theoretical results for finite populations of nodes,
because in practice all populations are finite.
Asymptotic consistency results can be obtained by allowing the number of neighborhoods $K$ to grow without bound.
If there exists a universal constant $C > 0$ such that $|\mA_k| < C$ ($k = 1, 2, \dots$),
then the main idea described in Section \ref{sec:idea} along with the concentration results in Section \ref{sec:probability} imply that $\bthetah$ is a consistent estimator of $\bthetas$ with rate of convergence $K^{1/2}$.
As the units of statistical analysis are neighborhoods,
the rate $K^{1/2}$ resembles the rate in classical statistical problems where the rate is the square root of the sample size,
albeit with two notable differences:
first,
the units are subsets of nodes (neighborhoods) rather than nodes or edges;
and,
second,
the sizes of units are not identical,
but the size of the largest unit is a constant multiple of the size of the smallest.
Last,
but not least,
it is possible to obtain asymptotic consistency results when $K$ grows and the neighborhoods grow with $K$,
which implies that $\norm{\mA}_\infty$ grows with $K$.
Then,
as long as $K$ grows faster than $\norm{\mA}_\infty^6\, \log\, \norm{\mA}_\infty$ in the sense that $K\, /\, (\norm{\mA}_\infty^6\, \log\, \norm{\mA}_\infty) \tends \infty$,
the probability of event $\widehat\btheta\, \in\, \mB(\bthetas, \epsilon)$ tends to $1$.

\subsection{$M$-estimators, correct and incorrect model specifications}
\label{sec:model.estimation}

The concentration and consistency results for maximum likelihood estimators in Section \ref{sec:curved} are special cases of more general results for $M$-estimators.
To demonstrate,
we introduce a natural class of $M$-estimators in Appendix \ref{assumptions},
which includes both likelihood- and moment-based estimators,
and present concentration results in Appendix \ref{mainresults} along with an application to misspecified models with omitted covariate terms.
These results cover both correct and incorrect model specifications,
as the example with omitted covariate terms demonstrates.
Due to space restrictions,
we provide details in Appendix \ref{supplement.mestimators} \citep[see the supplement,][]{Sc15}.

\section{Extendability and subgraph-to-graph estimators}
\label{extendability}

A question that has been asked about exponential-family random graphs is whether it is possible to extend, in a well-defined sense, an exponential-family random graph with a given set of nodes to an exponential-family random graph with more nodes \citep{ShRi11,CrDe15,LaRiSa17}.
We show that multilevel structure helps extend an exponential-family random graph with a given set of neighborhoods to an exponential-family random graph with more neighborhoods (Section \ref{sec:extendability}) and hence facilitates subgraph-to-graph estimation (Section \ref{incompletemle}).
The importance of these results lies in the fact that subgraph-to-graph estimation for exponential-family random graphs is believed to be difficult \citep[e.g.,][]{ShRi11},
but our results demonstrate that additional structure facilitates it.

\subsection{Extendability}
\label{sec:extendability}

While many exponential-family random graphs with a given set of nodes cannot be extended to exponential-family random graphs with more nodes \citep{ShRi11,CrDe15,LaRiSa17},
an exponential-family random graph with a given set of neighborhoods can be extended to an exponential-family random graph with more neighborhoods.

To demonstrate,
consider a population graph $(\bX_{\mL}, \bY_{\mL})$ with a set of neighborhoods $\mL = \{\mA_1, \dots, \mA_L\}$,
where $\bX_{\mL} \in \mX_{\mL}$ and $\bY_{\mL} \in \mY_{\mL}$ denote the sequences of within- and between-neighborhood edge variables based on the set of neighborhoods $\mL$,
respectively.
As before,
assume that $\bX_{\mL}$ is governed by an exponential family with countable support $\mX_{\mL}$ and local dependence,
with neighborhood-dependent natural parameters $\eta_{\mA,i}(\btheta) = \eta_i(\btheta)$ and sufficient statistics $s_{\mA,i}(\bx_{\mA})$ ($i = 1, \dots, \ppp_{\mA}$, $\mA \in \mL$).
Therefore,
the exponential family can be reduced to an exponential family with natural parameter vector
\be
\label{extendable.eta}
\bta(\btheta)
\= \left(\eta_1(\btheta), \dots, \eta_{\diminf}(\btheta)\right)
\ee
and sufficient statistic vector
\beno
s(\bx_{\mL})
\= \left(s_1(\bx_{\mL}), \dots, s_{\diminf}(\bx_{\mL})\right),
\ee
where $s_i(\bx_{\mL}) = \sum_{\mA \in \mL} s_{\mA,i}(\bx_{\mA})$ ($i = 1, \dots, \diminf$) and $\diminf = \max_{\mA \in \mL} \ppp_{\mA}$.

Consider a subgraph $(\bX_{\mK}, \bY_{\mK})$ induced by a subset of neighborhoods $\mK \subseteq \mL$.
Then the subgraph $(\bX_{\mK}, \bY_{\mK})$ with subset of neighborhoods $\mK$ is extendable to the population graph $(\bX_{\mL}, \bY_{\mL})$ with set of neighborhoods $\mL \supset \mK$ as follows.

\begin{proposition}
\label{proposition.extendability}
Consider a full or non-full, curved exponential-family random graph with 
set of neighborhoods $\mL$, 
countable support $\mbX_{\mL}$,
and local dependence.
Assume that,
for all $\by_{\mL} \in \mY_{\mL}$,
\beno
\mbP(\bY_{\mL} = \by_{\mL}) 
\= \dprod_{\mbC \in \mL,\; \mbD \in \mL,\; \mbC \neq \mbD} \mbP(\bY_{\mbC,\mbD} = \by_{\mbC,\mbD}),
\ee
where $\bY_{\mbC,\mbD} = (Y_{i,j})_{i\in\mbC,\, j\in\mbD}$.
Then, 
for all $\btheta\in\bTheta \subseteq \{\btheta \in \mR^{\qqq}:\; \psi_{\mL}(\bta(\btheta)) < \infty\}$,
all $\mK \subseteq \mL$,
and all $\bx_{\mK} \in \mX_{\mK}$ and $\by_{\mK} \in \mY_{\mK}$, 
\beno
&& \mbP_{\bta(\btheta)}(\bX_{\mK} = \bx_{\mK},\; \bY_{\mK} = \by_{\mK},\; \bX_{\mL \setminus \mK} \in \mX_{\mL \setminus \mK},\; \bY_{\mL \setminus \mK} \in \mbY_{\mL \setminus \mK})\s
\\
\= \mbP_{\bta(\btheta)}(\bX_{\mK} = \bx_{\mK},\; \bY_{\mK} = \by_{\mK}),
\ee
where 
\beno
\psi_{\mL}(\bta(\btheta)) 
\,=\; \dsum_{\mA\in\mL} \psi_{\mA}(\bta(\btheta)) 
\;=\; \dsum_{\mA\in\mL} \log\dsum_{\bx_{\mA} \in \mX_{\mA}} \exp\left(\langle\bta(\btheta),\, s_{\mA}(\bx_{\mA})\rangle\right)\, \nu_{\mA}(\bx_{\mA}).
\ee
The marginal density of a subgraph $\bx_{\mK} \in \mX_{\mK}$ of $\bx_{\mL} \in \mX_{\mL}$ induced by $\mK \subseteq \mL$ is an exponential-family density with support $\mX_{\mK}$ and local dependence:
\beno
\dsum_{\bx_{\mL\setminus\mK} \in \mbX_{\mL\setminus\mK}} p_{\bta(\btheta)}(\bx_{\mL})
= p_{\bta(\btheta)}(\bx_{\mK})
= \exp\left(\left\langle\bta(\btheta), s(\bx_{\mK})\right\rangle - \psi_{\mK}(\bta(\btheta))\right) \nu_{\mK}(\bx_{\mK}),
\ee
where $\psi_{\mK}(\bta(\btheta)) = \sum_{\mA\in\mK} \psi_{\mA}(\bta(\btheta))$,
$\bta(\btheta) = (\eta_1(\btheta), \dots, \eta_{\diminf}(\btheta))$,
$s(\bx_{\mK}) = (\sum_{\mA \in \mK} s_{\mA,1}(\bx_{\mA}), \dots, \sum_{\mA \in \mK} s_{\mA,\diminf}(\bx_{\mA}))$,
and $\nu_{\mK}(\bx_{\mK}) = \prod_{\mA \in \mK} \nu_{\mA}(\bx_{\mA})$.
\end{proposition}

\s

Thus,
in the above-mentioned sense,
the exponential-family random graph induced by a subset of neighborhoods $\mK$ can be extended to the exponential-family random graph with set of neighborhoods $\mL \supset \mK$.
A more restrictive result was proved by \citet[][Theorem 1]{ScHa13}.

\subsection{Subgraph-to-graph estimators}
\label{incompletemle}

The extendability of exponential-family random graphs with multilevel structure discussed in Section \ref{sec:extendability} facilitates subgraph-to-graph estimation.

To demonstrate,
let $\mL$ be the set of neighborhoods of the population graph and assume that $\bx_{\mL} \in \mX_{\mL}$ was generated by an exponential family with countable support $\mX_{\mL}$ and local dependence.
Suppose that it is infeasible to observe $\bx_{\mL} \in \mX_{\mL}$,
but it is feasible to sample a subset of neighborhoods $\mK \subseteq \mL$ and collect data on the subgraphs induced by $\mK \subseteq \mL$.
We assume henceforth that the sampling design is ignorable in the sense of \citet{Ru76} and \citet{HaGi09},
i.e.,
the probability of observing subgraphs does not depend on the unobserved subgraphs.
A simple example is a sampling design that samples neighborhoods at random and collects data on the subgraphs induced by the sampled neighborhoods.

By Proposition \ref{proposition.extendability} and the ignorability of the sampling design \citep{Ru76,HaGi09},
the\linebreak 
observed-data likelihood function based on the observed subgraph $\bx_{\mK} \in \mX_{\mK}$ of $\bx_{\mL} \in \mX_{\mL}$ is proportional to
\beno
\dsum_{\bx_{\mL\setminus\mK}\, \in\, \mbX_{\mL\setminus\mK}}\, p_{\bta(\btheta)}(\bx_{\mL})
\= p_{\bta(\btheta)}(\bx_{\mK}),
&& \bx_{\mL} \in \mX_{\mL},
\ee
where $\bta(\btheta)$ is of the form \eqref{extendable.eta}.
In other words,
maximum likelihood estimation can be based on $p_{\bta(\btheta)}(\bx_{\mK})$.
Motivated by the same considerations we outlined in Section \ref{mle},
we therefore consider an estimating function of the form
\beno
g_{\mK}(\btheta;\, \widehat{\bmu_{\mK}(\bta(\bthetas))}) 
\= \norm{\widehat{\bmu_{\mK}(\bta(\bthetas))} - \bmu_{\mK}(\bta(\btheta))}_2,
\ee
where $\widehat{\bmu_{\mK}(\bta(\bthetas))} = s(\bx_{\mK})$, 
$\bmu_{\mK}(\bta(\btheta)) = \mbE_{\bta(\btheta)}\, s(\bX_{\mK})$,
and $s(\bx_{\mK})$ is defined in Proposition \ref{proposition.extendability}.
The data-generating parameter vector $\bthetas$ of the population graph can hence be estimated by the estimator $\bthetah_{\mK}$ based on the observed subgraph $\bx_{\mK} \in \mX_{\mK}$ of $\bx_{\mL} \in \mX_{\mL}$:
\beno
\bthetah_{\mK}
\= \left\{\btheta \in \bTheta:\;\; g_{\mK}(\btheta;\; \widehat{\bmu_{\mK}(\bta(\bthetas))})\;\; =\;\; \inf\limits_{\dot\btheta\in\bTheta}\; g_{\mK}(\dot\btheta;\; \widehat{\bmu_{\mK}(\bta(\bthetas))})\right\}.
\ee
The following concentration result shows that,
with high probability,
the estimator $\bthetah_{\mK}$ based on the observed subgraph induced by $\mK \subseteq \mL$ is close to the data-generating parameter vector $\bthetas$ of the population graph as long as the number of sampled neighborhoods $|\mK|$ is large relative to $\norm{\mL}_\infty = \max_{\mA \in \mL} |\mA|$ and $\diminf = \max_{\mA \in \mL} \ppp_{\mA}$.

\begin{theorem}
\label{incomplete}
Consider a full or non-full, curved exponential-family random graph with set of neighborhoods $\mL$, countable support $\mbX_{\mL}$, and local dependence.
Let
\beno
\bTheta &\subseteq& \{\btheta\, \in\, \mR^{\qqq}:\, \psi_{\mL}(\bta(\btheta)) < \infty\}.
\ee
Assume that $\bthetas \in \interior(\bTheta)$ and that,
for all $\epsilon > 0$ small enough so that $\mathscr{B}(\bthetas, \epsilon) \subseteq \interior(\bTheta)$,
there exists $\gamma(\epsilon) > 0$ such that,
for all $\btheta \in \bTheta \setminus \mB(\bthetas, \epsilon)$,
we have $\bta(\btheta) \in \fullspace \setminus \mB(\bta(\bthetas),\, \gamma(\epsilon))$.
In addition,
assume that there exist $\delta(\epsilon) > 0$ and $A > 0$ such that,
for all $\mK \subseteq \mL$ and all $\bta(\btheta) \in \fullspace \setminus \mathscr{B}(\bta(\bthetas), \gamma(\epsilon))$,
\be
\label{nonfullid8}
\norm{\bmu_{\mK}(\bta(\bthetas)) - \bmu_{\mK}(\bta(\btheta))}_2
\;\geq\; \delta(\epsilon)\, \dsum_{\mA \in \mK} {|\mA| \choose 2}^\alpha
\mbox{ for some } 0 \leq \alpha \leq 1
\ee
and,
for all $\mK \subseteq \mL$ and all $(\bx_1, \bx_2) \in \mX_{\mK} \times \mX_{\mK}$,
\be
\label{smoothness8}
\left\norm{s(\bx_1) - s(\bx_2)\right}_\infty
\lte A\, d(\bx_1, \bx_2)\, \norm{\mK}_\infty,
\ee
where $\norm{\mK}_\infty = \max_{\mA \in \mK} |\mA|$.
Then,
for all $\epsilon > 0$ small enough so that $\mB(\bthetas, \epsilon) \subseteq \interior(\bTheta)$,
there exist $\kappa(\epsilon) > 0$ and $C > 0$ such that
\be
\nonumber
\mbP\left(\widehat\btheta_{\mK}\; \in\; \mB(\bthetas, \epsilon)\right)
\;\geq\; 1 - 2\, \exp\left(- \dfrac{\kappa(\epsilon)^2\, C\, |\mK|}{\diminf^{}\, \norm{\mL}_\infty^{4\, (2 - \alpha)}} + \log \diminf\right).
\ee
If the exponential family is full,
then $\bthetah_{\mK}$ is unique in the event $\widehat\btheta_{\mK}\, \in\, \mB(\bthetas, \epsilon)$.
\end{theorem}

\s

Theorem \ref{incomplete} shows that there are costs associated with observing a subset of neighborhoods $\mK \subseteq \mL$ rather than the whole set of neighborhoods $\mL$ of the population graph:
The probability of event $\widehat\btheta_{\mK}\; \not\in\; \mB(\bthetas, \epsilon)$ decays with the number of sampled neighborhoods $|\mK|$ and is hence lowest when the whole set of neighborhoods $\mL$ of the population graph is sampled.

As a specific example,
consider the main example of Section \ref{sec:curved}:
curved\linebreak 
exponential-family random graphs with within-neighborhood edge and geometrically weighted edgewise shared partner terms.

\s

\begin{corollary}
\label{corollary.curved.incomplete}
Consider a curved exponential-family random graph with set of neighborhoods $\mL$, countable support $\mbX_{\mL}$, and local dependence induced by within-neighborhood edge and geometrically weighted edgewise shared partner terms.
Let $\bTheta = \mR \times (0, 1)$ and assume that $\bthetas \in \interior(\bTheta)$.
Then all conditions of Theorem \ref{incomplete} are satisfied and hence,
for all $\epsilon > 0$ small enough so that $\mB(\bthetas, \epsilon) \subseteq \interior(\bTheta)$,
there exist $\kappa(\epsilon) > 0$ and $C > 0$ such that 
\be
\nonumber
\mbP\left(\widehat\btheta_{\mK}\; \in\; \mB(\bthetas, \epsilon)\right)
\;\geq\; 1 - 2\, \exp\left(- \dfrac{\kappa(\epsilon)^2\, C\, |\mK|}{\norm{\mL}_\infty^6} + \log \norm{\mL}_\infty\right),
\ee
provided $|\mA| \geq 4$ ($\mA \in \mL$) and $|\mK| \geq 2$.
\end{corollary}

\com {\em ``Bad" subsets of neighborhoods $\mK \subseteq \mL$.}
\label{com.bad.event}
Since the neighborhoods need not have the same size,
it is natural to ask whether it is possible to sample a ``bad" subset of neighborhoods $\mK$ with too small or too large neighborhoods,
which could make it challenging to estimate some of the parameters.
However,
the assumptions of Theorem \ref{corollary.curved.incomplete} rule out ``bad" subsets of neighborhoods $\mK$,
for two reasons.
First,
while some neighborhoods may be larger than others,
Theorem \ref{corollary.curved.incomplete} assumes that the neighborhoods are of the same order of magnitude,
as defined in Section \ref{sec:model}.
In other words,
the neighborhoods have similar sizes.
Second,
the conditions of Theorem \ref{incomplete} assume that the model satisfies identifiability and smoothness conditions for all possible subsets of neighborhoods $\mK \subseteq \mL$.
Corollary \ref{corollary.curved.incomplete} shows that, 
in the special case of curved exponential-family random graphs with within-neighborhood edge and geometrically weighted edgewise shared partner terms,
the identifiability conditions require $|\mA| \geq 4$ for all neighborhoods $\mA \in \mL$ of the population graph.
Thus,
no neighborhood can be too small,
and no neighborhood can be too large,
because all neighborhoods are of the same order of magnitude.
As a consequence,
under the stated assumptions,
it is impossible to sample a ``bad" subset of neighborhoods $\mK$ with too small or too large neighborhoods.

\section{Comparison with existing consistency results}
\label{literature}

To compare our concentration and consistency results to existing consistency results,
we focus on\linebreak
exponential-family random graphs with dependent edges.
It is worth noting that there are consistency results for exponential-family random graphs with independence assumptions---see,
e.g.,
\citet[][]{DiChSl11,RiPeFi13,KrKo14}, and \citet{YaLeZh11,Yaetal18}---but such independence assumptions may not be satisfied in applications,
as discussed in Section \ref{sec:introduction}.

Concerning exponential-family random graphs with dependent edges,
\citet{ShRi11} showed that maximum likelihood estimators of natural parameters of fixed dimension are consistent provided exponential-family random graphs satisfy strong extendability or projectability assumptions.
However, 
those\linebreak
projectability assumptions rule out dependencies induced by transitivity and many other interesting network phenomena.
\citet{XiNe11} reported consistency results under weak dependence assumptions,
but did not give any example of an exponential-family random graph with dependent edges that satisfies those assumptions.
\citet{Mu13} showed that consistent estimation of the so-called two-star model is possible,
but those results have not been extended to other exponential-family random graphs.
In addition,
\citet{ShRi11},
\citet{XiNe11},
and \citet{Mu13} focus on consistency results for estimators of natural parameter vectors whose dimensions do not depend on the number of nodes.
By contrast,
we advance the statistical theory of exponential-family random graphs by providing the first concentration and consistency results that cover
\bi
\item a wide range of exponential-family random graphs with dependence among edges induced by transitivity and other interesting network phenomena;
\item curved exponential-family random graphs with dependent edges and parameter vectors whose dimension depends on the number of nodes (Section \ref{mle});
\item maximum likelihood and $M$-estimators (Section \ref{mle} and Appendix \ref{supplement.mestimators});
\item correct and incorrect model specifications (Section \ref{mle} and Appendix \ref{supplement.mestimators});
\item subgraph-to-graph estimators (Section \ref{extendability}).
\ei
These results underscore the importance of additional structure:
It is the additional structure in the form of multilevel structure that facilitates these results.

\section{Simulation results}
\label{sec:simulations}

\begin{figure}[t]
\begin{center}
\includegraphics[scale=.281]{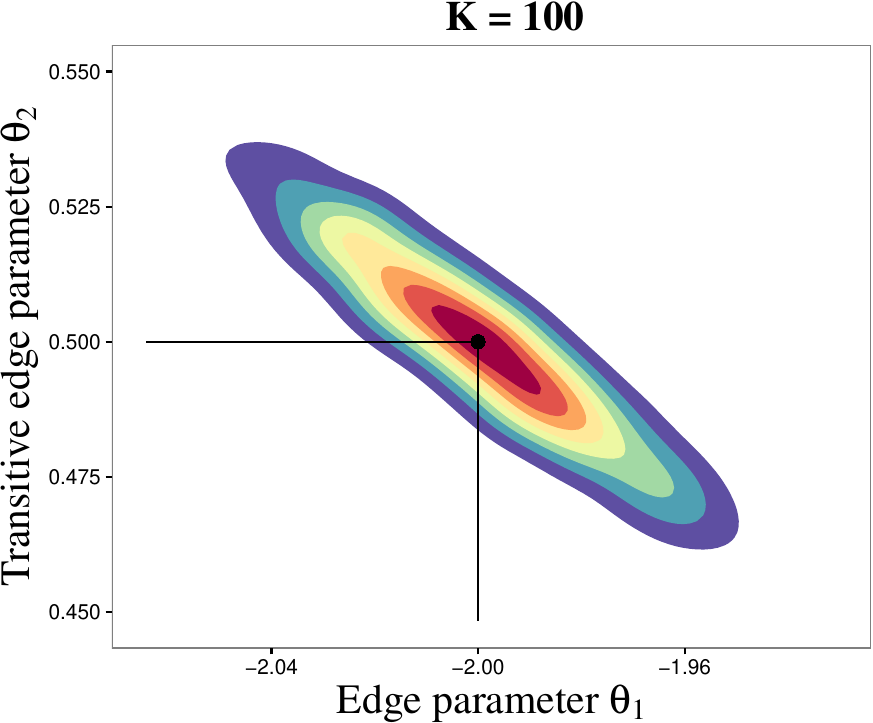}
\includegraphics[scale=.281]{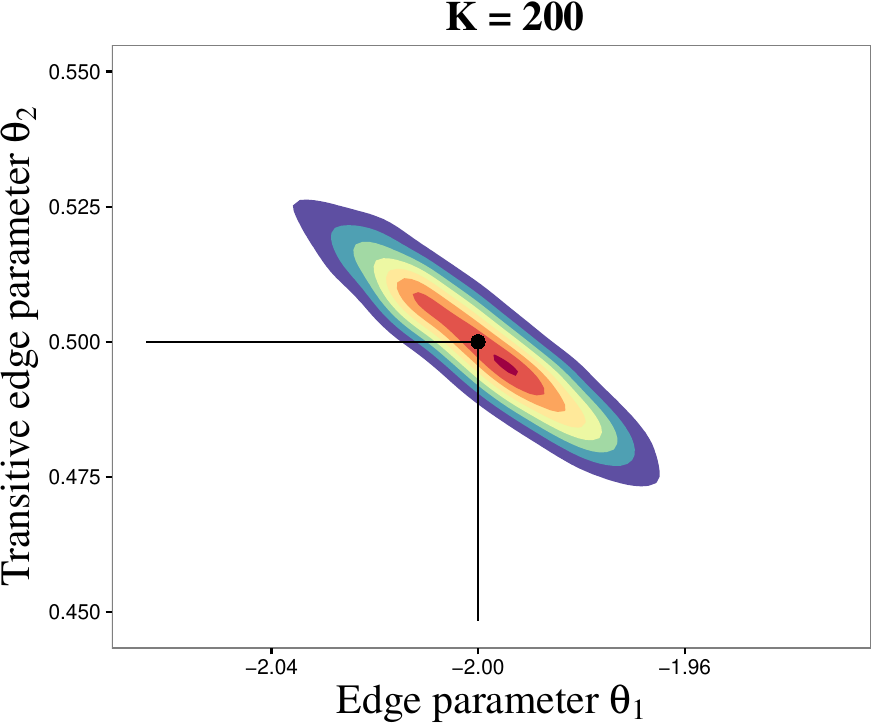}
\includegraphics[scale=.281]{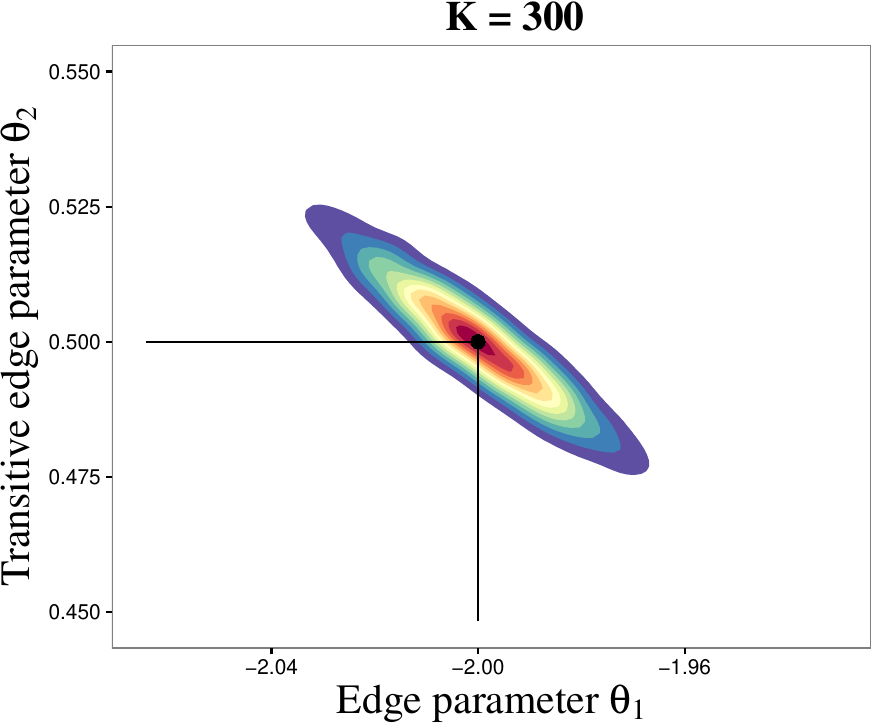}
\caption{
\label{fig:homo}
1,000 estimates of the exponential-family random graph with within-neighborhood edge and transitive edge terms,
where each graph consists of $K = 100$, $200$, and $300$ neighborhoods of size $50$ with natural parameter vectors $\bta_k(\btheta) = (\theta_1, \theta_2)$.
The horizontal and vertical lines indicate the coordinates of the data-generating parameter vector $\bthetas = (\theta_1^\star, \theta_2^\star)$.
}
\end{center}
\end{figure}

To shed light on the finite-graph properties of maximum likelihood estimators,
we generated data from the canonical and curved exponential-family random graphs mentioned in Section \ref{sec:curved}.
We used R package hergm \citep{ScLu15} to generate 1,000 graphs from each model and estimated the data-generating parameter vector by Monte Carlo maximum likelihood estimators \citep{HuGoHa08}.

We first consider canonical exponential-family random graphs with support $\mbX = \{0, 1\}^{\sum_{k=1}^K \binom{|\mA_k|}{2}}$ and local dependence induced by within-neighborhood edge and transitive edge terms \citep{HuKrSc12}.
Within-neighborhood edge and transitive edge terms correspond to neighborhood-dependent natural parameters $\eta_{k,1}(\btheta) = \theta_1$ and\linebreak 
$\eta_{k,2}(\btheta) = \theta_2$ and sufficient statistics $s_{k,1}(\bx_k)$ and $s_{k,2}(\bx_k)$ given by
\be
\nonumber
s_{k,1}(\bx_k)
&=& \dsum_{i\in\mA_k\, <\, j\in\mA_k} x_{i,j}\s
\\
s_{k,2}(\bx_k)
&=& \dsum_{i\in\mA_k\, <\, j\in\mA_k} x_{i,j}\, \max\limits_{h \in \mA_k,\, h \neq i,j} x_{i,h} \; x_{j,h},
\ee
where $k = 1, \dots, K$.
It is worth noting that the number of transitive edges is not the same as the number of triangles.
A discussion of the model along with concentration results for maximum likelihood estimators can be found in Appendix \ref{sec:canonical} \citep[see the supplement,][]{Sc15}.
Figure \ref{fig:homo} shows 1,000 estimates of the exponential-family random graph with within-neighborhood edge and transitive edge terms,
where each graph consists of $K = 100$, $200$, and $300$ neighborhoods of size $50$ with natural parameter vectors $\bta_k(\btheta) = (\theta_1, \theta_2)$ ($k = 1, \dots, K$).
The figure suggests that the probability mass of estimators becomes more and more concentrated in a neighborhood of the data-generating parameters as the number of neighborhoods $K$ increases from $100$ to $300$, 
demonstrating that the concentration results in Section \ref{mle} are manifest when $K$ is in the low hundreds and $\norm{\mA}_\infty = 50$.

\begin{figure}[t]
        \begin{center}
                \includegraphics[scale=.281]{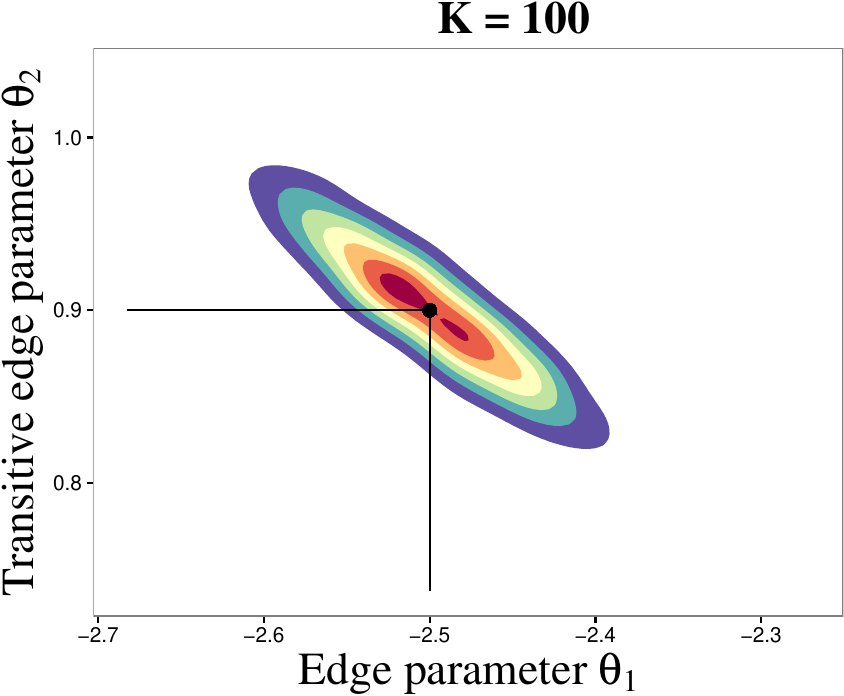}
                \includegraphics[scale=.281]{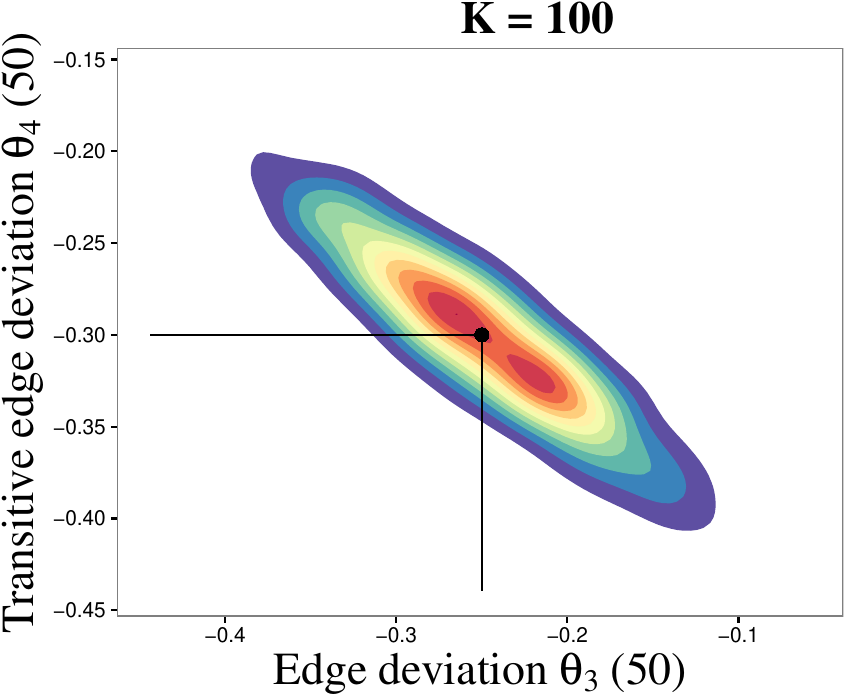}
                \includegraphics[scale=.281]{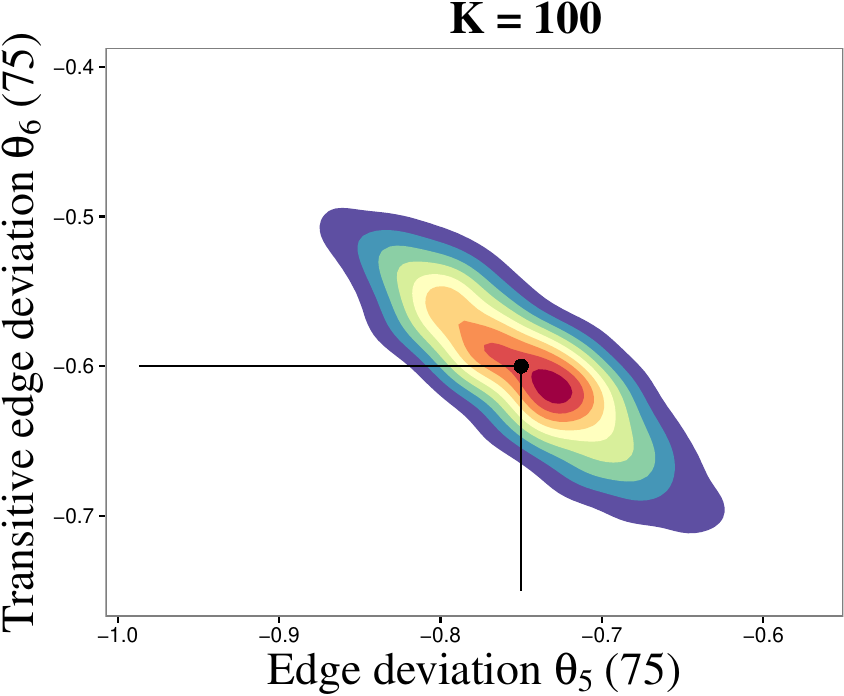}
                \caption{
                        \label{fig:hetero_lin_dep}
1,000 estimates of the exponential-family random graph with within-neighborhood edge and transitive edge terms,
where each graph consists of $33$ neighborhoods of size $25$ with natural parameter vectors $\bta_k(\btheta) = (\theta_1, \theta_2)$,
$34$ neighborhoods of size $50$ with natural parameter vectors $\bta_k(\btheta) = (\theta_1 + \theta_3, \theta_2 + \theta_4)$,
and $33$ neighborhoods of size $75$ with natural parameter vectors $\bta_k(\btheta) = (\theta_1 + \theta_5, \theta_2 + \theta_6)$.
The horizontal and vertical lines indicate the coordinates of the data-generating parameter vector $\bthetas = (\theta_1^\star, \theta_2^\star, \theta_3^\star, \theta_4^\star, \theta_5^\star, \theta_6^\star)$.
             }	
        \end{center}
\end{figure}
Figure \ref{fig:hetero_lin_dep} sheds light on the performance of a simple form of a size-dependent parameterization that allows small and large neighborhoods to have different parameters.
We consider exponential-family random graphs with within-neighborhood edge and transitive edge terms,
where each graph consists of $33$ neighborhoods of size $25$ with natural parameter vectors $\bta_k(\btheta) = (\theta_1, \theta_2)$ ($k = 1, \dots, 33$), 
$34$ neighborhoods of size $50$ with natural parameter vectors $\bta_k(\btheta) = (\theta_1 + \theta_3, \theta_2 + \theta_4)$ ($k = 34, \dots, 67$), 
and $33$ neighborhoods of size $75$ with natural parameter vectors $\bta_k(\btheta) = (\theta_1 + \theta_5, \theta_2 + \theta_6)$ ($k = 68, \dots, 100$). 
Figure \ref{fig:hetero_lin_dep} demonstrates that the estimates of the baseline edge and transitive edge parameters $\theta_1$ and $\theta_2$ tend to be closer to the data-generating parameters than the deviation parameters $\theta_3, \theta_4, \theta_5$, and $\theta_6$,
because $\theta_1$ and $\theta_2$ are estimated from all neighborhoods whereas $\theta_3, \theta_4, \theta_5$, and $\theta_6$ are estimated from a subset of neighborhoods.

\begin{figure}[t]
\begin{center}
\includegraphics[scale=.225]{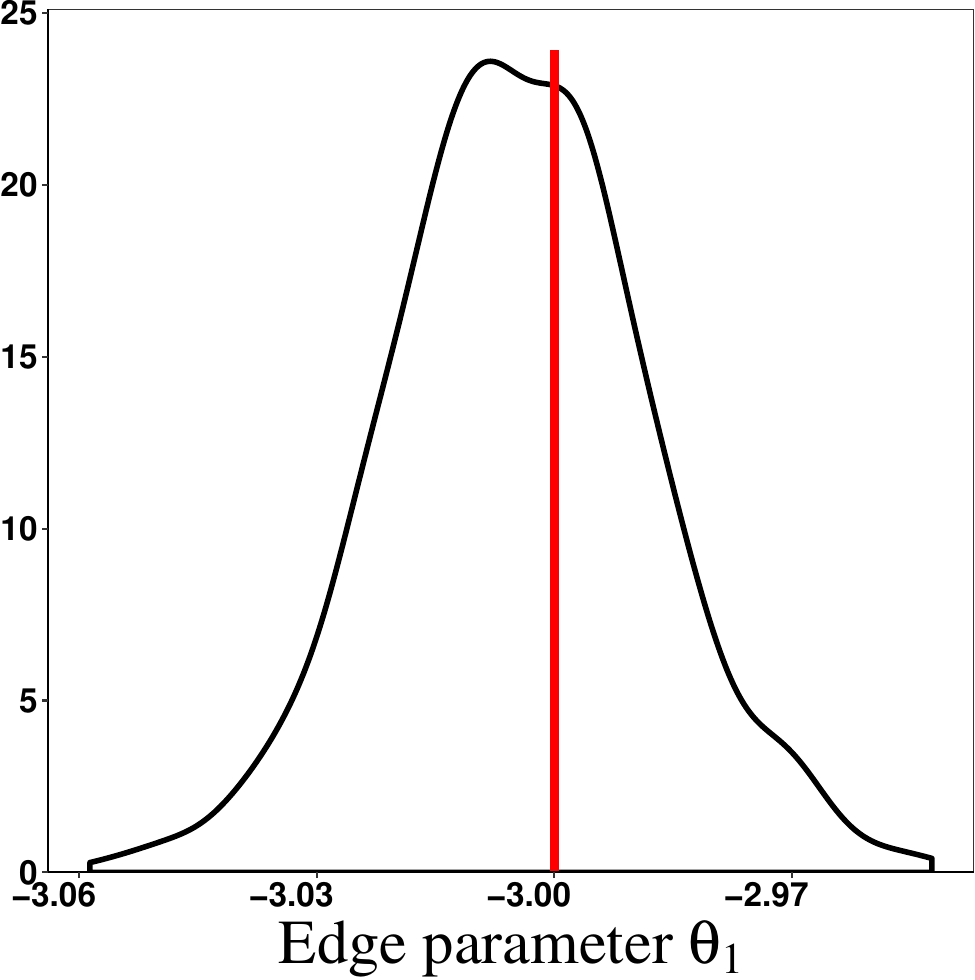}
\hspace{.25cm}
\includegraphics[scale=.225]{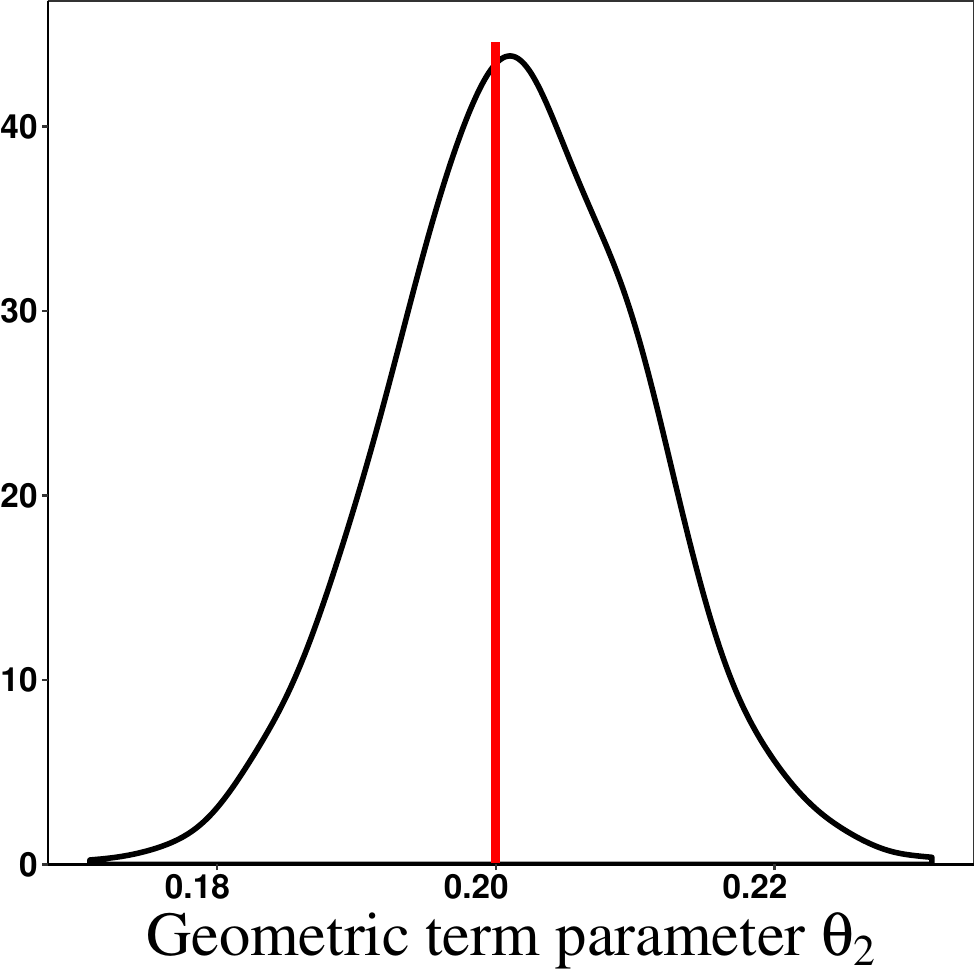}
\caption{
\label{fig:curved}
1,000 estimates of the curved exponential-family random graph with within-neighborhood edge and geometrically weighted edgewise shared partner terms,
where each graph consists of $K = 100$ neighborhoods of size $50$ with natural parameter vectors $\bta_k(\btheta) = (\theta_1,\, \eta_{k,1}(\theta_2),\, \dots,\, \eta_{k,48}(\theta_2))$.
The vertical lines indicate the coordinates of the data-generating parameter vector $\bthetas = (\theta_1^\star, \theta_2^\star)$.
}
\end{center}
\end{figure}

Figure \ref{fig:curved} shows 1,000 estimates of the curved exponential-family random graph with within-neighborhood edge and geometrically weighted edgewise shared partner terms described in Section \ref{sec:curved}.
Each graph consists of $K = 100$ neighborhoods of size $50$ with natural parameter vectors $\bta_k(\btheta) = (\theta_1,\, \eta_{k,1}(\theta_2),\, \dots,\, \eta_{k,48}(\theta_2))$,
where $\theta_1$ is the natural parameter of the edge term and $\eta_{k,1}(\theta_2), \dots, \eta_{k,48}(\theta_2)$ are the natural parameters of the geometrically weighted edgewise shared partner term ($k = 1, \dots, 100$).
The figure shows that the probability mass of the estimators is concentrated in a small neighborhood of the data-generating parameters.

\section{Discussion}
\label{sec:discussion}

We have taken constructive steps to demonstrate that statistical inference for exponential-family random graphs with dependence among edges induced by transitivity and other interesting network phenomena is possible,
provided additional structure in the form of multilevel structure is available.
The theoretical results reported here underscore the importance of additional structure.
In practice,
many other forms of additional structure exist and could be used for the purpose of facilitating statistical inference for exponential-family random graphs (e.g., other forms of multilevel structure or spatial structure).

Last,
but not least,
while we have focused here on theoretical results showing that multilevel structure facilitates statistical inference,
it is worth noting that multilevel structure has computational benefits as well:
The contributions of neighborhoods to estimating functions---e.g., the expectations of within-neighborhood sufficient statistics---may be computed or approximated by exploiting parallel computing on computing clusters.

\section*{Acknowledgements}

We acknowledge support from the National Science Foundation (NSF awards DMS-1513644 and DMS-1812119).

\section*{Supplementary materials}

All results are proved in the supplement\ghoster{ \citep{Sc15}}.
In addition,
the supplement contains concentration results for $M$-estimators.
These results cover correct and incorrect model specifications.

\newpage

\bibliographystyle{asa1}

\bibliography{base}

\pagebreak

\makeatletter

\setcounter{page}{1}

\setcounter{section}{0}

\title{Supplement:\\ \longtitle}

\begin{center}
{\normalfont\textsc{By Michael Schweinberger And Jonathan Stewart}}\s
\\
{\normalfont\em Rice University}
\end{center}

\maketitle

\begin{appendix}

\s

We first present an additional application of concentration results for maximum likelihood estimators in Appendix \ref{sec:canonical} and discuss concentration results for $M$-estimators in Appendix \ref{supplement.mestimators}.
We then prove the main concentration results for random graphs with dependent edges in Appendix \ref{sec:proofs.concentration} and the main concentration results for maximum likelihood estimators,
$M$-estimators,
and subgraph-to-graph estimators in Appendices \ref{mle.proofs}, 
\ref{sec:proofs.consistency},
and \ref{appendix.incomplete},
respectively.
Auxiliary lemmas are proved in Appendix \ref{m.id}.

\section{Concentration: maximum likelihood estimators}
\label{sec:canonical}

In addition to the application to curved exponential-family random graphs in Section \ref{sec:curved},
we present here an application to canonical exponential-family random graphs that induces dependence among edges through transitivity.

We consider canonical exponential-family random graphs with support $\mbX = \{0, 1\}^{\sum_{k=1}^K \binom{|\mA_k|}{2}}$ and local dependence induced by within-neighborhood edge and transitive edge terms \citep{HuKrSc12}.
Within-neighborhood edge and transitive edge terms correspond to neighborhood-dependent natural parameters $\eta_{k,1}(\btheta) = \theta_1$ and\linebreak 
$\eta_{k,2}(\btheta) = \theta_2$ and sufficient statistics $s_{k,1}(\bx_k)$ and $s_{k,2}(\bx_k)$ given by
\be
\nonumber
s_{k,1}(\bx_k)
&=& \dsum_{i\in\mA_k\, <\, j\in\mA_k} x_{i,j}\s
\\
s_{k,2}(\bx_k)
&=& \dsum_{i\in\mA_k\, <\, j\in\mA_k} x_{i,j}\, \max\limits_{h \in \mA_k,\, h \neq i,j} x_{i,h} \; x_{j,h},
\ee
where $k = 1, \dots, K$.

\begin{corollary}
\label{corollary.transitive1}
Consider an exponential-family random graph with\linebreak 
within-neighborhood edge and transitive edge terms.
Let $\bthetas \in \interior(\bTheta)$,
where $\bTheta = \mR \times \mR^+$.
Then conditions \eqref{nonfullid} and \eqref{smoothness4} of Theorem \ref{nonfull} are satisfied and hence,
for all $\epsilon > 0$,
there exist $\kappa(\epsilon) > 0$ and $C > 0$ such that $\widehat\btheta$ exists, is unique, and
\be
\nonumber
\mbP\left(\widehat\btheta\; \in\; \mB(\bthetas, \epsilon)\right)
\gte 1 - 4\, \exp\left(- \dfrac{\kappa(\epsilon)^2\; C\; K}{\norm{\mA}_\infty^4}\right),
\ee
provided $|\mA_k| \geq 3$ ($k = 1, \dots, K$) and $K \geq 2$.
\end{corollary}

\s

\com {\em Transitive edge terms versus triangle terms.}
\label{transitive.edge.triangle}
Transitive edge terms are not the same as triangle terms.
As pointed out in Section \ref{sec:curved},
there are many models of transitivity,
some of which are well-posed while others are ill-posed \citep[e.g., the triangle model with edge and triangle terms,][]{Jo99,Ha03p,Sc09b,BaBrSl11,ChDi11}.
As the curved exponential-family random graphs in Section \ref{sec:curved},
canonical exponential-family random graphs with within-neighborhood edge and transitive edge terms constrain the dependence among edges induced by transitivity to neighborhoods.
In addition,
the model ensures that the added value of additional triangles decreases provided $\theta_2 > 0$.
In fact, 
for each pair of nodes,
the value added by the first triangle to the log odds of the conditional probability of an edge is $\theta_2 > 0$,
whereas the added value of additional triangles is $0$.
Thus,
the model has similar properties as the curved exponential-family random graphs discussed in Section \ref{sec:curved}.

\section{Concentration results: $M$-estimators}
\label{supplement.mestimators}

The concentration and consistency results for maximum likelihood estimators in Section \ref{sec:curved} are special cases of more general results for $M$-estimators.
To demonstrate,
we introduce a natural class of $M$-estimators in Appendix \ref{assumptions},
which includes both likelihood- and moment-based estimators,
and present concentration results in Appendix \ref{mainresults} along with an application to misspecified models with omitted covariate terms.
These results cover both correct and incorrect model specifications,
as the example with omitted covariate terms demonstrates.

\subsection{$M$-estimators}
\label{assumptions}

A natural class of $M$-estimators,
which includes both\linebreak 
likelihood- and moment-based estimators,
can be constructed as follows.

Let $b: \mbX \mapsto \mR^{\rrr}$ be a vector of statistics.
These statistics might be
\bi
\item the sufficient statistics of the data-generating exponential family;
\item the sufficient statistics of an exponential family other than the data-generating exponential family,
which implies that the model is misspecified;
\item statistics motivated by computational considerations.
\ei
An example is presented in Corollary \ref{consistency.curved2},
where we consider a misspecified exponential family with omitted covariate terms.

A natural extension of the class of likelihood-based estimating functions in Section \ref{mle} is given by estimating functions of the form
\be
\label{mm}
g(\btheta;\, b(\bx))
\= \norm{b(\bx) - \bbeta(\btheta)}_2,
& \btheta \in \bTheta,
\ee
which are approximations of
\beno
g(\btheta;\, \mbE\, b(\bX))
\= \norm{\mbE\, b(\bX) - \bbeta(\btheta)}_2,
& \btheta \in \bTheta,
\ee
provided $\mbE\, b(\bX)$ exists.
Here,
$\mbE\, b(\bX)$ is the expectation of $b(\bX)$ under the data-generating exponential-family distribution.
The vector $\bbeta: \bTheta \mapsto \mR^{\rrr}$ is a vector-valued function of $\btheta \in \bTheta$ defined by $\bbeta(\btheta) = \mbE_{\btheta}\, b(\bX)$,
where $\mbE_{\btheta}\, b(\bX)$ is the expectation of $b(\bX)$ under a distribution parameterized by $\btheta \in \bTheta$,
provided $\mbE_{\btheta}\, b(\bX)$ exists.
The distribution may not belong to the data-generating exponential family or any other exponential family,
i.e.,
the model may be misspecified.

Estimating functions of the form \eqref{mm} cover both likelihood- and moment-based estimating functions:
\bi
\item Likelihood-based estimating functions: The choice $b(\bx) = s(\bx)$ gives\linebreak 
$g(\btheta;\, s(\bx)) = \norm{s(\bx) - \mbE_{\btheta}\, s(\bX)}_2$,
which is based on the gradient of the loglikelihood function with respect to the natural parameter vector of the data-generating exponential family 
and is therefore based on moments of the sufficient statistics of the data-generating exponential family.
\item Moment-based estimating functions: The choice $b(\bx) \neq s(\bx)$ gives\linebreak 
$g(\btheta;\, b(\bx)) = \norm{b(\bx) - \mbE_{\btheta}\, b(\bX)}_2$,
which is based on moments of statistics other than the sufficient statistics of the data-generating exponential family.
\ei

$M$-estimators based on estimating functions of the form \eqref{mm} are defined by
\beno
\bthetah
\= \left\{\btheta \in \bTheta:\;\; g(\btheta;\; b(\bx))\;\; =\;\; \inf\limits_{\dot\btheta\in\bTheta}\; g(\dot\btheta;\; b(\bx))\right\},
\ee
which are estimators of
\beno
\bthetad
\= \left\{\btheta \in \bTheta:\;\; g(\btheta;\; \mbE\, b(\bX))\;\; =\;\; \inf\limits_{\dot\btheta\in\bTheta}\; g(\dot\btheta;\; \mbE\, b(\bX))\right\}.
\ee
If,
given an observation $\bx$ of a random graph $\bX$,
the estimating function\linebreak 
$g(\btheta;\, b(\bx))$ is close to $g(\btheta;\, \mbE\, b(\bX))$,
then we would expect the minimizers $\bthetah$ and $\bthetad$ of $g(\btheta; b(\bx))$ and $g(\btheta; \mbE\, b(\bX))$ to be close under suitable conditions.
To show that $\bthetah$ is close to $\bthetad$ with high probability,
we make the following assumptions.
In the following,
$\BBB$ denotes the convex hull of the set $\{b(\bx): \bx \in \mbX\}$. 
\begin{enumerate}
\item[[C.1]\hspace{-.15cm}] The expectation $\mbE\, b(\bX) \in \rint(\BBB)$ exists and there exists $A > 0$ such that,
for all $(\bx_1, \bx_2) \in \mX\times\mX$,
\beno
\left\norm{b(\bx_1) - b(\bx_2)\right}_\infty
\lte A\, d(\bx_1, \bx_2)\, \norm{\mA}_\infty.
\ee
\item[[C.2]\hspace{-.15cm}] The parameter space $\bTheta$ is an open subset of $\mR^{\qqq}$,
the expectation $\bbeta(\btheta) = \mbE_{\btheta}\, b(\bX)$ exists for all $\btheta \in \bTheta$,
and there exists a unique element $\btheta_0 \in \bTheta$ such that $\bbeta(\btheta_0) = \mbE\, b(\bX) \in \rint(\BBB)$.
\item[[C.3]\hspace{-.15cm}] For all $\epsilon > 0$ small enough so that $\mathscr{B}(\bthetad, \epsilon) \subseteq \bTheta$,
there exists $\delta(\epsilon) > 0$ such that,
for all $\btheta \in \bTheta \setminus \mathscr{B}(\bthetad, \epsilon)$,
\beno
\norm{\bbeta(\btheta_0) - \bbeta(\btheta)}_2
&\geq& \delta(\epsilon)\; \dsum_{k=1}^K {|\mA_k| \choose 2}^\alpha
& \mbox{for some}
& 0 \leq \alpha \leq 1.
\ee
\end{enumerate}

Conditions [C.1]---[C.3] are verified in Corollary \ref{consistency.curved2}.
Conditions [C.1] and [C.3] resemble the conditions of Theorem \ref{nonfull} and are verified in Corollaries \ref{corollary.curved} and \ref{corollary.transitive1} in the special case $b(\bx) = s(\bx)$.
Condition [C.2] is a moment-matching condition:
It assumes that the family of distributions parameterized by $\btheta \in \bTheta$,
which may not be the data-generating exponential family,
is able to match the first moment $\mbE\, b(\bX)$ of $b(\bX)$ under the data-generating exponential-family distribution.
The unique parameter vector $\btheta_0 \in \bTheta$ that matches the moment $\mbE\, b(\bX)$ is the unique minimizer of $g(\btheta;\, \mbE\, b(\bX))$,
i.e.,
$\btheta_0 = \argmin_{\btheta\in\bTheta}\, g(\btheta;\, \mbE\, b(\bX))$.

\subsection{Concentration results}
\label{mainresults}

We establish concentration results for the class of $M$-estimators introduced in Section \ref{assumptions}.

To do so,
we need to show that the probability mass of statistic vector $b(\bX)$ concentrates around its expectation $\mbE\, b(\bX)$ under the data-generating exponential-family distribution.

\begin{proposition}
\label{concentration.b}
Consider an exponential family with countable support $\mbX$ and local dependence.
Let $b: \mbX \mapsto \mR^{\rrr}$ and assume that condition [C.1] is satisfied.
Then there exists $C > 0$ such that,
for all deviations of the form $t = \delta\, \sum_{k=1}^K {|\mA_k| \choose 2}^\alpha$ with $\delta > 0$ and $0 \leq \alpha \leq 1$,
\beno
\mbP\left(\norm{b(\bX) - \mbE\, b(\bX)}_2\, \geq\, t\right)
\lte 2\, \exp\left(- \dfrac{\delta^2\; C\; K}{\rrr\, \norm{\mA}_\infty^{4\, (2 - \alpha)}} + \log \rrr\right).
\ee
\end{proposition}

\s

Proposition \ref{concentration.b} paves the way for concentration results for the class of $M$-estimators introduced in Section \ref{assumptions}.
The following concentration result shows that $M$-estimators $\bthetah = \argmin_{\btheta\in\bTheta}\, g(\btheta;\, b(\bX))$ based on estimating functions of the form \eqref{mm} are close to $\bthetad = \argmin_{\btheta\in\bTheta}\, g(\btheta;\, \mbE\, b(\bX))$ with high probability as long as the neighborhoods and $\rrr$ are small relative to the number of neighborhoods $K$.

\begin{theorem}
\label{theorem:parameters.correct}
Consider an exponential-family random graph with countable support $\mbX$ and local dependence.
Assume that conditions [C.1]---[C.3] are satisfied.
Then,
for all $\epsilon > 0$ small enough so that $\mB(\bthetad, \epsilon) \subseteq \interior(\bTheta)$,
there exist $\kappa(\epsilon) > 0$ and $C > 0$ such that 
\be
\nonumber
\mbP\left(\widehat\btheta\; \in\; \mB(\bthetad, \epsilon)\right)
\;\geq\; 1 - 2\, \exp\left(- \dfrac{\kappa(\epsilon)^2\, C\, K}{\rrr\, \norm{\mA}_\infty^{4\, (2 - \alpha)}} + \log \rrr\right).
\ee
\end{theorem}

\s

Conditions [C.1]---[C.3] of Theorem \ref{theorem:parameters.correct} are verified by Corollary \ref{consistency.curved2}.
It is worth noting that $\bthetah$ may not be unique,
but an argument along the lines of Section \ref{sec:curved} shows that,
with high probability, 
the minimizers of estimating function \eqref{mm} do not give rise to global minima that are separated by large distances,
under the stated assumptions.

\s

{\em Example: misspecified exponential family with omitted covariate term.}
To demonstrate,
we consider an extension of the exponential-family random graph with within-neighborhood edge and transitive edge terms,
as described in Appendix \ref{sec:canonical}.
The extended model includes an additional same-attribute edge term with natural parameters $\eta_{k,3}(\btheta) = \theta_3$ and sufficient statistics $s_{k,3}(\bx_k)$,
\beno
s_{k,3}(\bx_k)
&=& \dsum_{i\in\mA_k\, <\, j\in\mA_k} x_{i,j}\, \one(c_i = c_j),
\ee
where $\one(c_i = c_j)$ is an indicator function,
which is $1$ if $c_i = c_j$ and is $0$ otherwise,
and $c_i \in \{C_1, \dots, C_H\}$ ($C_h \in \mR$,\, $h = 1, \dots, H$,\, $H \geq 2$) is a categorical attribute of node $i \in \mA_k$ ($k = 1, \dots, K$).
Same-attribute edge terms are popular in applications and capture excesses in the expected number of edges among nodes with the same attribute \citep[e.g.,][]{HuGoHa08}:
e.g.,
students may prefer to befriend other students of the same region of origin.

Suppose that researchers are unaware that the attributes $c_i$ of nodes $i$ are important predictors of edges and estimate the edge and transitive edge parameters $\theta_1$ and $\theta_2$ based on the misspecified exponential family with natural parameter vector $\btheta = (\theta_1, \theta_2)$ and sufficient statistic vector $b(\bx) = (s_1(\bx), s_2(\bx))$,
where $s_1(\bx) = \sum_{k=1}^K s_{k,1}(\bx_k)$ and $s_2(\bx) = \sum_{k=1}^K s_{k,2}(\bx_k)$ are defined in Appendix \ref{sec:canonical}.
In other words,
suppose that statistical inference is based on the estimating function
\beno
\label{mm.example}
g(\btheta;\; b(\bx))
\= \norm{b(\bx) - \mbE_{\btheta}\, b(\bX)}_2,
& \btheta \in \bTheta_0,
\ee
where $\bThetad = \mR \times \mR^+$ is the parameter space of the misspecified exponential family,
which is a two-dimensional subspace of the three-dimensional parameter space $\bThetas = \mR \times \mR^+ \times \mR$ of the data-generating exponential family.
The following concentration result shows that the $M$-estimator $\bthetah = \argmin_{\btheta\in\bTheta_0}\, g(\btheta;\; b(\bX))$ is close to 
$\bthetad = \argmin_{\btheta\in\bTheta_0}\, g(\btheta;\; \mbE\, b(\bX))$ with high probability as long as the neighborhoods are small relative to the number of neighborhoods $K$.

\s

\begin{corollary}
\label{consistency.curved2}
Consider an exponential-family random graph with within-neighborhood edge, transitive edge, and same-attribute edge terms.
Suppose that an observation of the random graph is generated by $\bthetas \in \bThetas$.
Then all conditions of Theorem \ref{theorem:parameters.correct} are satisfied and hence,
for all $\epsilon > 0$ small enough so that $\mB(\bthetad, \epsilon) \subseteq \bThetad$,
there exist $\kappa(\epsilon) > 0$ and $C > 0$ such that $\widehat\btheta$ exists, is unique, and 
\be
\nonumber
\mbP\left(\widehat\btheta\; \in\; \mB(\bthetad, \epsilon)\right)
&\geq& 1 - 4\, \exp\left(- \dfrac{\kappa(\epsilon)^2\, C\, K}{\norm{\mA}_\infty^4}\right),
\ee
provided $|\mA_k| \geq 3$ ($k = 1, \dots, K$) and $K \geq 2$.
\end{corollary}

\s

It is worth noting that $\bthetad \in \mR \times \mR^+$ and $\bthetas \in \mR \times \mR^+ \times \mR$ do not have the same dimension and cannot be identical,
but the distribution parameterized by $\bthetad$ is as close as possible to the data-generating distribution parameterized by $\bthetas$ in terms of Kullback-Leibler divergence:
by construction,
$\bthetad$ minimizes 
$\norm{\mbE\, b(\bX) - \mbE_{\btheta}\, b(\bX)}_2$\linebreak
$= \norm{\mbE\, \nabla_{\btheta} \log p_{\btheta}(\bX)}_2$ and hence maximizes $\mbE\, \log p_{\btheta}(\bX)$,
where\linebreak
$p_{\btheta}(\bx) \propto \exp(\langle\btheta,\, b(\bx)\rangle)$ denotes the density of $\bx \in \mbX$ under the misspecified exponential-family distribution with natural parameter vector $\btheta = (\theta_1, \theta_2)$ and sufficient statistic vector $b(\bx) = (s_1(\bx), s_2(\bx))$.
Therefore,
$\bthetad$ satisfies
\be
\label{kl}
\bthetad
= \argmax\limits_{\btheta\, \in\, \bThetad} \mbE\, \log p_{\btheta}(\bX) - \mbE\, \log p_{\bthetas}(\bX)
= \argmin\limits_{\btheta\, \in\, \bThetad} KL(\mbP_{\bthetas}; \mbP_{\btheta}),
\ee
where $KL(\mbP_{\bthetas}; \mbP_{\btheta})$ is the Kullback-Leibler divergence from $\mbP_{\bthetas}$ to $\mbP_{\btheta}$ and the expectations $\mbE\, \log p_{\btheta}(\bX)$ and $\mbE\, \log p_{\bthetas}(\bX)$ are with respect to $\mbP_{\bthetas}$.
Owing to the dependence of within-neighborhood edge variables and sufficient statistics,
it is not straightforward to bound $\norm{\bthetas - \bthetad}_2$,
but---by the properties of maximum likelihood estimation---we are assured that the distribution parameterized by $\bthetad$ is as close as possible to the distribution parameterized by the data-generating parameter vector $\bthetas$ in terms of Kullback-Leibler divergence,
as shown by \eqref{kl}.

\section{Proofs: Concentration results for random graphs with dependent edges}
\label{sec:proofs.concentration}

We prove the main concentration results of Sections \ref{sec:probability} and \ref{sec:model.estimation} and Appendix \ref{supplement.mestimators},
Propositions \ref{proposition.concentration},
\ref{existence.eta.mle},
and
\ref{concentration.b}.

\pproof \ref{proposition.concentration}.
By assumption,
$\mbE\, f(\bX)$ exists.
We are interested in deviations of the form $|f(\bX) - \mbE\, f(\bX)| \geq t$, 
where $t > 0$.
Choose any $t > 0$ and let $\mathscr{X} = \{\bx \in \mX:\, |f(\bx) - \mbE\, f(\bX)| \geq t\}$.
Since within-neighborhood edges do not depend on between-neighborhood edges,
\beno
\mbP(\bX \in \mathscr{X},\, \bY \in \mY)
\hide{
\= \mbP(g(\bX) \in \mathscr{X}) 
}
\= \mbP(\bX \in \mathscr{X}).
\ee
In the following,
we denote by $\mbP$ a probability measure on $(\mX, \mS)$ with densities of the form \eqref{local.ergm},
where $\mS$ is the power set of the countable set $\mX$.
Keep in mind that $\bX = (\bX_k)_{k=1}^K$ denotes the sequence of within-neighborhood edge variables,
where $\bX_k = (X_{i,j})_{i\in\mA_k\, <\, j\in\mA_k}$.
In an abuse of notation,
we denote the elements of the sequence of edge variables $\bX$ by $X_1, \dots, X_w$ with sample spaces $\mX_1, \dots, \mX_w$,
respectively,
where $w = \sum_{k=1}^K {|\mA_k| \choose 2}$ is the number of within-neighborhood edge variables.
Let $\bX_{i:j} = (X_i, \dots, X_j)$ be a subsequence of edge variables with sample space $\mX_{i:j}$,
where $i \leq j$.
By applying Theorem 1.1 of \citet[][]{KoRa08} to $\norm{f}_{\lip}$-Lipschitz functions $f: \mX \mapsto \mbR$ defined on the countable set $\mX$,
\beno
\mbP(|f(\bX) - \mbE\, f(\bX)| \;\geq\; t) 
\lte 2\, \exp\left(- \dfrac{t^2}{2\, w\, \norm{\Phi}_\infty^2\, \norm{f}_{\lip}^2}\right),
\ee
where $\Phi$ is the $w \times w$-upper triangular matrix with entries
\beno
\phi_{i,j} \=
\begin{cases}
\varphi_{i,j} & \mbox{if } i < j\\
1 & \mbox{if } i = j\\
0 & \mbox{if } i > j
\end{cases}
\ee
and
\beno
\norm{\Phi}_\infty
\= \max\limits_{1 \leq i \leq w} \left|1 + \dsum_{j=i+1}^w \varphi_{i,j}\right|.
\ee
The coefficients $\varphi_{i,j}$ are known as mixing coefficients and are defined by 
\beno
\label{mixing.coefficients}
\varphi_{i,j} 
\equiv \sup\limits_{\substack{\bx_{1:i-1} \in \mbX_{1:i-1}\\(\uu,\, \vv) \in \mX_i\times\mX_i}} \varphi_{i,j}(\bx_{1:i-1}, \uu, \vv)
\= \sup\limits_{\substack{\bx_{1:i-1} \in \mbX_{1:i-1}\\(\uu,\, \vv) \in \mX_i\times\mX_i}} \norm{\pi_{\uu} - \pi_{\vv}}_\tv,
\ee
where $\norm{\pi_{\uu} - \pi_{\vv}}_\tv$ is the total variation distance between the distributions $\pi_{\uu}$ and $\pi_{\vv}$ 
given by
\beno
\pi_{\uu}
\equiv \pi(\bx_{j:w} \mid \bx_{1:i-1}, \uu) 
= \mbP(\bX_{j:w} = \bx_{j:w} \mid \bX_{1:i-1} = \bx_{1:i-1}, X_i = \uu)
\ee
and
\beno
\pi_{\vv}
\equiv \pi(\bx_{j:w} \mid \bx_{1:i-1}, \vv) 
= \mbP(\bX_{j:w} = \bx_{j:w} \mid \bX_{1:i-1} = \bx_{1:i-1}, X_i = \vv).
\ee
Since the support 
of $\pi_{\uu}$ and $\pi_{\vv}$ is countable,
\beno
\norm{\pi_{\uu} - \pi_{\vv}}_\tv
\hide{
\= \dfrac12\, \norm{\pi_{\uu} - \pi_{\vv}}_1
}
\= \dfrac12 \dsum\limits_{\bx_{j:w} \in \mbX_{j:w}} |\pi(\bx_{j:w} \mid \bx_{1:i-1}, \uu) - \pi(\bx_{j:w} \mid \bx_{1:i-1}, \vv)|.
\ee
An upper bound on $\norm{\Phi}_\infty$ can be obtained by bounding the mixing coefficients $\varphi_{i,j}$ as follows.
Consider any pair of edge variables $X_i$ and $X_j$.
If $X_i$ and $X_j$ involve nodes in more than one neighborhood,
the mixing coefficient $\varphi_{i,j}$ vanishes by the local dependence induced by exponential families of the form \eqref{local.ergm}.
If the pair of nodes corresponding to $X_i$ and the pair of nodes corresponding to $X_j$ belong to the same neighborhood,
the mixing coefficient $\varphi_{i,j}$ can be bounded as follows:
\beno
\varphi_{i,j}(\bx_{1:i-1}, \uu, \vv) 
\hide{
\;=\; \norm{\pi - \pi^\prime}_\tv
}
\;=\; \dfrac12 \dsum\limits_{\bx_{j:w} \in \mbX_{j:w}} |\pi(\bx_{j:w} \mid \bx_{1:i-1}, \uu) - \pi(\bx_{j:w} \mid \bx_{1:i-1}, \vv)|\s
\hide{
\\
\lte \dfrac12 \dsum\limits_{x_{j:m} \in \mbX_{j:m}} \left(|\pi(x_{j:m} \mid \bx_{1:i-1}, \uu)| + |\pi(x_{j:m} \mid \bx_{1:i-1}, \vv)|\right)\s
}
\\
\hspace{0cm}\leq\; \dfrac12 \dsum\limits_{\bx_{j:w} \in \mbX_{j:w}} \pi(\bx_{j:w} \mid \bx_{1:i-1}, \uu) + \dfrac12 \dsum\limits_{\bx_{j:w} \in \mbX_{j:w}} \pi(\bx_{j:w} \mid \bx_{1:i-1}, \vv)
\;=\; 1,
\ee
because $\pi_{\uu}$ and $\pi_{\vv}$ are conditional probability mass functions with countable support $\mbX_{j:w}$.
We note that the upper bound is not sharp,
but it has the advantage that it covers a wide range of dependencies within neighborhoods.
Thus,
\beno
\norm{\Phi}_\infty
\= \max\limits_{1 \leq i \leq w} \left|1 + \dsum_{j=i+1}^w \varphi_{i,j}\right|
\lte \dis{\norm{\mA}_\infty \choose 2},
\ee
because each edge variable $X_i$ can depend on at most ${\norm{\mA}_\infty \choose 2} - 1$ other edge variables corresponding to pairs of nodes belonging to the same pair of neighborhoods.
Therefore,
there exists $C > 0$ such that,
for all $t > 0$,
\beno
\mbP(|f(\bX) - \mbE\, f(\bX)|\, \geq\, t) 
\hide{
\lte 2\, \exp\left(- \dfrac{t^2}{2\, m\, \norm{\Phi}_\infty^2\, \norm{f}_{\lip}^2}\right)\s
\\
}
\,\leq\, 2\, \exp\left(- \dfrac{t^2}{C\, \sum_{k=1}^K {|\mA_k| \choose 2}\, \norm{\mAs}_\infty^4\, \norm{f}_{\lip}^2}\right),
\ee
where $\norm{\mA}_\infty \geq 1$ because all neighborhoods $\mA_k$ are non-empty and $\norm{f}_{\lip} > 0$ by assumption.

\pproof \ref{existence.eta.mle}.
By assumption,
the neighborhood-dependent natural parameters $\eta_{k,i}(\btheta)$ are of the form $\eta_{k,i}(\btheta) = \eta_i(\btheta)$ ($i = 1, \dots, \ppp_k$, $k = 1, \dots, K$).
Therefore,
the exponential family can be reduced to an exponential family with natural parameter vector
\beno
\bta(\btheta)
\= \left(\eta_1(\btheta), \dots, \eta_{\diminf}(\btheta)\right)
\ee
and sufficient statistic vector
\beno
s(\bx)
\= \left(s_1(\bx), \dots, s_{\diminf}(\bx)\right),
\ee
where $\diminf = \max_{1 \leq k \leq K} \ppp_k$.
Here, 
the sufficient statistics $s_i(\bx)$ are sums of within-neighborhood sufficient statistics $s_{k,i}(\bx_k)$ ($k = 1, \dots, K$):
\beno
s_i(\bx)
\= \dsum_{k=1}^K s_{k,i}(\bx_k),
& i = 1, \dots, \diminf.
\ee
Observe that $\widehat{\bmu(\bta(\bthetas))} = s(\bX)$ and that $\bmu(\bta(\bthetas)) = \mbE_{\bta(\bthetas)}\, s(\bX) = \mbE\, s(\bX) \in \rint(\mM)$ exists, 
because $\bta(\bthetas) \in \interior(\fullspace)$ \citep[][Theorem 2.2, pp.\ 34--35]{Br86}.
We bound the probability of deviations of $s(\bX)$ from $\mbE\, s(\bX)$ in terms of the $\ell_\infty$- and $\ell_2$-norm below.

\s

{\bf $\ell_\infty$-norm.}
For all $\delta > 0$,
\beno
&& \mbP\left(\norm{\widehat{\bmu(\bta(\bthetas))} - \bmu(\bta(\bthetas))}_\infty\; \geq\; \delta\; \dsum_{k=1}^K {|\mA_k| \choose 2}^\alpha\right)\s
\\
\= \mbP\left(\norm{s(\bX) - \mbE\, s(\bX)}_\infty\; \geq\; \delta\; \dsum_{k=1}^K {|\mA_k| \choose 2}^\alpha\right)\s
\\
\lte \mbP\left(\dis\bigcup_{i=1}^{\diminf} \left(\left|s_i(\bX) - \mbE\, s_i(\bX)\right| \;\geq\; \delta\; \dsum_{k=1}^K {|\mA_k| \choose 2}^\alpha\right)\right).
\ee
By condition \eqref{smoothness6} of Proposition \ref{existence.eta.mle},
there exists $A > 0$ such that the Lipschitz coefficient of $s_i: \mX \mapsto \mR$ with respect to the Hamming metric $d: \mX \times \mX \mapsto \mbR_0^+$ is bounded above by $\norm{s_i}_{\lip} \leq A\; \norm{\mA}_\infty$ ($i = 1, \dots, \diminf$).
\hide{
because,
for all $(\bx_1, \bx_2) \in \mX\times\mX$,
\beno
|f_i(\bx_1) - f_i(\bx_2)|
\;\leq\; \left\norm{(\mJ|_{\btheta=\bthetas})^\top\, (s(\bx_1) - s(\bx_2))\right}_\infty
\;\leq\; A_4\, d(\bx_1, \bx_2)\, \norm{\mA}_\infty
\ee
by condition [C.2].
}
Thus,
by a union bound over the $\diminf$ components of $s(\bX)$ and by applying Proposition \ref{proposition.concentration} to deviations of size $t = \delta\; \sum_{k=1}^K {|\mA_k| \choose 2}^\alpha$,
we obtain,
for all $\delta > 0$,
\beno
&& \mbP\left(\dis\bigcup_{i=1}^{\diminf} \left(\left|s_i(\bX) - \mbE\, s_i(\bX)\right| \;\geq\; \delta\; \dsum_{k=1}^K {|\mA_k| \choose 2}^\alpha\right)\right)\s
\hide{
\;\leq\; \dsum_{i=1}^{\diminf} \mbP\left(\left|g_i(\bthetas;\, \bX)\right| \;\geq\; \epsilon\; \dsum_{k=1}^K {|\mA_k| \choose 2}^\alpha\right)\s
}
\\
\lte 2\, \exp\left(- \dfrac{\delta^2\, \left(\sum_{k=1}^K {|\mA_k| \choose 2}^\alpha\right)^2}{B\, \sum_{k=1}^K {|\mA_k| \choose 2}\, \norm{\mAs}_\infty^4\, \norm{\mAs}_\infty^2} + \log \diminf\right),
\ee
where $B > 0$.
Since all neighborhoods are of the same order of magnitude,
there exist $C_1 > 0$ and $C_2 > 0$ such that $C_1\, \norm{\mA}_\infty \leq |\mA_k| \leq C_2\, \norm{\mA}_\infty$.
Thus,
there exists $C_3 > 0$ such that the term $(\sum_{k=1}^K {|\mA_k| \choose 2}^\alpha)^2$ in the numerator of the exponent can be bounded below by
\beno
\left(\dsum_{k=1}^K \dis{|\mA_k| \choose 2}^{\alpha} \right)^2
\gte C_3\; K^2\; \norm{\mA}_\infty^{4\, \alpha},
\ee
and there exists $C_4 > 0$ such that the term $\sum_{k=1}^K {|\mA_k| \choose 2}$ in the denominator of the exponent can be bounded above by
\beno
\dsum_{k=1}^K \dis{|\mA_k| \choose 2}
\lte C_4\, K\; \norm{\mA}_\infty^2.
\ee
As a result,
there exists $C > 0$ such that,
for all $\delta > 0$,
\beno
&& \mbP\left(\norm{\widehat{\bmu(\bta(\bthetas))} - \bmu(\bta(\bthetas))}_\infty\; \geq\; \delta\; \dsum_{k=1}^K {|\mA_k| \choose 2}^\alpha\right)\s
\\
\lte 2 \exp\left(- \dfrac{\delta^2\; C\; K}{\norm{\mA}_\infty^{4\, (2 - \alpha)}} + \log \diminf\right).
\ee

\s

{\bf $\ell_2$-norm.}
The same argument used above shows that,
for all $\delta > 0$,
\hide{
\beno
&& \mbP\left(\norm{\widehat{\bmu(\bta(\bthetas))} - \bmu(\bta(\bthetas))}_1 \geq \delta\; \dsum_{k=1}^K {|\mA_k| \choose 2}^\alpha\right)\s
\\
\lte
\hide{
\mbP\left(\diminf\; \left\norm{s(\bX) - \mbE\, s(\bX)\right}_\infty \;\geq\; \epsilon\right)\s
\\
\=
}
\mbP\left(\norm{\widehat{\bmu(\bta(\bthetas))} - \bmu(\bta(\bthetas))}_\infty\, \geq\, \dfrac{\delta\; \sum_{k=1}^K {|\mA_k| \choose 2}^\alpha}{\diminf}\right)\s
\\
\lte 2\, \exp\left(- \dfrac{\delta^2\; C\; K}{\diminf^2\, \norm{\mA}_\infty^{4\, (2 - \alpha)}} + \log \diminf\right)
\ee
and
}
\beno
&& \mbP\left(\norm{\widehat{\bmu(\bta(\bthetas))} - \bmu(\bta(\bthetas))}_2\, \geq\, \delta\; \dsum_{k=1}^K {|\mA_k| \choose 2}^\alpha\right)\s
\\
\lte
\hide{
\mbP\left(\sqrt{\diminf}\; \left\norm{s(\bX) - \mbE\, s(\bX)\right}_\infty \;\geq\; \epsilon\right)\s
\\
\=
}
\mbP\left(\norm{\widehat{\bmu(\bta(\bthetas))} - \bmu(\bta(\bthetas))}_\infty\, \geq\, \dfrac{\delta\; \sum_{k=1}^K {|\mA_k| \choose 2}^\alpha}{\sqrt{\diminf}}\right)\s
\\
\lte 2\, \exp\left(- \dfrac{\delta^2\; C\; K}{\diminf\, \norm{\mA}_\infty^{4\, (2 - \alpha)}} + \log \diminf\right).
\ee

\pproof \ref{concentration.b}.
Condition [C.1] implies that $\mbE\, b(\bX)$ exists.
Thus,
for all $\delta > 0$,
\beno
&& \mbP\left(\left\norm{b(\bX) - \mbE\, b(\bX)\right}_2\; \geq\; \delta\; \dsum_{k=1}^K {|\mA_k| \choose 2}^\alpha\right)\s
\\
\lte \mbP\left(\left\norm{b(\bX) - \mbE\, b(\bX)\right}_\infty \geq \dfrac{\delta\; \sum_{k=1}^K {|\mA_k| \choose 2}^\alpha}{\sqrt{\rrr}}\right)\s
\\
\lte \mbP\left(\dis\bigcup_{i=1}^{\rrr} \left(\left|b_i(\bX) - \mbE\, b_i(\bX)\right| \;\geq\; \dfrac{\delta\; \sum_{k=1}^K {|\mA_k| \choose 2}^\alpha}{\sqrt{\rrr}}\right)\right).
\ee
By condition [C.1],
there exists $A > 0$ such that the Lipschitz coefficient of $b_i: \mX \mapsto \mR$ with respect to the Hamming metric $d: \mX \times \mX \mapsto \mbR_0^+$ is bounded above by $\norm{b_i}_{\lip} \leq A\; \norm{\mA}_\infty$ ($i = 1, \dots, \rrr$).
Thus,
by a union bound over the $\rrr$ components of $b(\bX)$ and by applying Proposition \ref{proposition.concentration} to deviations of size $t = \delta\; \sum_{k=1}^K {|\mA_k| \choose 2}^\alpha\, /\, \sqrt{\rrr}$,
we obtain,
for all $\delta > 0$,
\beno
&& \mbP\left(\dis\bigcup_{i=1}^{\rrr} \left(\left|b_i(\bX) - \mbE\, b_i(\bX)\right| \;\geq\; \dfrac{\delta\; \sum_{k=1}^K {|\mA_k| \choose 2}^\alpha}{\sqrt{\rrr}}\right)\right)
\\
\lte 2\, \exp\left(- \dfrac{\delta^2\, \left(\sum_{k=1}^K {|\mA_k| \choose 2}^\alpha\right)^2}{B\, \rrr\, \sum_{k=1}^K {|\mA_k| \choose 2}\, \norm{\mAs}_\infty^4\, \norm{\mA}_\infty^2} + \log \rrr\right),
\ee
where $B > 0$.
We showed in the proof of Proposition \ref{existence.eta.mle} that there exist $C_1 > 0$ and $C_2 > 0$ such that 
$(\sum_{k=1}^K {|\mA_k| \choose 2}^\alpha)^2 \geq C_1\; K^2\; \norm{\mA}_\infty^{4\, \alpha}$
and
$\sum_{k=1}^K {|\mA_k| \choose 2} \leq C_2\, K\; \norm{\mA}_\infty^2$.
As a consequence,
there exists $C > 0$ such that,
for all $\delta > 0$,
\beno
&& \mbP\left(\norm{b(\bX) - \mbE\, b(\bX)}_2\; \geq\ \delta\; \dsum_{k=1}^K {|\mA_k| \choose 2}^\alpha\right)\s
\\
\lte 2\, \exp\left(- \dfrac{\delta^2\; C\; K}{\rrr\, \norm{\mA}_\infty^{4\, (2 - \alpha)}} + \log \rrr\right).
\ee

\s

\section{Proofs: Concentration results for maximum likelihood estimators}
\label{mle.proofs}

We prove the main concentration results of Section \ref{mle},
Theorem \ref{nonfull} along with Corollaries \ref{corollary.curved} and \ref{corollary.transitive1}.
Auxiliary lemmas are proved in Appendix \ref{m.id}.

\ttproof {\bf \ref{nonfull}}.
By assumption,
$\bthetas \in \interior(\bTheta)$.
Observe that $\bthetas \in \interior(\bTheta)$ implies $\bmu(\bta(\bthetas)) \in \rint(\mM)$,
because the maps $\bta: \interior(\bTheta) \mapsto \interior(\fullspace)$ and $\bmu: \interior(\fullspace) \mapsto \rint(\mM)$ are one-to-one: 
the map $\bta: \interior(\bTheta) \mapsto \interior(\fullspace)$ is one-to-one by assumption while the map $\bmu: \interior(\fullspace) \mapsto \rint(\mM)$ is one-to-one by classic exponential-family theory \citep[][Theorem 3.6, p.\ 74]{Br86}.

We are interested in estimators of the form
\beno
\widehat\btheta
\= \left\{\btheta \in \bTheta:\;\; \norm{s(\bx) - \bmu(\bta(\btheta))}_2\;\; =\;\; \inf\limits_{\dot\btheta\in\bTheta}\; \norm{s(\bx) - \bmu(\bta(\dot\btheta))}_2\right\},
\ee
where
\beno
\bTheta
&\subseteq& \{\btheta \in \mR^{\qqq}:\; \psi(\bta(\btheta)) < \infty\}.
\ee
The following proof covers both full and non-full exponential families,
including curved exponential families.
We first discuss the existence of $\bthetah$ in full and non-full, curved exponential families and then bound the probability of event $\bthetah \in \mB(\bthetas, \epsilon)$.

\s

{\bf Existence of $\bthetah$: full and non-full, curved exponential families.}
For any given $\omega > 0$,
let
\beno
\mM(\omega)
&=& \left\{\bmu' \in \mM:\; \norm{\bmu' - \bmu(\bta(\bthetas))}_2\; <\; \omega\; \dsum_{k=1}^K {|\mA_k| \choose 2}^\alpha\right\}
\ee
be the subset of $\bmu' \in \mM$ close to the data-generating mean-value parameter vector $\bmu(\bta(\bthetas)) \in \rint(\mM)$ in the sense that $\norm{\bmu' - \bmu(\bta(\bthetas))}_2\, <\, \omega\, \sum_{k=1}^K {|\mA_k| \choose 2}^\alpha$.
Choose $\omega > 0$ small enough so that $\mM(2\, \omega) \subseteq \rint(\mM)$ and let
\beno
\mbG
\= \left\{\bx \in \mbX:\;\, \norm{s(\bx) - \bmu(\bta(\bthetas))}_2\; <\; \omega\; \dsum_{k=1}^K {|\mA_k| \choose 2}^\alpha\right\}
\ee
be the subset of $\bx  \in \mbX$ such that $s(\bx) \in \mM(\omega) \subseteq \rint(\mM)$.
In the following,
we show that in the event $\mbG$ the set $\bthetah$ is non-empty,
i.e.,
$\bthetah$ exists.
To see that,
let $\mM(\bTheta)$ be the set of mean-value parameter vectors induced by $\bTheta$,
i.e.,
\beno
\mM(\bTheta)
&=& \left\{\bmu' \in \rint(\mM): \mbox{ there exists } \btheta \in \bTheta \mbox{ such that } \bmu(\bta(\btheta)) = \bmu'\right\}.
\ee
If the exponential family is full,
then $\mM(\bTheta) = \rint(\mM)$,
otherwise $\mM(\bTheta) \subset \rint(\mM)$,
because non-full exponential families exclude some natural parameter vectors along with the corresponding mean-value parameter vectors.
\hide{
Since the set $\bTheta$ is non-empty and contains at least one element---the element $\bthetas \in \interior(\bTheta)$---and the maps $\bta: \interior(\bTheta) \mapsto \interior(\fullspace)$ and $\bmu: \interior(\fullspace) \mapsto \rint(\mM)$ are one-to-one \citep[][Theorem 3.6, p.\ 74]{Br86},
the set $\mM(\bTheta) \subseteq \rint(\mM)$ is non-empty and contains at least one element---the element $\bmu(\bta(\bthetas)) \in \mM(\bTheta) \subseteq \rint(\mM)$.
}
To show that the set $\bthetah$ is non-empty in the event $\mbG$,
note that,
for all $\bx \in \mbG$ and hence all $s(\bx) \in \mM(\omega) \subseteq \rint(\mM)$,
there exists at least one element $\bmu(\bta(\bthetah)) \in \mM(\bTheta) \subseteq \rint(\mM)$ such that
\beno
\norm{s(\bx) - \bmu(\bta(\bthetah))}_2
\lte \norm{s(\bx) - \bmu(\bta(\bthetas))}_2,
\ee
because the data-generating mean-value parameter vector $\bmu(\bta(\bthetas)) \in \mM(\bTheta) \subseteq \rint(\mM)$ is known to be an element of $\mM(\bTheta) \subseteq \rint(\mM)$ and any minimizer of $\norm{s(\bx) - \bmu(\bta(\btheta))}_2$ cannot be farther from $s(\bx) \in \mM(\omega) \subseteq \rint(\mM)$ than the data-generating mean-value parameter vector $\bmu(\bta(\bthetas)) \in \mM(\bTheta) \subseteq \rint(\mM)$.
In addition,
for each element $\bmu(\bta(\bthetah)) \in \mM(\bTheta) \subseteq \rint(\mM)$,
we have 
\beno
\norm{\bmu(\bta(\bthetah)) - \bmu(\bta(\bthetas))}_2
\lte \norm{s(\bx) - \bmu(\bta(\bthetah))}_2 + \norm{s(\bx) - \bmu(\bta(\bthetas))}_2\s\s
\\
\lte 2\, \norm{s(\bx) - \bmu(\bta(\bthetas))}_2,
\ee
which implies that each element $\bmu(\bta(\bthetah)) \in \mM(\bTheta) \subseteq \rint(\mM)$ is contained in $\mM(2\, \omega) \subseteq \rint(\mM)$;
note that $\omega > 0$ was chosen small enough so that $\mM(2\, \omega) \subseteq \rint(\mM)$.
Therefore,
for all $\bx \in \mbG$,
the set $\bthetah$ contains at least one element.
If there exists a unique element $\btheta \in \bTheta$ such that $\bmu(\bta(\btheta)) = s(\bx) \in \mM(\omega) \subseteq \rint(\mM)$,
then the set $\bthetah$ contains one and only one element.
In particular,
in full exponential families with $\bta(\btheta) = \btheta$,
for all $\bx \in \mbG$,
there exists a unique element $\btheta \in \bTheta$ such that $\bmu(\bta(\btheta)) = s(\bx) \in \mM(\omega) \subseteq \rint(\mM)$ \citep[][Theorem 3.6, p.\ 74]{Br86}.
Therefore,
in full exponential families $\bthetah$ exists and is unique in the event $\mbG$,
whereas in non-full exponential families $\bthetah$ exists but may not be unique in the event $\mbG$.

\s

{\bf Bounding the probability of event $\bthetah \in \mB(\bthetas, \epsilon)$: full and non-full, curved exponential families.}
By the identifiability conditions of Theorem \ref{nonfull},
for all $\epsilon > 0$ small enough so that $\mathscr{B}(\bthetas, \epsilon) \subseteq \interior(\bTheta)$,
there exists $\gamma(\epsilon) > 0$ such that,
for all $\btheta \in \bTheta \setminus \mB(\bthetas, \epsilon)$,
we have $\bta(\btheta) \in \fullspace \setminus \mB(\bta(\bthetas),\, \gamma(\epsilon))$. 
In addition,
there exists $\delta(\epsilon) > 0$ such that,
for all $\bta(\btheta) \in \fullspace \setminus \mB(\bta(\bthetas), \gamma(\epsilon))$,
\beno
\norm{\bmu(\bta(\btheta)) - \bmu(\bta(\bthetas))}_2
\gte \delta(\epsilon)\; \dsum_{k=1}^K {|\mA_k| \choose 2}^\alpha
\;\mbox{ for some }\; 0 \leq \alpha \leq 1.
\ee
Therefore,
\beno
\norm{\bmu(\bta(\btheta)) - \bmu(\bta(\bthetas))}_2
&<& \delta(\epsilon)\; \dsum_{k=1}^K {|\mA_k| \choose 2}^\alpha
\ee
implies $\bta(\btheta) \in \mB(\bta(\bthetas),\, \gamma(\epsilon))$,
which in turn implies $\btheta \in \mB(\bthetas, \epsilon)$.

To bound the probability of event $\bthetah \in \mB(\bthetas, \epsilon)$,
we exploit the fact that,
for all $\bx \in \mbG$ and hence all $s(\bx) \in \mM(\omega) \subseteq \rint(\mM)$,
the set $\bthetah$ is non-empty and,
for each element of the set $\bthetah$,
\beno
\label{triangle0}
\norm{\bmu(\bta(\bthetah)) - \bmu(\bta(\bthetas))}_2
\lte \norm{s(\bx) - \bmu(\bta(\bthetah))}_2 + \norm{s(\bx) - \bmu(\bta(\bthetas))}_2\s\s
\\
\lte 2\, \norm{s(\bx) - \bmu(\bta(\bthetas))}_2,
\ee
which implies that $\bmu(\bta(\bthetah)) \in \mM(2\, \omega) \subseteq \rint(\mM)$.
Thus,
the probability of event 
\beno
\norm{\bmu(\bta(\bthetah)) - \bmu(\bta(\bthetas))}_2\; <\; \delta(\epsilon)\; \dsum_{k=1}^K {|\mA_k| \choose 2}^\alpha &\cap& \mbG
\ee
can be bounded from below by bounding the probability of event 
\beno
2\, \norm{s(\bx) - \bmu(\bta(\bthetas))}_2\; <\; \delta(\epsilon)\; \dsum_{k=1}^K {|\mA_k| \choose 2}^\alpha &\cap& \mbG,
\ee
which, by definition of event $\mbG$, is equivalent to bounding the probability of event
\beno
\norm{s(\bx) - \bmu(\bta(\bthetas))}_2\; <\; \min\left(\dfrac{\delta(\epsilon)}{2},\; \omega\right)\; \dsum_{k=1}^K {|\mA_k| \choose 2}^\alpha.
\ee
These results show that the probability of event $\bthetah \in \mB(\bthetas, \epsilon)$ can be bounded from below as follows:
\beno
&& \mbP\left(\bthetah\; \in\; \mB(\bthetas, \epsilon)\; \subseteq\; \mbox{int}(\bTheta)\right)
\;\geq\; \mbP\left(\bthetah\; \in\; \mB(\bthetas, \epsilon)\; \subseteq\; \mbox{int}(\bTheta)\; \cap\; \mbG\right)\s
\\
\gte \mbP\left(\norm{\bmu(\bta(\bthetah)) - \bmu(\bta(\bthetas))}_2\; <\; \delta(\epsilon)\; \dsum_{k=1}^K {|\mA_k| \choose 2}^\alpha\; \cap\; \mbG\right)\s
\\
\= 
\hide{
\mbP\left(\norm{s(\bX) - \bmu(\bta(\bthetas))}_2\; <\; \dfrac{\delta(\epsilon)}{2}\; \dsum_{k=1}^K {|\mA_k| \choose 2}^\alpha\; \cap\; \mbG\right)\s
\\
\= 
}
\mbP\left(\norm{s(\bX) - \bmu(\bta(\bthetas))}_2\; <\; \min\left(\dfrac{\delta(\epsilon)}{2},\; \omega\right)\; \dsum_{k=1}^K {|\mA_k| \choose 2}^\alpha\right).
\ee
To bound the probability of the event on the right-hand side of the inequality above,
note that the smoothness condition \eqref{smoothness4} of Theorem \ref{nonfull} ensures that the smoothness condition \eqref{smoothness6} of Proposition \ref{existence.eta.mle} is satisfied.
Therefore,
Proposition \ref{existence.eta.mle} can be invoked to show that there exists $C > 0$ such that
\beno
&& \mbP\left(\norm{s(\bX) - \bmu(\bta(\bthetas))}_2\; \geq\; \min\left(\dfrac{\delta(\epsilon)}{2},\; \omega\right)\; \dsum_{k=1}^K {|\mA_k| \choose 2}^\alpha\right)\s
\\
\lte 2\, \exp\left(- \dfrac{\kappa(\epsilon)^2\, C\, K}{\diminf\, \norm{\mA}_\infty^{4\, (2 - \alpha)}} + \log \diminf\right),
\ee
where $\kappa(\epsilon) = \min(\delta(\epsilon)\, /\, 2,\; \omega) > 0$.
Collecting results shows that 
\beno
\mbP\left(\bthetah\; \in\; \mB(\bthetas, \epsilon)\; \subseteq\; \mbox{int}(\bTheta)\right)
\gte 1 - 2\, \exp\left(- \dfrac{\kappa(\epsilon)^2\, C\, K}{\diminf\, \norm{\mA}_\infty^{4\, (2 - \alpha)}} + \log \diminf\right).
\ee

\ccproof \ref{corollary.curved}.
Corollary \ref{corollary.curved} follows from Theorem \ref{nonfull},
because all conditions of Theorem \ref{nonfull} are satisfied.
First,
by Lemma \ref{lemma:gw_grad_ident} in Appendix \ref{appendix.mle.id.curved},
the map $\bta: \interior(\bTheta) \mapsto \interior(\fullspace)$ is one-to-one and,
for all $\epsilon > 0$ small enough so that $\mathscr{B}(\bthetas, \epsilon) \subseteq \interior(\bTheta)$,
there exists $\gamma(\epsilon) > 0$ such that,
for all $\btheta \in \bTheta \setminus \mB(\bthetas, \epsilon)$,
we have $\bta(\btheta) \in \fullspace \setminus \mB(\bta(\bthetas),\, \gamma(\epsilon))$.
Second,
condition \eqref{nonfullid} of Theorem \ref{nonfull} is satisfied with $\alpha = 3/4$ by Lemma \ref{lemma:theorem_2_ident} in Appendix \ref{appendix.mle.id.curved} provided $|\mA_k| \geq 4$ ($k = 1, \dots, K$).
Last,
but not least,
condition \eqref{smoothness4} of Theorem \ref{nonfull} is satisfied,
because changing an edge cannot change the number of within-neighborhood edges by more than $1$ and the number of within-neighborhood connected pairs of nodes with $i$ shared partners by more than $2\, (\norm{\mA}_\infty - 2) + 1$ ($i = 1, \dots, |\mA_k|-2$, $k = 1, \dots, K$).
Thus,
by Theorem \ref{nonfull},
for all $\epsilon > 0$ small enough so that $\mB(\bthetas, \epsilon) \subseteq \interior(\bTheta)$,
there exist $\kappa(\epsilon) > 0$ and $C > 0$ such that 
\be
\nonumber
\mbP\left(\widehat\btheta\; \in\; \mB(\bthetas, \epsilon)\; \subseteq\; \interior(\bTheta)\right)
\gte 1 - 2\, \exp\left(- \dfrac{\kappa(\epsilon)^2\, C\, K}{\diminf\, \norm{\mA}_\infty^{4\, (2 - \alpha)}} + \log \diminf\right)\s
\\
\gte 1 - 2\, \exp\left(- \dfrac{\kappa(\epsilon)^2\, C\, K}{\norm{\mA}_\infty^6} + \log\, \norm{\mA}_\infty\right),
\ee
where we used the fact that $\diminf = \max_{1 \leq k \leq K} \ppp_k = \norm{\mA}_\infty - 1$ and $\alpha = 3/4$.

\ccproof \ref{corollary.transitive1}.
Corollary \ref{corollary.transitive1} follows from Theorem \ref{nonfull},
because conditions \eqref{nonfullid} and \eqref{smoothness4} of Theorem \ref{nonfull} are satisfied.
Condition \eqref{nonfullid} is satisfied with $\alpha = 1$ by Lemma \ref{edge.transitive1} in Appendix \ref{appendix.mle.id.curved} provided $|\mA_k| \geq 3$ ($k = 1, \dots, K$).
Condition \eqref{smoothness4} is satisfied, 
because changing an edge cannot change the number of within-neighborhood edges by more than $1$ and the number of within-neighborhood transitive edges by more than $2\, (\norm{\mA}_\infty - 2) + 1$.
Thus,
Theorem \ref{nonfull} can be invoked to conclude that,
for all $\epsilon > 0$ small enough so that $\mB(\bthetas, \epsilon) \subseteq \interior(\bTheta)$,
there exist $\kappa(\epsilon) > 0$,
$C_1 > 0$, 
and $C_2 > 0$ such that
\be
\nonumber
\mbP\left(\bthetah \in \mB(\bthetas, \epsilon) \subseteq \interior(\bTheta)\right)
\geq\; 1 - 2\, \exp\left(- \dfrac{\kappa(\epsilon)^2\; C_1\; K}{\qqq\; \norm{\mA}_\infty^{4\, (2 - \alpha)}} + \log \qqq\right)\s
\\
\hspace{4.5cm}=\; 1 - 4\, \exp\left(- \dfrac{\kappa(\epsilon)^2\; C_2\; K}{\norm{\mA}_\infty^4}\right),
\ee
where we used the fact that $\qqq = 2$ and $\alpha = 1$.

\s

\section{Proofs: Concentration results for $M$-estimators}
\label{sec:proofs.consistency}

We prove the main concentration result of Appendix \ref{supplement.mestimators},
Theorem \ref{theorem:parameters.correct} along with Corollary \ref{consistency.curved2}.

\ttproof \ref{theorem:parameters.correct}.
Theorem \ref{theorem:parameters.correct} can be proved along the same lines as Theorem \ref{nonfull}.

By condition [C.1],
the expectation $\mbE\, b(\bX)$ of $b(\bX)$ under the data-generating exponential-family distribution exists.
Let $\BBB$ be the convex hull of the set $\{b(\bx): \bx \in \mbX\}$ and,
given any $\omega > 0$,
let
\beno
\BBB(\omega)
\= \left\{\bbeta' \in \BBB:\; \norm{\bbeta' - \mbE\, b(\bX)}_2\, <\, \omega\; \dsum_{k=1}^K {|\mA_k| \choose 2}^\alpha\right\}
\ee
be the subset of $\bbeta' \in \BBB$ close to the expectation $\mbE\, b(\bX)$ of $b(\bX)$ under the data-generating exponential-family distribution in the sense that $\norm{\bbeta' - \mbE\, b(\bX)}_2\, <\, \omega\, \sum_{k=1}^K {|\mA_k| \choose 2}^\alpha$.
Choose $\omega > 0$ small enough so that $\BBB(2\, \omega) \subseteq \rint(\BBB)$ and let
\beno
\mbG
&=& \left\{\bx \in \mbX:\; \norm{b(\bx) - \mbE\, b(\bX)}_2 < \omega\, \dsum_{k=1}^K {|\mA_k| \choose 2}^\alpha\right\}
\ee
be the subset of $\bx  \in \mbX$ such that $b(\bx) \in \BBB(\omega) \subseteq \rint(\BBB)$.

We are interested in estimators of the form
\beno
\widehat\btheta
\= \left\{\btheta \in \bTheta:\;\; \norm{b(\bx) - \bbeta(\btheta)}_2\;\; =\;\; \inf\limits_{\dot\btheta\in\bTheta}\; \norm{b(\bx) - \bbeta(\dot\btheta)}_2\right\},
\ee
where $\bTheta$ is an open subset of $\mR^{\qqq}$ and the expectation $\bbeta(\btheta) = \mbE_{\btheta}\, b(\bX)$ exists for all $\btheta \in \bTheta$ by condition [C.2].

In the following,
we do not assume that the family of distributions parameterized by $\btheta \in \bTheta$ is an exponential family and we do not assume that the set of expectations $\{\mbE_{\btheta}\, b(\bX),\, \btheta \in \bTheta\}$ covers the whole set $\rint(\BBB)$.
In other words,
the family may not be able to match all possible $b(\bx) \in \rint(\BBB)$ in the sense that there exists $\btheta \in \bTheta$ such that $\mbE_{\btheta}\, b(\bX) = b(\bx) \in \rint(\BBB)$.
The critical assumption exploited by the following proof is that the family can match the expectation $\mbE\, b(\bX)$ of $b(\bX)$ under the data-generating exponential-family distribution in the sense that there exists one and only one element $\btheta_0 \in \bTheta$ such that $\mbE_{\btheta_0}\, b(\bX) = \mbE\, b(\bX) \in \rint(\BBB)$,
along with the fact that $b(\bX)$ falls with high probability into a small ball centered at $\mbE_{\btheta_0}\, b(\bX) = \mbE\, b(\bX) \in \rint(\BBB)$ under suitable conditions.
We first discuss the existence of $\bthetah$ and then bound the probability of event $\bthetah \in \mB(\btheta_0, \epsilon)$.

\s

{\bf Existence of $\bthetah$.}
We show that in the event $\mathscr{G}(\omega)$ the set $\bthetah$ is non-empty,
i.e.,
$\bthetah$ exists.
To show that the set $\bthetah$ is non-empty in the event $\mathscr{G}(\omega)$,
note that,
for all $\bx \in \mathscr{G}(\omega)$ and hence all $b(\bx) \in \BBB(\omega) \subseteq \rint(\BBB)$,
there exists at least one element $\btheta \in \bTheta$ such that
\beno
\norm{b(\bx) - \bbeta(\btheta)}_2
\lte \norm{b(\bx) - \bbeta(\btheta_0)}_2
\= \norm{b(\bx) - \mbE\, b(\bX)}_2,
\ee
because,
by condition [C.2],
there exists $\btheta_0 \in \bTheta$ such that $\bbeta(\btheta_0) = \mbE\, b(\bX) \in \BBB(\omega) \subseteq \rint(\BBB)$ and any minimizer of $\norm{b(\bx) - \bbeta(\btheta)}_2$ cannot be farther from $b(\bx) \in \BBB(\omega) \subseteq \rint(\BBB)$ than $\bbeta(\btheta_0) = \mbE\, b(\bX) \in \BBB(\omega) \subseteq \rint(\BBB)$.
Therefore,
the set $\bthetah$ is non-empty in the event $\mathscr{G}(\omega)$.
In addition,
for each element of the set $\bthetah$,
we have
\beno
\norm{\bbeta(\bthetah) - \mbE\, b(\bX)}_2
\lte \norm{b(\bx) - \bbeta(\bthetah)}_2 + \norm{b(\bx) - \mbE\, b(\bX)}_2\s\s
\\
\lte 2\, \norm{b(\bx) - \mbE\, b(\bX)}_2,
\ee
which implies that $\bbeta(\bthetah) \in \BBB(2\, \omega) \subseteq \rint(\BBB)$;
note that $\omega > 0$ was chosen small enough so that $\BBB(2\, \omega) \subseteq \rint(\BBB)$.
Therefore,
for all $\bx \in \mathscr{G}(\omega)$ and hence all $b(\bx) \in \BBB(\omega) \subseteq \rint(\BBB)$,
the set $\bthetah$ contains at least one element.
If the set $\bthetah$ contains more than one element,
then all elements of the set $\bthetah$ map to expectations $\bbeta(\bthetah) = \mbE_{\bthetah}\; b(\bX) \in \BBB(2\, \omega) \subseteq \rint(\BBB)$ that have the same $\ell_2$-distance from $b(\bx) \in \BBB(\omega) \subseteq \rint(\BBB)$ by construction of estimating function $\norm{b(\bx) - \bbeta(\btheta)}_2$.

\s

{\bf Bounding the probability of event $\bthetah \in \mB(\btheta_0, \epsilon)$.}
To bound the probability of event $\bthetah \in \mB(\btheta_0, \epsilon)$,
note that,
by condition [C.2], 
the expectation $\bbeta(\btheta) = \mbE_{\btheta}\, b(\bX)$ exists for all $\btheta \in \bTheta$ and,
by condition [C.3],
for all $\epsilon > 0$ small enough so that $\mathscr{B}(\btheta_0, \epsilon) \subseteq \bTheta$,
there exists $\delta(\epsilon) > 0$ such that,
for all $\btheta \in \bTheta \setminus \mathscr{B}(\btheta_0, \epsilon)$,
\beno
\norm{\bbeta(\btheta) - \bbeta(\btheta_0)}_2
\= \norm{\bbeta(\btheta) - \mbE\, b(\bX)}_2
\gte \delta(\epsilon)\; \dsum_{k=1}^K {|\mA_k| \choose 2}^\alpha
\ee
for some $0 \leq \alpha \leq 1$.
Therefore,
\beno
\norm{\bbeta(\btheta) - \bbeta(\btheta_0)}_2
\= \norm{\bbeta(\btheta) - \mbE\, b(\bX)}_2
&<& \delta(\epsilon)\; \dsum_{k=1}^K {|\mA_k| \choose 2}^\alpha
\ee
implies $\btheta \in \mB(\btheta_0, \epsilon)$.

To bound the probability of event $\bthetah \in \mB(\btheta_0, \epsilon)$,
we exploit the fact that,
for all $\bx \in \mathscr{G}(\omega)$ and hence all $b(\bx) \in \BBB(\omega) \subseteq \rint(\BBB)$,
the set $\bthetah$ is non-empty and,
for each element of the set $\bthetah$,
\beno
\label{triangle0}
\norm{\bbeta(\bthetah) - \mbE\, b(\bX)}_2
\lte \norm{b(\bx) - \bbeta(\bthetah)}_2 + \norm{b(\bx) - \mbE\, b(\bX)}_2\s\s
\\
\lte 2\, \norm{b(\bx) - \mbE\, b(\bX)}_2,
\ee
which implies that $\bbeta(\bthetah) \in \BBB(2\, \omega) \subseteq \rint(\BBB)$,
as pointed out above.
Thus,
the probability of event
\beno
\norm{\bbeta(\bthetah) - \mbE\, b(\bX)}_2\; <\; \delta(\epsilon)\; \dsum_{k=1}^K {|\mA_k| \choose 2}^\alpha &\cap& \mathscr{G}(\omega)
\ee
can be bounded from below by bounding the probability of event
\beno
2\, \norm{b(\bx) - \mbE\, b(\bX)}_2\; <\; \delta(\epsilon)\; \dsum_{k=1}^K {|\mA_k| \choose 2}^\alpha &\cap& \mathscr{G}(\omega),
\ee
which, by definition of event $\mathscr{G}(\omega)$, is equivalent to bounding the probability of event
\beno
\norm{b(\bx) - \mbE\, b(\bX)}_2\; <\; \min\left(\dfrac{\delta(\epsilon)}{2},\, \omega\right)\; \dsum_{k=1}^K {|\mA_k| \choose 2}^\alpha.
\ee
These results show that the probability of event $\bthetah \in \mB(\btheta_0, \epsilon)$ can be bounded from below as follows:
\beno
&& \mbP\left(\bthetah\; \in\; \mB(\btheta_0, \epsilon)\; \subseteq\; \bTheta\right)
\;\geq\; \mbP\left(\bthetah\; \in\; \mB(\btheta_0, \epsilon)\; \subseteq\; \bTheta\; \cap\; \mathscr{G}(\omega)\right)\s
\\
\gte \mbP\left(\norm{\bbeta(\bthetah) - \mbE\, b(\bX)}_2\; <\; \delta(\epsilon)\; \dsum_{k=1}^K {|\mA_k| \choose 2}^\alpha\; \cap\; \mathscr{G}(\omega)\right)\s
\\
\gte \mbP\left(\norm{b(\bX) - \mbE\, b(\bX)}_2\; <\; \min\left(\dfrac{\delta(\epsilon)}{2},\, \omega\right)\; \dsum_{k=1}^K {|\mA_k| \choose 2}^\alpha\right).
\ee
To bound the probability of the event on the right-hand side of the inequality above,
Proposition \ref{concentration.b} can be invoked to show that there exists $C > 0$ such that
\beno
&& \mbP\left(\norm{b(\bX) - \mbE\, b(\bX)}_2\; \geq\; \min\left(\dfrac{\delta(\epsilon)}{2},\, \omega\right)\; \dsum_{k=1}^K {|\mA_k| \choose 2}^\alpha\right)\s
\\
\lte 2\, \exp\left(- \dfrac{\kappa(\epsilon)^2\, C\, K}{\rrr\, \norm{\mA}_\infty^{4\, (2 - \alpha)}} + \log \rrr\right),
\ee
where $\kappa(\epsilon) = \min(\delta(\epsilon)\, /\, 2,\; \omega) > 0$.
Collecting results shows that
\beno
\mbP\left(\bthetah\; \in\; \mB(\btheta_0, \epsilon)\; \subseteq\; \bTheta\right)
\gte 1 - 2\, \exp\left(- \dfrac{\kappa(\epsilon)^2\, C\, K}{\rrr\, \norm{\mA}_\infty^{4\, (2 - \alpha)}} + \log \rrr\right).
\ee

\hide{

By condition [C.1],
$\mbE\, b(\bX)$ exists.
Let
\beno
\mbG
&=& \left\{\bx \in \mbX:\; b(\bx) \in \rint(\BBB),\; \norm{b(\bx) - \mbE\, b(\bX)}_2 < \omega\, \dsum_{k=1}^K {|\mA_k| \choose 2}^\alpha\right\}
\ee
be the subset of $\bx  \in \mbX$ such that $b(\bx) \in \BBB(\omega) \subseteq \rint(\BBB)$,
where $\omega > 0$ is identical to the constant $\omega$ in condition [C.2].

By condition [C.3],
for all $\epsilon > 0$ small enough so that $\mB(\bthetad, \epsilon) \subseteq \interior(\bTheta)$,
there exists $\deltaepsilon > 0$ such that,
for all $\btheta \in \bTheta \setminus \mB(\bthetad, \epsilon)$,
\beno
g(\btheta;\, \mbE\, b(\bX)) - g(\bthetad;\, \mbE\, b(\bX))
\;=\; \norm{\mbE\, b(\bX) - \bbeta(\btheta)}_2
\;\geq\; \delta(\epsilon)\, \dsum_{k=1}^K {|\mA_k| \choose 2}^\alpha,
\ee
where we used the fact that $g(\bthetad;\, \mbE\, b(\bX)) = 0$.
Observe that $g(\bthetad;\, \mbE\, b(\bX)) = 0$ follows from the assumption that there exists $\btheta \in \interior(\bTheta)$ such that $\bbeta(\btheta) = \mbE\, b(\bX)$ and hence there exists $\btheta \in \interior(\bTheta)$ such that $g(\btheta;\, \mbE\, b(\bX)) = 0$,
which implies that the minimizer $\bthetad$ of $g(\btheta;\, \mbE\, b(\bX))$ must satisfy $g(\bthetad;\, \mbE\, b(\bX)) = 0$;
note that the minimizer $\bthetad$ of $g(\btheta;\, \mbE\, b(\bX))$ exists and is unique by condition [C.2].
Thus,
$g(\btheta;\, \mbE\, b(\bX)) - g(\bthetad;\, \mbE\, b(\bX))$ $< \delta(\epsilon)\, \sum_{k=1}^K {|\mA_k| \choose 2}^\alpha$ implies $\btheta\, \in\, \mB(\bthetad, \epsilon)\, \subseteq\, \interior(\bTheta)$.
We want to show that,
when $\bthetah$ exists,
then
\beno
g(\bthetah;\, \mbE\, b(\bX)) - g(\bthetad;\, \mbE\, b(\bX))
&<& \delta(\epsilon)\, \dsum_{k=1}^K {|\mA_k| \choose 2}^\alpha 
\ee
implies $\bthetah\, \in\, \mB(\bthetad, \epsilon)\, \subseteq\, \interior(\bTheta)$,
so that the probability of event $\bthetah\, \in\, \mB(\bthetad, \epsilon)\, \subseteq\, \interior(\bTheta)$ can be bounded by bounding the probability of event\linebreak 
$g(\bthetah;\, \mbE\, b(\bX)) - g(\bthetad;\, \mbE\, b(\bX)) < \delta(\epsilon)\, \sum_{k=1}^K {|\mA_k| \choose 2}^\alpha$.
To do so,
note that,
by condition [C.2],
the $M$-estimator $\bthetah$ exists and is unique in the event $\mbG$.
We can therefore bound the probability of event $\bthetah\, \in\, \mB(\bthetad, \epsilon)\, \subseteq\, \interior(\bTheta)$ as follows:
\beno
\mbP(\bthetah\; \in\; \mB(\bthetad, \epsilon)\; \subseteq\; \interior(\bTheta))
\gte \mbP(\bthetah\; \in\; \mB(\bthetad, \epsilon)\; \subseteq\; \interior(\bTheta)\; \cap\; \mbG).
\ee
To bound the probability of event $\bthetah\, \in\, \mB(\bthetad, \epsilon)\, \subseteq\, \interior(\bTheta)\, \cap\, \mbG$,
observe that,
by condition [C.3],
the event
\beno
g(\bthetah;\, \mbE\, b(\bX)) - g(\bthetad;\, \mbE\, b(\bX))\; <\; \delta(\epsilon)\, \dsum_{k=1}^K {|\mA_k| \choose 2}^\alpha
&\cap& \mbG
\ee
implies the event $\bthetah\, \in\, \mB(\bthetad, \epsilon)\; \cap\; \mbG$.
As a consequence,
we can bound the probability of event $\bthetah\, \in\, \mB(\bthetad, \epsilon)\, \subseteq\, \interior(\bTheta)\, \cap\, \mbG$ as follows: 
\beno
&& \mbP(\bthetah\, \in\, \mB(\bthetad, \epsilon)\, \subseteq\, \interior(\bTheta)\, \cap\, \mbG)\s
\\
\gte \mbP\left(g(\bthetah;\, \mbE\, b(\bX)) - g(\bthetad;\, \mbE\, b(\bX))\, <\, \delta(\epsilon)\; \dsum_{k=1}^K {|\mA_k| \choose 2}^\alpha\; \cap\; \mbG\right).
\ee
To bound the probability of the event on the right-hand side of the inequality above,
define
\beno
\mE
= \left\{\bx \in \mbX:\; g(\bthetah;\, \mbE\, b(\bX)) - g(\bthetad;\, \mbE\, b(\bX)) < \delta(\epsilon)\; \dsum_{k=1}^K {|\mA_k| \choose 2}^\alpha\right\}
\ee
and note that the probability of the event on the right-hand side of the inequality above is equal to $\mbP(\mE\, \cap\, \mbG)$.
We bound the probability of event $\mE\, \cap\, \mbG$ as follows:
\beno
\mbP(\mE\, \cap\, \mbG)
\gte 1 - \mbP(\comp{\mE}) - \mbP(\comp{\mbG})\s
\\
\gte 
1 - \mbP(\comp{\mE}\, \cap\, \mbG) - 2\; \mbP(\comp{\mbG}),
\ee
where $\comp{\mE}$ and $\comp{\mbG}$ indicate the complement of events $\mE$ and $\mbG$ in $\mbX$,
respectively.

\s

{\bf Bounding the probability of event } $\comp{\mE} \cap \mbG${\bf.}
The probability of event $\comp{\mE}\, \cap\, \mbG$ can be bounded by observing that,
for all $\bx \in \mbG$,
the $M$-estimator $\bthetah$ exists, 
is unique, 
and satisfies
\beno
g(\bthetah;\, b(\bx)) - g(\bthetad;\, b(\bx))
\lte 0,
\ee
which implies that $g(\bthetah;\, \mbE\, b(\bX)) - g(\bthetad;\, \mbE\, b(\bX))$ is bounded above by
\hide{
\beno
&& g(\bthetad;\, \mbE\, b(\bX)) + [g(\bthetad;\, b(\bx)) - g(\bthetad;\, \mbE\, b(\bX))]\s
\\
&-& g(\bthetah;\, \mbE\, b(\bX)) - [g(\bthetah;\, b(\bx)) - g(\bthetah;\, \mbE\, b(\bX))]
\lte 0.
\ee
Therefore,
for all $\bx \in \mbG$,
$g(\bthetad;\, \mbE\, b(\bX)) - g(\bthetah;\, \mbE\, b(\bX))$ is bounded above by
}
\beno
&& g(\bthetah;\, \mbE\, b(\bX)) - g(\bthetad;\, \mbE\, b(\bX))\s
\\
\hide{
\lte [g(\bthetah;\, b(\bx)) - g(\bthetah;\, \mbE\, b(\bX))] - [g(\bthetad;\, b(\bx)) - g(\bthetad;\, \mbE\, b(\bX))]\s
\\
\lte |[g(\bthetah;\, b(\bx)) - g(\bthetah;\, \mbE\, b(\bX))] - [g(\bthetad;\, b(\bx)) - g(\bthetad;\, \mbE\, b(\bX))]|\s
\\
\lte |[g(\bthetah;\, b(\bx)) - g(\bthetah;\, \mbE\, b(\bX))| + |g(\bthetad;\, b(\bx)) - g(\bthetad;\, \mbE\, b(\bX))|\s
\\
}
\lte |g(\bthetah;\, b(\bx)) - g(\bthetah;\, \mbE\, b(\bX))| + |g(\bthetad;\, b(\bx)) - g(\bthetad;\, \mbE\, b(\bX))|.
\ee
To bound the two terms on the right-hand side,
observe that,
for any $\btheta \in \interior(\bTheta)$,
$g(\btheta;\; b(\bx)) - g(\btheta;\; \mbE\, b(\bX))$ is given by
\beno
\label{basicequality}
g(\btheta;\; b(\bx)) - g(\btheta;\; \mbE\, b(\bX))
\,=\, \norm{b(\bx) - \mbE_{\btheta}\, b(\bX)}_2 - \norm{\mbE\, b(\bX) - \mbE_{\btheta}\, b(\bX)}_2.
\ee
By the triangle inequality,
the term $\norm{b(\bx) - \mbE_{\btheta}\, b(\bX)}_2$ is bounded above by
\beno
\norm{b(\bx) - \mbE_{\btheta}\, b(\bX)}_2
\lte \norm{b(\bx) - \mbE\, b(\bX)}_2\s
+ \norm{\mbE\, b(\bX) - \mbE_{\btheta}\, b(\bX)}_2,
\ee
which implies that $g(\btheta;\; b(\bx)) - g(\btheta;\; \mbE\, b(\bX))$ is bounded above by
\beno
\label{basicequality}
g(\btheta;\; b(\bx)) - g(\btheta;\; \mbE\, b(\bX))
\lte \norm{b(\bx) - \mbE\, b(\bX)}_2.
\ee
Along the same lines,
the term $\norm{\mbE\, b(\bX) - \mbE_{\btheta}\, b(\bX)}_2$ can be bounded above by
\beno
\norm{\mbE\, b(\bX) - \mbE_{\btheta}\, b(\bX)}_2
\lte \norm{b(\bx) - \mbE\, b(\bX)}_2 + \norm{b(\bx) - \mbE_{\btheta}\, b(\bX)}_2,
\ee
which implies that $g(\btheta;\; b(\bx)) - g(\btheta;\; \mbE\, b(\bX))$ is bounded below by
\beno
g(\btheta;\; b(\bx)) - g(\btheta;\; \mbE\, b(\bX))
\gte -\norm{b(\bx) - \mbE\, b(\bX)}_2.
\ee
Thus,
for all $\btheta \in \interior(\bTheta)$ and all $\bx \in \mbG$,
\beno
|g(\btheta;\, b(\bx)) - g(\btheta;\, \mbE\, b(\bX))|
\lte \norm{b(\bx) - \mbE\, b(\bX)}_2.
\ee
As a result,
for all $\btheta \in \interior(\btheta)$,
including any fixed $\bthetad \in \interior(\bTheta)$ and any fixed $\bthetah \in \interior(\btheta)$ based on fixed $\bx \in \mbG$,
\beno
g(\bthetah;\; \mbE\, b(\bX)) - g(\bthetad;\; \mbE\, b(\bX))
\hide{
\\
\lte |g(\bthetad;\, b(\bx)) - g(\bthetad;\, \mbE\, b(\bX))| + |g(\bthetah;\, b(\bx)) - g(\bthetah;\, \mbE\, b(\bX))|\s
\\
}
\lte 2\; \norm{b(\bx) - \mbE\, b(\bX)}_2
\ee
and hence 
\beno
\mbP\left(g(\bthetah;\, \mbE\, b(\bX)) - g(\bthetad;\, \mbE\, b(\bX))\; \geq\; \delta(\epsilon)\; \dsum_{k=1}^K {|\mA_k| \choose 2}^\alpha\; \cap\; \mbG\right)\s
\\
\leq\; \mbP\left(\norm{b(\bX) - \mbE\, b(\bX)}_2\; \geq\; \dfrac{\delta(\epsilon)}{2}\; \dsum_{k=1}^K {|\mA_k| \choose 2}^\alpha\; \cap\; \mbG\right)\s
\\
\leq\; \mbP\left(\norm{b(\bX) - \mbE\, b(\bX)}_2\; \geq\; \dfrac{\delta(\epsilon)}{2}\; \dsum_{k=1}^K {|\mA_k| \choose 2}^\alpha\right).
\ee
By assumption,
smoothness condition [C.1] is satisfied,
hence Proposition \ref{concentration.b} can be invoked.
Proposition \ref{concentration.b} shows that there exists $C_1 > 0$ such that
\beno
\label{compl.e.g}
\mbP\left(g(\bthetah;\, \mbE\, b(\bX)) - g(\bthetad;\, \mbE\, b(\bX)) \geq \delta(\epsilon)\, \dsum_{k=1}^K {|\mA_k| \choose 2}^\alpha\, \cap\, \mbG\right)\s
\\
\leq\; 2\, \exp\left(- \dfrac{\delta(\epsilon)^2\, C_1\, K}{\rrr\, \norm{\mA}_\infty^{4\, (2 - \alpha)}} + \log \rrr\right).
\ee

{\bf Bounding the probability of event} $\comp{\mbG}$.
The probability of event\linebreak 
$\comp{\mbG}$ can be bounded by Proposition \ref{concentration.b},
which shows that there exists $C_2 > 0$ such that the probability of event $\comp{\mbG}$ satisfies
\beno
\label{compl.g}
\mbP(\comp{\mbG})
\= \mbP\left(\norm{b(\bX) - \mbE\, b(\bX)}_2\; \geq\; \omega\, \dsum_{k=1}^K {|\mA_k| \choose 2}^\alpha\right)\s
\\
\lte 2\, \exp\left(- \dfrac{\omega^2\, C_2\, K}{\rrr\, \norm{\mA}_\infty^{4\, (2 - \alpha)}} + \log \rrr\right).
\ee

{\bf Conclusion.}
Collecting results,
we conclude that,
for all $\epsilon > 0$ small enough so that $\mB(\bthetad, \epsilon)\, \subseteq\, \interior(\bTheta)$,
there exist $\delta(\epsilon) > 0$ and $C_3 = \min(C_1, C_2) > 0$ such that
\beno
\mbP(\bthetah \in \mB(\bthetad, \epsilon) \subseteq \interior(\bTheta))
\,\geq\, 1 - 6\, \exp\left(- \dfrac{\min(\omega^2, \delta(\epsilon)^2)\, C_3\, K}{\rrr\, \norm{\mA}_\infty^{4\, (2 - \alpha)}} + \log \rrr\right).
\ee

}

\ccproof \ref{consistency.curved2}.
Corollary \ref{consistency.curved2} follows from Theorem \ref{theorem:parameters.correct},
because all conditions of Theorem \ref{theorem:parameters.correct} are satisfied.
Condition [C.1] is satisfied,
because $\mbE\, b(\bX) \in \rint(\BBB)$ follows from $\bthetas \in \interior(\bThetas)$ and because changing an edge cannot change the number of within-neighborhood edges by more than $1$ and the number of within-neighborhood transitive edges by more than $2\, (\norm{\mA}_\infty - 2) + 1$.
Conditions [C.2] and [C.3] follow from classic exponential-family theory \citep{Br86},
because $\btheta \in \bThetad$ is the natural parameter vector and $\bbeta(\btheta) \in \rint(\BBB)$ is the mean-value parameter vector of the misspecified exponential family.
Condition [C.2] is satisfied,
because 
$\bbeta(\btheta) = \mbE_{\btheta}\, b(\bX)$ exists for all $\btheta \in \bThetad$ by classic exponential-family theory \citep[][Theorem 2.2, pp.\ 34--35]{Br86}
and there exists a unique element $\bthetad \in \bThetad$ such that $\bbeta(\bthetad) = \mbE_{\bthetad}\, b(\bX) = \mbE\, b(\bX)$,
as the map $\bbeta: \interior(\bThetad) \mapsto \rint(\BBB)$ is one-to-one by classic exponential-family theory \citep[e.g., Theorem 3.6,][p.\ 74]{Br86}.
Condition [C.3] is satisfied with $\alpha = 1$ by Lemma \ref{edge.transitive.attribute} in Appendix \ref{appendix.m.id} provided $|\mA_k| \geq 3$ ($k = 1, \dots, K$).
Since all conditions of Theorem \ref{theorem:parameters.correct} are satisfied,
we can invoke Theorem \ref{theorem:parameters.correct} to conclude that,
for all $\epsilon > 0$ small enough so that $\mB(\bthetad, \epsilon) \subseteq \bThetad$,
there exist $\kappa(\epsilon) > 0$, $C_1 > 0$, and $C_2 > 0$ such that the $M$-estimator $\widehat\btheta$ exists, is unique, and is contained in $\mB(\bthetad, \epsilon)$ with probability
\be
\nonumber
\mbP\left(\widehat\btheta\; \in\; \mB(\bthetad, \epsilon)\; \subseteq\; \bThetad\right)
\gte 1 - 2 \exp\left(- \dfrac{\kappa(\epsilon)^2\, C_1\, K}{\rrr\, \norm{\mA}_\infty^{4\, (2 - \alpha)}} + \log \rrr\right)\s
\\
\= 1 - 4\, \exp\left(- \dfrac{\kappa(\epsilon)^2\, C_2\, K}{\norm{\mA}_\infty^4}\right),
\ee
where we used the fact that $\rrr = 2$ and $\alpha = 1$.
Observe that uniqueness follows from the fact that there is a unique $\btheta \in \bThetad$ such that $\mbE_{\btheta}\, b(\bX) = b(\bx) \in \rint(\BBB)$ for all possible $b(\bx) \in \rint(\BBB)$,
because $\btheta \in \bThetad$ is the natural parameter vector and $\bbeta(\btheta) \in \rint(\BBB)$ is the mean-value parameter vector of the misspecified exponential family and therefore the map $\bbeta: \interior(\bThetad) \mapsto \rint(\BBB)$ is one-to-one \citep[e.g., Theorem 3.6,][p.\ 74]{Br86}.
 
\section{Proofs: Concentration results for subgraph-to-graph estimators}
\label{appendix.incomplete}

We prove the main concentration result of Section \ref{extendability},
Theorem \ref{incomplete} along with Proposition \ref{proposition.extendability} and Corollary \ref{corollary.curved.incomplete}.

\pproof \ref{proposition.extendability}.
By the factorization properties implied by local dependence,
we have,
for all $\btheta\in\bTheta \subseteq \{\btheta \in \mR^{\qqq}:\; \psi_{\mL}(\bta(\btheta)) < \infty\}$,
all $\bx_{\mL} \in \mX_{\mL}$,
and all $\mK \subseteq \mL$,
\be
\label{factorization}
&& \dsum_{\bx_{\mL \setminus \mK} \in \mX_{\mL \setminus \mK}} p_{\bta(\btheta)}(\bx_{\mL}) \s
\\
\= \dsum_{\bx_{\mL \setminus \mK} \in \mX_{\mL \setminus \mK}} \exp\left(\langle\bta(\btheta),\, s(\bx_{\mL})\rangle - \psi_{\mL}(\bta(\btheta))\right)\, \nu_{\mL}(\bx_{\mL})\s
\\
\= \dsum_{\bx_{\mL \setminus \mK} \in \mX_{\mL \setminus \mK}}\; \dprod_{\mA \in \mL} \exp\left(\langle\bta(\btheta),\, s(\bx_{\mA})\rangle - \psi_{\mA}(\bta(\btheta))\right)\, \nu_{\mA}(\bx_{\mA})\s
\\
\= \dprod_{\mA \in \mK} \exp\left(\langle\bta(\btheta),\, s(\bx_{\mA})\rangle - \psi_{\mA}(\bta(\btheta))\right)\, \nu_{\mA}(\bx_{\mA})\s
\\
\= \exp\left(\langle\bta(\btheta),\, s(\bx_{\mK})\rangle - \psi_{\mK}(\bta(\btheta))\right)\, \nu_{\mK}(\bx_{\mK})\s
\\
\= p_{\bta(\btheta)}(\bx_{\mK}),
\ee
where we exploited the fact that, by local dependence, 
$\nu_{\mL}(\bx_{\mK})$ satisfies
\beno
\nu_{\mL}(\bx_{\mL})
\= \dprod_{\mA \in \mL} \nu_{\mA}(\bx_{\mA}),
\ee
while $\psi_{\mL}(\bta(\btheta))$ satisfies
\beno
\psi_{\mL}(\bta(\btheta))
\= \dsum_{\mA \in \mL} \psi_{\mA}(\bta(\btheta))
\ee
and
\beno
\psi_{\mA}(\bta(\btheta))
\= \dsum_{\bx_{\mA} \in \mX_{\mA}} \exp\left(\langle\bta(\btheta),\, s(\bx_{\mA})\rangle\right)\, \nu_{\mA}(\bx_{\mA}),
& \mA \in \mL.
\ee
We note that the natural parameter vector $\bta(\btheta)$ takes the form
\beno
\bta(\btheta) 
\= (\eta_1(\btheta), \dots, \eta_{\diminf}(\btheta)),
\ee
while the sufficient statistic vector takes the form
\beno
s(\bx_{\mK}) 
\= (s_1(\bx_{\mK}), \dots, s_{\diminf}(\bx_{\mK})),
\ee
where the sufficient statistics $s_i(\bx_{\mK})$ are based on the sufficient statistics $s_{\mA,i}(\bx_{\mA})$ of neighborhoods $\mA \in \mK$:
\beno
s_i(\bx_{\mK}) 
\= \dsum_{\mA \in \mK} s_{\mA,i}(\bx_{\mA}),
& i = 1, \dots, \diminf.
\ee
In addition,
note that
\beno
\nu_{\mK}(\bx_{\mK})
\= \dprod_{\mA \in \mK} \nu_{\mA}(\bx_{\mA})
\ee
and that
\beno
\psi_{\mK}(\bta(\btheta))
\= \dsum_{\mA \in \mK} \psi_{\mA}(\bta(\btheta))
\ee
is finite provided $\psi_{\mL}(\bta(\btheta)) = \sum_{\mA \in \mL} \psi_{\mA}(\bta(\btheta))$ is finite.

To prove that,
for all $\bx_{\mK} \in \mX_{\mK}$ and all $\by_{\mK} \in \mY_{\mK}$,
\beno
&& \mbP_{\bta(\btheta)}(\bX_{\mK} = \bx_{\mK},\; \bY_{\mK} = \by_{\mK},\; \bX_{\mL \setminus \mK} \in \mX_{\mL \setminus \mK},\; \bY_{\mL \setminus \mK} \in \mbY_{\mL \setminus \mK})\s
\\
\= \mbP_{\bta(\btheta)}(\bX_{\mK} = \bx_{\mK},\; \bY_{\mK} = \by_{\mK}),
\ee
observe that the independence of within- and between-neighborhood subgraphs implies that
\be
\label{basic}
& \mbP_{\bta(\btheta)}(\bX_{\mK} = \bx_{\mK},\; \bY_{\mK} = \by_{\mK},\; \bX_{\mL \setminus \mK} \in \mX_{\mL \setminus \mK},\; \bY_{\mL \setminus \mK} \in \mbY_{\mL \setminus \mK})\s
\\
&=\; \mbP_{\bta(\btheta)}(\bX_{\mK} = \bx_{\mK},\; \bX_{\mL \setminus \mK} \in \mX_{\mL \setminus \mK})\; \mbP(\bY_{\mK} = \by_{\mK},\; \bY_{\mL \setminus \mK} \in \mbY_{\mL \setminus \mK}).
\ee
The first term on the right-hand side of \eqref{basic} can be dealt with by using \eqref{factorization},
which implies that
\beno
\mbP_{\bta(\btheta)}(\bX_{\mK} = \bx_{\mK},\; \bX_{\mL \setminus \mK} \in \mX_{\mL \setminus \mK})
\= \mbP_{\bta(\btheta)}(\bX_{\mK} = \bx_{\mK}).
\ee
The second term on the right-hand side of \eqref{basic} can be dealt with by using the assumption that the between-neighborhood subgraphs are independent:
\beno
\mbP(\bY_{\mK} = \by_{\mK},\; \bY_{\mL \setminus \mK} \in \mbY_{\mL \setminus \mK})
\,=\, \mbP(\bY_{\mK} = \by_{\mK})\; \mbP(\bY_{\mL \setminus \mK} \in \mbY_{\mL \setminus \mK})
\,=\, \mbP(\bY_{\mK} = \by_{\mK}).
\ee
Collecting results gives
\beno
&& \mbP_{\bta(\btheta)}(\bX_{\mK} = \bx_{\mK},\; \bY_{\mK} = \by_{\mK},\; \bX_{\mL \setminus \mK} \in \mX_{\mL \setminus \mK},\; \bY_{\mL \setminus \mK} \in \mbY_{\mL \setminus \mK})\s
\\
\= \mbP_{\bta(\btheta)}(\bX_{\mK} = \bx_{\mK})\; \mbP(\bY_{\mK} = \by_{\mK})
\;=\; \mbP_{\bta(\btheta)}(\bX_{\mK} = \bx_{\mK},\; \bY_{\mK} = \by_{\mK}). 
\ee

\ttproof \ref{incomplete}.
By Proposition \ref{proposition.extendability},
we have,
for all $\btheta\in\bTheta \subseteq \{\btheta \in \mR^{\qqq}:\; \psi_{\mL}(\bta(\btheta)) < \infty\}$,
all $\mK \subseteq \mL$,
and all $\bx_{\mK} \in \mX_{\mK}$ and all $\by_{\mK} \in \mY_{\mK}$,
\beno
&& \mbP_{\bta(\btheta)}(\bX_{\mK} = \bx_{\mK},\; \bY_{\mK} = \by_{\mK},\; \bX_{\mL \setminus \mK} \in \mX_{\mL \setminus \mK},\; \bY_{\mL \setminus \mK} \in \mbY_{\mL \setminus \mK})\s
\\
\= \mbP_{\bta(\btheta)}(\bX_{\mK} = \bx_{\mK},\; \bY_{\mK} = \by_{\mK})
\ee
and,
for all $\bx_{\mL} \in \mX_{\mL}$, 
\beno
\dsum_{\bx_{\mL\setminus\mK}\, \in\, \mbX_{\mL\setminus\mK}}\, p_{\bta(\btheta)}(\bx_{\mL})
\= p_{\bta(\btheta)}(\bx_{\mK}),
\ee
where
\beno
p_{\bta(\btheta)}(\bx_{\mK})
\= \exp\left(\langle\bta(\btheta),\, s(\bx_{\mK})\rangle - \psi_{\mK}(\bta(\btheta))\right) \nu_{\mK}(\bx_{\mK}).
\ee
Here, 
$\bta(\btheta)$ and $s(\bx_{\mK})$ have dimension $\diminf = \max_{\mA \in \mL} \ppp_{\mA}$.
We note that $\psi_{\mK}(\bta(\btheta))$ satisfies
\beno
\psi_{\mK}(\bta(\btheta))
\= \dsum_{\mA \in \mK} \psi_{\mA}(\bta(\btheta)),
\ee
where 
\beno
\psi_{\mA}(\bta(\btheta))
\= \dsum_{\bx_{\mA} \in \mX_{\mA}} \exp\left(\langle\bta(\btheta),\, s(\bx_{\mA})\rangle\right)\, \nu_{\mA}(\bx_{\mA}),
& \mA \in \mK.
\ee
In addition,
note that $\psi_{\mK}(\bta(\btheta)) = \sum_{\mA \in \mK} \psi_{\mA}(\bta(\btheta))$ is finite provided $\psi_{\mL}(\bta(\btheta)) = \sum_{\mA \in \mL} \psi_{\mA}(\bta(\btheta))$ is finite.
In other words,
local dependence implies that $\bX_{\mK}$ is independent of $\bX_{\mL \setminus \mK}$ and the marginal density of $\bx_{\mK} \in \mX_{\mK}$ induced by $\mK \subseteq \mL$ is an exponential-family density with support $\mX_{\mK}$ and local dependence,
with natural parameter vector $\bta(\btheta)$ and sufficient statistic vector $s(\bx_{\mK})$.
\hide{
\citet[][Theorem 1]{ScHa13} proved that,
for all $L > K$ and all $\mathscr{X}_L \subseteq \mbX_L$ and $\mathscr{Y}_L \subseteq \mbY_L$,
\beno
&& \mbP(\bX_K \in \mathscr{X}_K,\; \bY_K \in \mathscr{Y}_K)\s
\\
\= \mbP(\bX_K \in \mathscr{X}_K,\; \bY_K \in \mathscr{Y}_K,\; \bX_{L \setminus K} \in \mbX_{L \setminus K},\; \bY_{L \setminus K} \in \mbY_{L \setminus K}),
\ee
where $\bX_{L \setminus K} \in \mathbb{X}_{L \setminus K}$ and $\bY_{L \setminus K} \in \mathbb{Y}_{L \setminus K}$ correspond to $\bX_{L} \in \mathbb{X}_{L}$ and $\bY_{L} \in \mathbb{Y}_{L}$ excluding $\bX_{K} \in \mathbb{X}_{K}$ and $\bY_{K} \in \mathbb{Y}_{K}$,
respectively.
}
As a result,
Theorem \ref{nonfull} can be applied to the exponential family with support $\mX_{\mK}$ and local dependence,
with natural parameter vector $\bta(\btheta)$ and sufficient statistic vector $s(\bx_{\mK})$.
Observe that the conditions of Theorem \ref{incomplete} ensure that all conditions of Theorem \ref{nonfull} are satisfied for all $\mK \subseteq \mL$.
Therefore,
for full exponential families,
Theorem \ref{nonfull} shows that,
for all $\epsilon > 0$ small enough so that $\mB(\bthetas, \epsilon) \subseteq \interior(\bTheta)$,
there exist $\kappa_1(\epsilon) > 0$ and $C_1 > 0$ such that $\widehat\btheta_{\mK}$ exists, is unique, and is contained in $\mB(\bthetas, \epsilon)$ with probability
\be
\nonumber
\mbP\left(\bthetah_{\mK} \in \mB(\bthetas, \epsilon) \subseteq \interior(\bTheta)\right)
\geq 1 - 2\, \exp\left(- \dfrac{\kappa_1(\epsilon)^2\; C_1\; |\mK|}{\qqq\, \norm{\mK}_\infty^{4\, (2 - \alpha)}} + \log \qqq\right),
\ee
where $\norm{\mK}_{\infty} = \max_{\mA \in \mK} |\mA|$.
For non-full, curved exponential families,
Theorem \ref{nonfull} shows that,
for all $\epsilon > 0$ small enough so that $\mB(\bthetas, \epsilon) \subseteq \interior(\bTheta)$,
there exist $\kappa_2(\epsilon) > 0$ and $C_2 > 0$ such that $\widehat\btheta_{\mK}$ exists and is contained in $\mB(\bthetas, \epsilon)$ with probability
\be
\nonumber
\mbP\left(\widehat\btheta_{\mK}\; \in\; \mB(\bthetas, \epsilon)\; \subseteq\; \interior(\bTheta)\right)
\;\geq\; 1 - 2\, \exp\left(- \dfrac{\kappa_2(\epsilon)^2\, C_2\, |\mK|}{\diminf\, \norm{\mK}_\infty^{4\, (2 - \alpha)}} + \log \diminf\right).
\ee
Thus,
there exist $\kappa(\epsilon) = \min(\kappa_1(\epsilon), \kappa_2(\epsilon)) > 0$ and $C = \min(C_1, C_2) > 0$ such that
\be
\nonumber
\mbP\left(\bthetah_{\mK} \in \mB(\bthetas, \epsilon) \subseteq \interior(\bTheta)\right)
\gte 1 - 2\, \exp\left(- \dfrac{\kappa(\epsilon)^2\, C\, |\mK|}{\diminf\, \norm{\mK}_\infty^{4\, (2 - \alpha)}} + \log \diminf\right)\s
\\
\gte 1 - 2\, \exp\left(- \dfrac{\kappa(\epsilon)^2\, C\, |\mK|}{\diminf\, \norm{\mL}_\infty^{4\, (2 - \alpha)}} + \log \diminf\right),
\ee
where we used the fact that $\norm{\mK}_{\infty} = \max_{\mA \in \mK} |\mA| \leq \max_{\mA \in \mL} |\mA| = \norm{\mL}_\infty$ and $\diminf = \max_{\mA \in \mL} \ppp_{\mA} = \qqq$ in full exponential families with natural parameter vectors of the form $\bta(\btheta) = \btheta$.

\s

\ccproof \ref{corollary.curved.incomplete}.
Corollary \ref{corollary.curved.incomplete} follows from Theorem \ref{incomplete},
because the conditions of Theorem \ref{incomplete} are satisfied.
The proof of Corollary \ref{corollary.curved.incomplete} resembles the proof of Corollary \ref{corollary.curved} and is therefore omitted.

\s

\section{Proofs: Auxiliary lemmas}
\label{m.id}

We prove auxiliary lemmas that are useful for verifying the conditions of Theorem \ref{nonfull} (Appendix \ref{appendix.mle.id.curved}) and Theorem \ref{theorem:parameters.correct} (Appendix \ref{appendix.m.id}).
We start with two lemmas that are useful for bounding expectations of sufficient statistics.
Throughout,
we write
\be
\nonumber
\Lambda(\bta)
&=& \min\limits_{i\in\mA_k\, <\, j\in\mA_k,\; k=1, \dots, K}\; \min\limits_{\bx_{-i,j}\, \in\, \mbX_{-i,j}} \, \mbP_{\bta}(X_{i,j} = 1 \mid \bX_{-i,j} = \bx_{-i,j})
\ee
and
\be
\nonumber
\Omega(\bta)
&=& \max\limits_{i\in\mA_k\, <\, j\in\mA_k,\; k=1, \dots, K}\; \max\limits_{\bx_{-i,j}\, \in\, \mbX_{-i,j}} \mbP_{\bta}(X_{i,j} = 1 \mid \bX_{-i,j} = \bx_{-i,j}),
\ee
where $\bta$ is the natural parameter vector of the exponential family under consideration and $\bx_{-i,j} \in \mbX_{-i,j}$ corresponds to $\bx \in \mbX$ excluding $x_{i,j} \in \mbX_{i,j}$ ($i\in\mA_k < j\in\mA_k$, $k = 1, \dots, K$).
Let $f: \mbX \mapsto \mR$ be a function of the random graph.
We denote by $\mbE_{\Lambda(\bta)}\, f(\bX)$ and $\mbE_{\Omega(\bta)}\, f(\bX)$ the expectations of $f: \mbX \mapsto \mR$ under the Bernoulli random graph model with probabilities $\Lambda(\btheta)$ and $\Omega(\btheta)$,
respectively,
provided $\mbE_{\Lambda(\bta)}\, f(\bX)$ and $\mbE_{\Omega(\bta)}\, f(\bX)$ exist.
Here,
the Bernoulli random graph model refers to the exponential-family random graph model which assumes that the edge variables $X_{i,j}$ are independent Bernoulli random variables with probabilities $\Lambda(\btheta)$ and $\Omega(\btheta)$,
respectively.

\s\s

\begin{lemma} 
\label{mean.bound} 
Consider a full or non-full, canonical or curved exponential family with support $\mbX = \{0, 1\}^{\sum_{k=1}^K \binom{|\mA_k|}{2}}$ and local dependence and natural parameter vector $\bta \in \interior(\fullspace)$.
Then there exists $C(\bta) \in [\Lambda(\bta),\, \Omega(\bta)]$ such that
\be
\nonumber
\mbE_{\bta}\; \dsum_{k=1}^K\; \dsum_{i\in\mA_k\, <\, j\in\mA_k} X_{i,j}
&=& C(\bta) \dsum_{k=1}^K \dis\binom{|\mA_k|}{2},
\ee
where $C(\bta) \in [\Lambda(\bta),\, \Omega(\bta)]$ denotes the probability of an edge under the\linebreak 
Bernoulli$(C(\bta))$ random graph model.
In addition, 
if $f: \mbX \mapsto \mR$ is a function of the random graph that is non-decreasing in the number of edges,
then 
\be 
\nonumber
\mbE_{\Lambda(\bta)}\, f(\bX)
&\leq& \mbE_{\bta}\, f(\bX) 
&\leq& \mbE_{\Omega(\bta)}\, f(\bX),
\ee 
where $\mbE_{\Lambda(\bta)}\, f(\bX)$ and $\mbE_{\Omega(\bta)}\, f(\bX)$ are the expectations of $f: \mbX \mapsto \mR$ under the Bernoulli random
 graph model with probabilities $\Lambda(\btheta)$ and $\Omega(\btheta)$,
respectively.
If $f: \mbX \mapsto \mR$ is a function of the random graph that is non-increasing in the number of edges,
then the inequalities are reversed.
\end{lemma}

\s

\llproof \ref{mean.bound}. 
To ease the presentation,
we consider $K = 1$ neighborhood and drop the subscript $k$ from all neighborhood-dependent quantities.
The extension to $K \geq 2$ neighborhoods is straightforward.
Observe that $\mbE_{\bta}\, X_{i,j} = \mbP_{\bta}(X_{i,j} = 1)$ can be written as
\be
\nonumber 
\mbP_{\bta}(X_{i,j} = 1)
\;=\; \dsum_{\bx_{-i,j}\, \in\, \mbX_{-i,j}} \mbP_{\bta}(X_{i,j} = 1 \, | \, \bX_{-i,j} = \bx_{-i,j})\; \mbP_{\bta}(\bX_{-i,j} = \bx_{-i,j}),
\ee
which implies that $\mbE_{\bta}\, X_{i,j}$ is bounded below by
\be
\nonumber
\dsum_{\bx_{-i,j}\, \in\, \mbX_{-i,j}}\, \Lambda(\bta)\, \mbP_{\bta}(\bX_{-i,j} = \bx_{-i,j})
&=& \Lambda(\bta)
\ee 
and bounded above by
\be
\nonumber
\dsum_{\bx_{-i,j}\, \in\, \mbX_{-i,j}}\, \Omega(\bta)\, \mbP_{\bta}(\bX_{-i,j} = \bx_{-i,j})
&=& \Omega(\bta).
\ee 
Since the expectation $\mbE_{\bta}\, X_{i,j}$ is contained in the convex set $[\Lambda(\bta),\, \Omega(\bta)]$,
there exists $C(\bta) \in [\Lambda(\bta), \Omega(\bta)]$ such that
\be
\nonumber
\mbE_{\bta}\, \dsum_{i\in\mA\, <\, j\in\mA} X_{i,j}
\= \dsum_{i\in\mA\, <\, j\in\mA} \mbE_{\bta}\, X_{i,j}
\= C(\bta) \dis\binom{|\mA|}{2}.
\ee 
In addition,
\citet[][Corollary 1, p.\ 306]{Bu10} proved that,
if $f: \mbX \mapsto \mR$  is a function of the random graph that is non-decreasing in the number of edges,
then
\be 
\nonumber
\mbE_{\Lambda(\bta)}\, f(\bX) 
&\leq& \mbE_{\bta}\, f(\bX) 
&\leq& \mbE_{\Omega(\bta)}\, f(\bX),
\ee
where the inequalities are reversed when $f: \mbX \mapsto \mR$ is a function of the random graph that is non-increasing in the number of edges.
The result follows from the recursive factorization of the probability mass function $\mbP_{\bta}(\bX = \bx)$.
Denote the elements of the sequence of within-neighborhood edge variables $\bX$ by $X_1, \dots, X_w$ and the corresponding sample spaces by $\mbX_1, \dots, \mbX_w$,
where  $w = {|\mA| \choose 2}$.
Then,
for all $\bx \in \mbX$,
we have
\beno
\mbP_{\bta}(\bX = \bx)
\= \dprod_{i=1}^w \mbP_{\bta}(X_i = x_i \mid X_j = x_j,\; j = 1, \dots, i-1),
\ee
where the conditional probabilities $\mbP_{\bta}(X_i = 1 \mid X_j = x_j,\; j = 1, \dots, i-1)$ are bounded by
\beno
&& \Lambda(\bta)
\= \min\limits_{1\, \leq\, i\, \leq\, w}\; \min\limits_{\bx_{-i}\, \in\, \mbX_{-i}} \, \mbP_{\bta}(X_{i} = 1 \mid \bX_{-i} = \bx_{-i})\s
\\
&&\lte \mbP_{\bta}(X_i = 1 \mid X_j = x_j,\, j = 1, \dots, i-1)\s
\\
&&\lte \max\limits_{1\, \leq\, i\, \leq\, w}\; \max\limits_{\bx_{-i}\, \in\, \mbX_{-i}} \mbP_{\bta}(X_i = 1 \mid \bX_{-i} = \bx_{-i})
\= \Omega(\bta),
\ee
where $\bx_{-i} \in \mbX_{-i}$ corresponds to $\bx \in \mbX$ excluding $x_{i} \in \mbX_{i}$ ($i = 1, \dots, w$).
Therefore,
$\mbP_{\bta}(\bX = \bx)$ can be bounded below and above by the probability mass function of the Bernoulli random graph model with probabilities $\Lambda(\bta)$ and $\Omega(\bta)$,
respectively.
As a result,
the cumulative distribution function and expectation of any function $f: \mbX \mapsto \mR$ of the random graph that is non-decreasing in the number of edges can be bounded by the corresponding cumulative distribution functions and expectations of the Bernoulli random graph model with probabilities $\Lambda(\bta)$ and $\Omega(\bta)$ \citep[][Corollary 1, p.\ 306]{Bu10}:
\be
\nonumber
\mbE_{\Lambda(\bta)}\, f(\bX)
&\leq& \mbE_{\bta}\, f(\bX)
&\leq& \mbE_{\Omega(\bta)}\, f(\bX),
\ee
where the inequalities are reversed when $f: \mbX \mapsto \mR$ is a function of the random graph that is non-increasing in the number of edges.

\s\s

\begin{lemma}
\label{tran_edge_alpha_exp}
Consider a full or non-full, canonical or curved exponential family with support $\mbX = \{0, 1\}^{\sum_{k=1}^K \binom{|\mA_k|}{2}}$ and local dependence and natural parameter vector $\bta \in \interior(\fullspace)$.
Let $f: \mbX \mapsto \mR$ be the number of within-neighborhood transitive edges defined by
\be
\nonumber
f(\bx)
&=& \dsum_{k=1}^{K}\; \dsum_{i\in\mA_k\, <\, j\in\mA_k}\, x_{i,j}\, \max\limits_{h \in \mA_k,\; h\, \neq\, i, j}\,  x_{i,h}\, x_{j,h},
& \bx \in \mbX,
\ee
which may or may not be a sufficient statistic of the exponential family under consideration.
If $|\mA_k| \geq 3$ ($k = 1, \dots, K$),
then there exists $C(\bta) > 0$ such that
\be
\nonumber
\mbE_{\bta}\, f(\bX)
&=& C(\bta)\, \dsum_{k=1}^K \dis\binom{|\mA_k|}{2},
\ee
where $C(\bta)$ satisfies,
for some $\lambda \in (0, 1)$,
\be
\nonumber
0
\;<\;
\lambda\; \Lambda(\bta)^3 + (1 - \lambda)\, \Omega(\bta)^3
&\leq& C(\bta)
&\leq& \lambda\; \Lambda(\bta) + (1 - \lambda)\, \Omega(\bta)
\;<\; 1.
\ee
\end{lemma}
\llproof \ref{tran_edge_alpha_exp}.
To ease the presentation,
we consider $K = 1$ neighborhood and drop the subscript $k$ from all neighborhood-dependent quantities.
The extension to $K \geq 2$ neighborhoods is straightforward.
Since the number of transitive edges $f: \mbX \mapsto \mR$ is a function of the random graph that is non-decreasing in the number of edges,
Lemma \ref{mean.bound} implies that the expectation $\mbE_{\bta}\, f(\bX)$ can be bounded by the expectations $\mbE_{\Lambda(\bta)}\, f(\bX)$ and $\mbE_{\Omega(\bta)}\, f(\bX)$ under the Bernoulli random graph model with edge probabilities $\Lambda(\bta)$ and $\Omega(\bta)$,
respectively:
\be
\nonumber
\mbE_{\Lambda(\bta)}\, f(\bX)
&\leq& \mbE_{\bta}\, f(\bX)
&\leq& \mbE_{\Omega(\bta)}\, f(\bX),
\ee
where $\Lambda(\bta)$ and $\Omega(\bta)$ are defined in the introduction of Appendix \ref{m.id} and all expectations exist,
because the sample space $\mX$ is finite.
The expectation $\mbE_{\Lambda(\bta)} \, f(\bX)$ can be written as
\be
\nonumber
\mbE_{\Lambda(\bta)}\, f(\bX)
\hide{
&=& \mbE_{\Lambda(\bta)}\, \dsum_{i\in\mA\, <\, j\in\mA}\, X_{i,j}\, \max\limits_{h \in \mA,\; h\, \neq\, i, j}\,  X_{i,h}\, X_{j,h}\s
\\
}
&=& \mbE_{\Lambda(\bta)} \left[\dsum_{i\in\mA\, <\, j\in\mA} X_{i,j}\, \max\limits_{h \in \mA,\; h\, \neq\, i, j}\,  X_{i,h}\, X_{j,h}\right]\s\s
\\
&=& \mbE_{\Lambda(\bta)} \left[\dsum_{i\in\mA\, <\, j\in\mA} X_{i,j}\; \one\left(\dsum_{h \in \mA,\; h\, \neq\, i, j}\,  X_{i,h}\, X_{j,h} \geq 1\right)\right]\s\s
\\
&=& \dsum_{i\in\mA\, <\, j\in\mA}\; \mbE_{\Lambda(\bta)} \left[X_{i,j}\; \one\left(\dsum_{h \in \mA,\; h\, \neq\, i, j}\,  X_{i,h}\, X_{j,h} \geq 1\right)\right],
\ee
where $\one(\sum_{h \in \mA,\; h\, \neq\, i, j}\,  x_{i,h}\, x_{j,h} \geq 1)$ is an indicator function,
which is $1$ if nodes $i$ and $j$ have one or more shared partners in neighborhood $\mA$ and is $0$ otherwise.
Here, 
we exploited the fact that the number of transitive edges in neighborhood $\mA$ is equal to the number of pairs of nodes with one or more edgewise shared partners.
To evaluate the expectation $\mbE_{\Lambda(\bta)}[X_{i,j}\; \one(\sum_{h \in \mA,\; h\, \neq\, i, j}\,  X_{i,h}\, X_{j,h} \geq 1)]$,
we take advantage of the independence of edge variables under the Bernoulli$(\Lambda(\bta))$ random graph model,
which implies that
\be
\nonumber
&& \mbE_{\Lambda(\bta)} \left[X_{i,j}\; \one\left(\dsum_{h \in \mA,\; h\, \neq\, i, j}\,  X_{i,h}\, X_{j,h} \geq 1\right)\right]\s
\\
\= \mbE_{\Lambda(\bta)} (X_{i,j})\;\; \mbE_{\Lambda(\bta)}\, \left[\one\left(\dsum_{h \in \mA,\; h\, \neq\, i, j}\,  X_{i,h}\, X_{j,h} \geq 1\right)\right]\s
\\
\hide{
= \Lambda(\bta)\, \mbP_{\Lambda(\bta)}\left(\max\limits_{h \in \mA,\; h\, \neq\, i, j} X_{i,h}\, X_{j,h} = 1\right)
}
\= \Lambda(\bta)\; \mbP_{\Lambda(\bta)}\left(\dsum_{h \in \mA,\; h\, \neq\, i, j} X_{i,h}\, X_{j,h}\, \geq\, 1\right).
\ee
In addition,
the independence of edge variables under the Bernoulli$(\Lambda(\bta))$ random graph model implies that,
for all $i \in \mA$ and all $j \in \mA$ such that $i \neq j$, 
the distribution of $\sum_{h \in \mA,\; h\, \neq\, i, j}\, X_{i,h}\, X_{j,h}$ is Binomial($|\mA| - 2$, $\Lambda(\bta)^2$).
Therefore,
\be
\nonumber
\mbP_{\Lambda(\bta)}\left(\dsum_{h \in \mA,\; h\, \neq\, i, j}\, X_{i,h}\, X_{j,h}\; \geq\; 1\right)
\hide{
&=& 1 - \mbP_{\Lambda(\bta)}\left(\dsum_{h \in \mA,\; h\, \neq\, i, j}\, X_{i,h}\, X_{j,h} = 0\right) \s
\\
}
&=& 1 - (1 - \Lambda(\bta)^2)^{|\mA| - 2},
\ee
which implies that
\be
\nonumber
\mbE_{\Lambda(\bta)}\, f(\bX)
&=& \Lambda(\bta)\, [1 - (1 - \Lambda(\bta)^2)^{|\mA| - 2}]\, \dis\binom{|\mA|}{2}.
\ee
Along the same lines,
it can be shown that the expectation $\mbE_{\Omega(\bta)}\, f(\bX)$ is given by
\be
\nonumber
\mbE_{\Omega(\bta)}\, f(\bX)
&=& \Omega(\bta)\, [1 - (1 - \Omega(\bta)^2)^{|\mA| - 2}]\, \dis\binom{|\mA|}{2}.
\ee
Let $C_1(\bta) = \Lambda(\bta)\, [1 - (1 - \Lambda(\bta)^2)^{|\mA| - 2}]$ and $C_2(\bta) = \Omega(\bta)\, [1 - (1 - \Omega(\bta)^2)^{|\mA| - 2}]$.
Observe that $\Lambda(\bta) \in (0, 1)$ implies that $C_1(\bta)$ is bounded below by $\Lambda(\bta)^3$ and bounded above by $\Lambda(\bta)$ for all $|\mA| \geq 3$.
Likewise,
$C_2(\bta)$ is bounded below by $\Omega(\bta)^3$ and bounded above by $\Omega(\bta)$ for all $|\mA| \geq 3$.
Since the expectation $\mbE_{\bta}\, f(\bX)$ is contained in the convex set $[C_1(\bta)\, \binom{|\mA|}{2},\; C_2(\bta)\, \binom{|\mA|}{2}]$,
there exists $\lambda \in (0, 1)$ such that
\be
\nonumber
\mbE_{\bta}\, f(\bX)
&=& \lambda\, C_1(\bta)\, \dis\binom{|\mA|}{2} + (1 - \lambda)\, C_2(\bta)\, \dis\binom{|\mA|}{2}
\hide{\s
\\
&=& (\lambda\, C_1(\bta) + (1 - \lambda)\, C_2(\bta))\, \dis\binom{|\mA|}{2}\s
\\
}
&=& C(\bta)\, \dis\binom{|\mA|}{2},
\ee
where $C(\bta)$ satisfies
\be
\nonumber
0
\;<\;
\lambda\; \Lambda(\bta)^3 + (1 - \lambda)\, \Omega(\bta)^3
&\leq& C(\bta)
&\leq& \lambda\; \Lambda(\bta) + (1 - \lambda)\, \Omega(\bta)
\;<\; 1.
\ee
\hide{
where $C(\bta)$ is bounded below by
\be
\nonumber
C(\bta)
&=& \lambda\, C_1(\bta) + (1 - \lambda)\, C_2(\bta)
&\geq& \lambda\; \Lambda(\bta)^3 + (1 - \lambda)\, \Omega(\bta)^3
&>& 0
\ee
and bounded above by
\be
\nonumber
C(\bta)
&=& \lambda\, C_1(\bta) + (1 - \lambda)\, C_2(\bta)
&\leq& \lambda\; \Lambda(\bta) + (1 - \lambda)\, \Omega(\bta)
&<& 1,
\ee
which implies that $C(\bta)$ satisfies
\be
\nonumber
0
\;<\;
\lambda\, \Lambda(\bta)^3 + (1 - \lambda)\, \Omega(\bta)^3
&\leq& C(\bta)
&\leq& \lambda\, \Lambda(\bta) + (1 - \lambda)\, \Omega(\bta)
\;<\; 1.
\ee
}

\subsection{Auxiliary lemmas: Theorem \ref{nonfull}}
\label{appendix.mle.id.curved}

We prove auxiliary lemmas that are useful for verifying the conditions of Theorem \ref{nonfull}.

\s

\begin{lemma}
\label{lemma:gw_grad_ident}
Consider a curved exponential-family random graph with\linebreak 
within-neighborhood edge and geometrically weighted edgewise shared partner terms as defined in Section \ref{sec:curved}.
Let $\bTheta = \mR \times (0, 1)$.
Then the map $\bta: \interior(\bTheta) \mapsto \interior(\fullspace)$ is one-to-one and,
for all $\epsilon > 0$ small enough so that $\mathscr{B}(\bthetas, \epsilon) \subseteq \interior(\bTheta)$,
there exists $\gamma(\epsilon) > 0$ such that,
for all $\btheta \in \bTheta \setminus \mB(\bthetas, \epsilon)$,
we have $\bta(\btheta) \in \fullspace \setminus \mB(\bta(\bthetas),\, \gamma(\epsilon))$ provided $|\mA_k| \geq 4$ ($k = 1, \dots, K$).
\end{lemma}

\llproof \ref{lemma:gw_grad_ident}.
It is straightforward to show that the map $\bta: \interior(\bTheta) \mapsto \interior(\fullspace)$ is one-to-one and that, 
for all $\epsilon > 0$ small enough so that $\mathscr{B}(\bthetas, \epsilon) \subseteq \interior(\bTheta)$,
there exists $\gamma(\epsilon) > 0$ such that,
for all $\btheta \in \bTheta \setminus \mB(\bthetas, \epsilon)$,
\beno
\label{summ}
\norm{\bta(\btheta) - \bta(\bthetas)}_2
\= \sqrt{\dsum_{i=1}^{\diminf} (\eta_i(\btheta) - \eta_i(\bthetas))^2}
\gte \gamma(\epsilon).
\ee
To see that,
note that $\diminf = \norm{\mA}_\infty - 1 \geq 3$ since $|\mA_k| \geq 4$ ($k = 1, \dots, K$).
Therefore,
\beno
\norm{\bta(\btheta) - \bta(\bthetas)}_2
\= \sqrt{\dsum_{i=1}^{\diminf} (\eta_i(\btheta) - \eta_i(\bthetas))^2}
\gte \sqrt{\dsum_{i=1}^3 (\eta_i(\btheta) - \eta_i(\bthetas))^2},
\ee
where
\beno
\eta_1(\btheta)
\= \theta_1\s
\\
\eta_2(\btheta)
\= 1\s
\\
\eta_3(\btheta)
\= 2 - \theta_2.
\ee
Thus,
for all $\bthetas \in \interior(\bTheta) = \mR \times (0, 1)$ and all $\bepsilon \in \mR^2$ such that $\bthetas + \bepsilon \in \interior(\bTheta) = \mR \times (0, 1)$,
\beno
(\eta_1(\bthetas+\bepsilon) - \eta_1(\bthetas))^2
\= (\thetas_1 + \delta_1 - \thetas_1)^2\; =\; \delta_1^2\s
\\
(\eta_2(\bthetas + \bepsilon) - \eta_2(\bthetas))^2
\= 0\s
\\
(\eta_3(\bthetas + \bepsilon) - \eta_3(\bthetas))^2
\= (2 - \thetas_2 - \delta_2 - 2 + \thetas_2)^2 \;=\; \delta_2^2,
\ee
which implies that
\beno
\norm{\bta(\bthetas + \bepsilon) - \bta(\bthetas)}_2
\,\geq\, \sqrt{\dsum_{i=1}^{3} (\eta_i(\bthetas+\bepsilon) - \eta_i(\bthetas))^2}
\,=\, \sqrt{\delta_1^2 + \delta_2^2}
\,=\, \norm{\bepsilon}_2.
\ee
Therefore,
the $\ell_2$-distance of $\bta(\bthetas + \bepsilon)$  from $\bta(\bthetas) \in \interior(\fullspace)$ in the natural parameter space $\fullspace$ is a strictly increasing function of the $\ell_2$-distance $\epsilon = \norm{\bepsilon}_2$ of $\bthetas + \bepsilon$ from $\bthetas \in \interior(\bTheta)$ in the parameter space $\bTheta$.
As a result,
for all $\epsilon > 0$,
there exists $\gamma(\epsilon) > 0$ such that,
for all $\btheta \in \bTheta \setminus \mB(\bthetas, \epsilon)$,
\beno
\norm{\bta(\bthetas) - \bta(\btheta)}_2
\= \sqrt{\dsum_{i=1}^{\diminf} (\eta_i(\bthetas) - \eta_i(\btheta))^2}
\gte \gamma(\epsilon).
\ee
\s

\begin{lemma}
\label{lemma:theorem_2_ident}
Consider a curved exponential-family random graph with\linebreak 
within-neighborhood edge and geometrically weighted edgewise shared partner terms as defined in Section \ref{sec:curved}.
Let $\bTheta = \mbR \times (0, 1)$ and assume that $\bthetas \in \interior(\bTheta)$.
Then,
for all $\epsilon > 0$ small enough so that $\mathscr{B}(\bthetas, \epsilon) \subseteq \interior(\bTheta)$,
there exists $\gamma(\epsilon) > 0$ such that,
for all $\btheta \in \bTheta \setminus \mB(\bthetas, \epsilon)$,
we have $\bta(\btheta) \in \fullspace \setminus \mB(\bta(\bthetas),\, \gamma(\epsilon))$.
In addition,
there exists $\delta(\epsilon) > 0$ such that,
for all $\bta(\btheta) \in \fullspace \setminus \mB(\bta(\bthetas), \gamma(\epsilon))$,
\be
\nonumber
\norm{\bmu(\bta(\bthetas)) - \bmu(\bta(\btheta))}_2
\gte \delta(\epsilon) \dsum_{k=1}^K {|\mA_k| \choose 2}^{3/4},
\ee
provided $|\mA_k| \geq 4$ ($k = 1, \dots, K$).
Therefore, 
identifiability condition \eqref{nonfullid} of Theorem \ref{nonfull} is satisfied with $\alpha = 3/4$ provided $|\mA_k| \geq 4$ ($k = 1, \dots, K$).
\end{lemma}

\llproof \ref{lemma:theorem_2_ident}.
To ease the presentation,
we consider $K = 1$ neighborhood and drop the subscript $k$ from all neighborhood-dependent quantities.
The extension to $K \geq 2$ neighborhoods is straightforward.
By Lemma \ref{lemma:gw_grad_ident},
for all $\epsilon > 0$ small enough so that $\mathscr{B}(\bthetas, \epsilon) \subseteq \interior(\bTheta)$,
there exists $\gamma(\epsilon) > 0$ such that,
for all $\btheta \in \bTheta \setminus \mB(\bthetas, \epsilon)$,
we have $\bta(\btheta) \in \fullspace \setminus \mB(\bta(\bthetas),\, \gamma(\epsilon))$.
Therefore,
it is enough to show that there exists $\delta(\epsilon) > 0$ such that,
for all $\bta(\btheta) \in \fullspace \setminus \mB(\bta(\bthetas), \gamma(\epsilon))$,
\beno
\norm{\bmu(\bta(\bthetas)) - \bmu(\bta(\btheta))}_2
&\geq& \delta(\epsilon)\, \dis\binom{|\mA|}{2}^\alpha
& \mbox{for some}
& 0 \leq \alpha \leq 1.
\ee
To do so,
observe that
\beno
&& \norm{\bmu(\bta(\bthetas)) - \bmu(\bta(\btheta))}_2\s
\\
\= \sqrt{(\mu_1(\bta(\bthetas)) - \mu_1(\bta(\btheta)))^2 + \dsum_{i=2}^{|\mA| - 1} (\mu_i(\bta(\bthetas)) - \mu_i(\bta(\btheta)))^2}
&>& 0,
\ee
where strict positivity follows from the fact that the maps $\bta: \interior(\bTheta) \mapsto \interior(\fullspace)$ and $\bmu: \interior(\fullspace) \mapsto \rint(\mM)$ are one-to-one and hence $\btheta \neq \bthetas$ implies $\norm{\bmu(\bta(\bthetas)) - \bmu(\bta(\btheta))}_2 > 0$;
note that the map $\bta: \interior(\bTheta) \mapsto \interior(\fullspace)$ is one-to-one by Lemma \ref{lemma:gw_grad_ident} while the map $\bmu: \interior(\fullspace) \mapsto \rint(\mM)$ is one-to-one by classic exponential-family theory \citep[][Theorem 3.6, p.\ 74]{Br86}.

We distinguish two cases,
the case $\mu_1(\bta(\bthetas)) \neq \mu_1(\bta(\btheta))$ and the case\linebreak 
$\mu_1(\bta(\bthetas)) = \mu_1(\bta(\btheta))$.
We note that $\btheta \neq \bthetas$ and $\mu_1(\bta(\bthetas)) = \mu_1(\bta(\btheta))$ imply
$\sum_{i=2}^{|\mA| - 1} (\mu_i(\bta(\bthetas)) - \mu_i(\bta(\btheta)))^2 > 0$,
because $\btheta \neq \bthetas$ implies $\norm{\bmu(\bta(\bthetas)) - \bmu(\bta(\btheta))}_2 > 0$.

\s

{\bf Case $\mu_1(\bta(\bthetas)) \neq \mu_1(\bta(\btheta))$.}
We have 
\beno
\norm{\bmu(\bta(\bthetas)) - \bmu(\bta(\btheta))}_2
\gte |\mu_1(\bta(\bthetas)) - \mu_1(\bta(\btheta))|
&>& 0.
\ee
Observe that $\mu_1(\bta(\btheta)) = \mbE_{\bta(\btheta)}\, s_1(\bX)$ is the expected number of edges in neighborhood $\mA$.
Therefore,
Lemma \ref{mean.bound} can be invoked to show that there exist $C_1(\bta(\bthetas)) > 0$ and $C_1(\bta(\btheta)) > 0$ such that
\beno
\label{exp_1}
|\mu_1(\bta(\bthetas)) - \mu_1(\bta(\btheta))|
\= |C_1(\bta(\bthetas)) - C_1(\bta(\btheta))|\, \dis\binom{|\mA|}{2},
\ee
where $\Lambda(\bta(\bthetas)) \leq C_1(\bta(\bthetas)) \leq \Omega(\bta(\bthetas))$ and $\Lambda(\bta(\btheta)) \leq C_1(\bta(\btheta)) \leq \Omega(\bta(\btheta))$ for all $\btheta \in \interior(\bTheta) = \mbR \times (0, 1)$.
We note that $\Lambda(\bta(\btheta))$ and $\Omega(\bta(\btheta))$ are defined in the introduction of Appendix \ref{m.id} and that $\Lambda(\bta(\btheta)) \geq [1 + \exp(-\theta_1)]^{-1} > 0$ while $\Omega(\bta(\btheta))$ satisfies $0 < \Lambda(\bta(\btheta)) \leq \Omega(\bta(\btheta)) < 1$ for all $\btheta \in \interior(\bTheta) = \mbR \times (0, 1)$.

\s

{\bf Case $\mu_1(\bta(\bthetas)) = \mu_1(\bta(\btheta))$.}
As pointed out above,
$\btheta \neq \bthetas$ and $\mu_1(\bta(\bthetas)) = \mu_1(\bta(\btheta))$ imply
$\sum_{i=2}^{|\mA| - 1} (\mu_i(\bta(\bthetas)) - \mu_i(\bta(\btheta)))^2 > 0$,
where\linebreak 
$\mu_i(\bta(\btheta)) = \mbE_{\bta(\btheta)}\, s_i(\bX)$ is the expected number of pairs of nodes with $i-1$ edgewise shared partners in neighborhood $\mA$ ($i = 2, \dots, |\mA| - 1$).
Bounding the term $\sum_{i=2}^{|\mA| - 1} (\mu_i(\bta(\bthetas)) - \mu_i(\bta(\btheta)))^2$ is more challenging than bounding the term $(\mu_1(\bta(\bthetas)) - \mu_1(\bta(\btheta)))^2$,
because the numbers of pairs of nodes with $i-1$ edgewise shared partners are neither non-decreasing nor non-increasing functions of the number of edges ($i = 2, \dots, |\mA|-1$).
Therefore,
Lemma \ref{mean.bound} cannot be applied to the expectations $\mu_i(\bta(\btheta)) = \mbE_{\bta(\btheta)}\, s_i(\bX)$ ($i = 2, \dots, |\mA| - 1$).
But it turns out to be possible to bound $\sum_{i=2}^{|\mA| - 1} (\mu_i(\bta(\bthetas)) - \mu_i(\bta(\btheta)))^2$ from below in terms of absolute deviations of the expected numbers of transitive edges under $\bta(\bthetas)$ and $\bta(\btheta)$.
The advantage of doing so is that the number of transitive edges is a non-decreasing function of the number of edges and hence Lemma \ref{mean.bound} can be applied via Lemma \ref{tran_edge_alpha_exp}.
To see that the expected numbers of pairs of nodes with one or more edgewise shared partners is related to the expected number of transitive edges,
note that the expected number of transitive edges in neighborhood $\mA$ can be written as
\beno
&& \mbE_{\bta(\btheta)}\; \dsum_{a\in\mA\, <\, b\in\mA}\; X_{a,b}\; \max\limits_{c \in \mA,\; c\, \neq\, a, b}\,  X_{a,c}\, X_{b,c}\s
\\
\= \mbE_{\bta(\btheta)}\; \dsum_{i=2}^{|\mA| - 1}\; \dsum_{a\in\mA\, <\, b\in\mA}\, X_{a,b}\; \one\left(\sum_{c \in \mA,\; c\, \neq\, a, b}\,  X_{a,c}\, X_{b,c} = i-1\right)\s
\\
\= \dsum_{i=2}^{|\mA| - 1}\; \mbE_{\bta(\btheta)} \dsum_{a\in\mA\, <\, b\in\mA}\, X_{a,b}\; \one\left(\sum_{c \in \mA,\; c\, \neq\, a, b}\,  X_{a,c}\, X_{b,c} = i-1\right)\s
\\
\= \dsum_{i=2}^{|\mA| - 1}\; \mu_i(\bta(\btheta)),
\ee
which shows that $\sum_{i=2}^{|\mA| - 1}\; \mu_i(\bta(\btheta))$ is equal to the expected number of transitive edges under $\bta(\btheta)$. 
To bound $\sum_{i=2}^{|\mA| - 1}\, (\mu_i(\bta(\bthetas)) - \mu_i(\bta(\btheta)))^2$ from below in terms of absolute deviations of the expected numbers of transitive edges under $\bta(\bthetas)$ and $\bta(\btheta)$,
it is convenient to operate in the natural parameter space $\fullspace = \mR^{|\mA| - 1}$ of the full exponential family with natural parameter vector $\bta$ and sufficient statistic vector $s(\bx)$.
Write $\bta \equiv \bta(\btheta)$ and $\btas \equiv \bta(\bthetas)$.
Then,
to bound the term $\sum_{i=2}^{|\mA| - 1}\, (\mu_i(\btas) - \mu_i(\bta))^2$ from below in terms of absolute deviations of the expected numbers of transitive edges under $\btas$ and $\bta$ while avoiding the trivial lower bound of $0$,
we note that it is possible to choose $\dot\bta \in \interior(\fullspace)$ such that
\beno
\norm{\bmu(\btas) - \bmu(\bta)}_2
\gte \norm{\bmu(\btas) - \bmu(\dot\bta)}_2
\ee 
and such that the expected number of transitive edges under $\dot\bta$ is not identical to the expected number of transitive edges under $\btas$.
To see that,
note that
\beno
\label{nubers}
\norm{\bmu(\bta)}_1 
\= \dsum_{i=1}^{|\mA| - 1} |\mu_i(\bta)|
\= \mu_1(\bta) + \dsum_{i=2}^{|\mA| - 1} \mu_i(\bta),
\ee
where $\mu_1(\bta)$ is the expected number of edges while $\sum_{i=2}^{|\mA| - 1} \mu_i(\bta)$ is the expected number of transitive edges under $\bta$.
Therefore,
if $\mu_1(\btas) = \mu_1(\bta)$ and $\norm{\bmu(\bta)}_1 = \norm{\bmu(\btas)}_1$,
then the expected numbers of transitive edges under $\bta$ and $\btas$ are identical.
Let
\beno
\mT
\= \{\bmu' \in \rint(\mM):\;\; \mu_1'\; =\; \mu_1(\btas),\;\; \norm{\bmu'}_1 = \norm{\bmu(\btas)}_1\}
\ee
be the set of mean-value parameter vectors for which the expected numbers of edges and the expected numbers of transitive edges are identical.
To see that it is possible to choose $\dot\bta \in \interior(\fullspace)$ such that $\mu_1(\dot\bta) = \mu_1(\btas)$ and $\norm{\bmu(\dot\bta)}_1 \neq \norm{\bmu(\btas)}_1$,
note that the set $\mT$ is a subset of the boundary of the $\ell_1$-ball with center $\bm{0} \in \mR^{|\mA|-1}$ and radius $\norm{\bmu(\btas)}_1$.
Suppose that the $\ell_2$-distance of $\bmu(\bta)$ from $\bmu(\btas)$ is equal to $\rho_1 > 0$;
note that the $\ell_2$-distance is strictly positive because $\bta \neq \btas$ implies $\bmu(\bta) \neq \bmu(\btas)$.
Observe that the $\ell_2$-ball $\mB(\bmu(\btas),\, \rho_1)$ with center $\bmu(\btas) \in \rint(\mM)$ and radius $\rho_1 > 0$ need not be contained in $\rint(\mM)$,
but,
owing to the fact that the set $\mM$ is convex by construction,
it is possible to construct a smaller $\ell_2$-ball $\mB(\bmu(\btas),\, \rho_2)$ with the same center $\bmu(\btas) \in \rint(\mM)$ but smaller radius $0 < \rho_2 < \rho_1$ such that the resulting $\ell_2$-ball $\mB(\bmu(\btas),\, \rho_2)$ is contained in $\rint(\mM)$.
The $\ell_1$-ball with center $\bm{0} \in \mR^{|\mA|-1}$ and radius $\norm{\bmu(\btas)}_1$ and the $\ell_2$-ball $\mB(\bmu(\btas),\, \rho_2) \subseteq \rint(\mM)$ with center $\bmu(\btas) \in \rint(\mM)$ and radius $\rho_2$ intersect,
but it is not too hard to see that it is possible to choose $\dot\bmu \in \rint(\mM)$ inside the $\ell_2$-ball $\mB(\bmu(\btas),\, \rho_2) \subseteq \rint(\mM)$ such that $\dot\mu_1 = \mu_1(\btas)$ and $\norm{\dot\bmu}_1 \neq \norm{\bmu(\btas)}_1$.
In addition,
due to the fact that the map $\bmu: \interior(\fullspace) \mapsto \rint(\mM)$ is a homeomorphism \citep[][Theorem 3.6, p.\ 74]{Br86},
there exists $\dot\bta \in \interior(\fullspace)$ such that $\bmu(\dot\bta) = \dot\bmu \in \mB(\bmu(\btas),\, \rho_2) \subseteq \rint(\mM)$,
which implies that $\mu_1(\dot\bta) = \mu_1(\btas)$ and $\norm{\bmu(\dot\bta)}_1 \neq \norm{\bmu(\btas)}_1$.
In conclusion,
there exists $\dot\bta \in \interior(\fullspace)$ such that $\mu_1(\dot\bta) = \mu_1(\btas)$ and $\norm{\bmu(\dot\bta)}_1 \neq \norm{\bmu(\btas)}_1$ and
\beno
\norm{\bmu(\btas) - \bmu(\bta)}_2
\gte \norm{\bmu(\btas) - \bmu(\dot\bta)}_2
\= \dsum_{i=2}^{|\mA| - 1} (\mu_i(\btas) - \mu_i(\dot\bta))^2
&>& 0.
\ee
As a consequence,
we can bound $\norm{\bmu(\btas) - \bmu(\dot\bta)}_2$ from below in terms of absolute deviations of the expected numbers of transitive edges under $\btas$ and $\dot\bta$ while avoiding the trivial lower bound of $0$.

To do so,
we first use the Cauchy-Schwarz inequality to obtain
\beno
&& \norm{\bmu(\btas) - \bmu(\dot\bta)}_2
\;=\; \sqrt{(\mu_1(\btas) - \mu_1(\dot\bta))^2 + \dsum_{i=2}^{|\mA| - 1} (\mu_i(\btas) - \mu_i(\dot\bta))^2}\s
\\
\gte \sqrt{\dsum_{i=2}^{|\mA| - 1} (\mu_i(\btas) - \mu_i(\dot\bta))^2}
\;\geq\; \dfrac{1}{\sqrt{|\mA|-2}}\; \dsum_{i=2}^{|\mA| - 1} \left|\mu_i(\btas) - \mu_i(\dot\bta))\right|
\ee
and then use the triangle inequality to obtain
\beno
\dsum_{i=2}^{|\mA| - 1} \left|\mu_i(\btas) - \mu_i(\dot\bta)\right|
\gte \left|\dsum_{i=2}^{|\mA| - 1} (\mu_i(\btas) - \mu_i(\dot\bta))\right|\s
\\
\= \left|\dsum_{i=2}^{|\mA| - 1} \mu_i(\btas) - \dsum_{i=2}^{|\mA| - 1} \mu_i(\dot\bta)\right|.
\ee
The sums $\sum_{i=2}^{|\mA| - 1} \mu_i(\btas)$ and $\sum_{i=2}^{|\mA| - 1} \mu_i(\dot\bta)$ are equal to the expected number of transitive edges under $\btas$ and $\dot\bta$,
respectively,
as shown above.
Lemma \ref{tran_edge_alpha_exp} can hence be applied to show that,
for all $\bta \in \interior(\fullspace)$,
there exists $C_2(\bta) > 0$ such that
\beno
\dsum_{i=2}^{|\mA|-1} \mu_i(\bta)
\,=\, \mbE_{\bta}\; \dsum_{a\in\mA\, <\, b\in\mA}\; X_{a,b}\; \max\limits_{c \in \mA,\; c\, \neq\, a, b}\,  X_{a,c}\, X_{b,c}
\,=\, C_2(\bta)\, \dis{|\mA| \choose 2},
\ee
where $C_2(\bta)$ satisfies
\beno
\lambda\; \Lambda(\bta)^3 + (1 - \lambda)\; \Omega(\bta)^3
\;\leq\; C_2(\btheta)
\;\leq\; \lambda\; \Lambda(\bta) + (1 - \lambda)\; \Omega(\bta),
\ee
provided $|\mA| \geq 3$.
Collecting results shows that
\beno
&& \norm{\bmu(\btas) - \bmu(\dot\bta)}_2
\;=\; \sqrt{(\mu_1(\btas) - \mu_1(\dot\bta))^2 + \dsum_{i=2}^{|\mA| - 1} (\mu_i(\btas) - \mu_i(\dot\bta))^2}\s
\\
\hide{
\gte \dfrac{1}{\sqrt{|\mA|-2}}\; \dsum_{i=2}^{|\mA| - 1} \left|\mu_i(\btas) - \mu_i(\dot\bta))\right|\s
\\
}
\gte \dfrac{1}{\sqrt{|\mA|-2}}\; \left|\dsum_{i=2}^{|\mA| - 1} \mu_i(\bta(\bthetas)) - \dsum_{i=2}^{|\mA| - 1} \mu_i(\bta(\dot\btheta))\right|\s
\\
\= \dfrac{1}{\sqrt{|\mA|-2}}\; \left|C_2(\btas) - C_2(\dot\bta)\right|\; \dis{|\mA| \choose 2}\s
\\
\gte \dfrac{1}{|\mA|^{1/4}\, (|\mA| - 1)^{1/4}}\; \left|C_2(\btas) - C_2(\dot\bta)\right|\; \dis{|\mA| \choose 2}\s
\\
\= \dfrac{1}{2^{1/4}}\; \left|C_2(\btas) - C_2(\dot\bta)\right|\; \dis{|\mA| \choose 2}^{3/4},
\ee
where $|C_2(\btas) - C_2(\dot\bta)| > 0$ by the choice of $\dot\bta$,
as explained above.

We can therefore conclude that,
for all $\btheta \in \bTheta \setminus \mB(\bthetas, \epsilon)$,
there exists a function $C(\bthetas, \btheta) > 0$ such that
\beno
\norm{\bmu(\bta(\bthetas)) - \bmu(\bta(\btheta))}_2
\gte C(\bthetas, \btheta)\, \dis{|\mA| \choose 2}^{3/4}.
\ee

\hide{
Last,
but not least,
we show that there exists $\delta(\epsilon) > 0$ such that,
for all $\bta(\btheta) \in \fullspace \setminus \mB(\bta(\bthetas), \gamma(\epsilon))$,
we have $\norm{\bmu(\bta(\bthetas)) - \bmu(\bta(\btheta))}_2 \geq \delta(\epsilon)\, {|\mA| \choose 2}^{3/4} > 0$.
To do so,
it is convenient to operate in the natural parameter space $\fullspace = \mR^{|\mA| - 1}$ of the full exponential family with natural parameter vector $\bta$ and sufficient statistic vector $s(\bx)$.
Write $\bta \equiv \bta(\btheta)$ and $\btas \equiv \bta(\bthetas)$.
By classic exponential-family theory, 
$\psi(\bta)$ is a strictly convex function of the natural parameter vector $\bta$ on the interior $\interior(\fullspace)$ of the natural parameter space $\fullspace$,
which is a convex set \citep[][Theorem 1.13, p.\ 19]{Br86}.
Thus,
the gradient $\nabla_{\bta}\, \psi(\bta) = \bmu(\bta)$ is a strictly increasing function of the natural parameter vector $\bta$ and $\norm{\bmu(\bta^\star) - \bmu(\bta)}_2$,
considered as a function of $\bta$ for fixed $\bta^\star \in \interior(\fullspace)$,
attains its minimum on the boundary of the ball $\mB(\bta^\star,\, \gamma(\epsilon))$.
Therefore,
$\norm{\bmu(\bta^\star) - \bmu(\bta)}_2$ is bounded away from $0$ by a function of $\epsilon > 0$ for all $\bta \in \fullspace \setminus \mB(\bta^\star,\, \gamma(\epsilon))$.
}
Last,
but not least,
we show that there exists $\delta(\epsilon) > 0$ such that,
for all $\bta(\btheta) \in \fullspace \setminus \mB(\bta(\bthetas), \gamma(\epsilon))$,
we have $\norm{\bmu(\bta(\bthetas)) - \bmu(\bta(\btheta))}_2 \geq \delta(\epsilon)\, {|\mA| \choose 2}^{3/4} > 0$.
To do so,
it is convenient to operate in the natural parameter space $\fullspace = \mR^{|\mA| - 1}$ of the full exponential family with natural parameter vector $\bta$ and sufficient statistic vector $s(\bx)$.
Write $\bta \equiv \bta(\btheta)$ and $\btas \equiv \bta(\bthetas)$.
By classic exponential family theory,
the map $\bmu : \interior(\fullspace) \mapsto \rint(\mbM)$ is a homeomorphism,
i.e.,
$\bmu: \interior(\fullspace) \mapsto \rint(\mbM)$ is one-to-one and continuous,
and so is its inverse $\bmu^{-1}: \rint(\mbM) \mapsto \interior(\fullspace)$ \citep[][Theorem 3.6, p.\ 74]{Br86}.
By the continuity of $\bmu^{-1}: \rint(\mbM) \mapsto \interior(\fullspace)$,
we know that,
for each $\epsilon > 0$ and each $\gamma(\epsilon) > 0$,
there exists $a(\epsilon) > 0$ such that $\bmu^{-1}(\mB(\bmu(\btas),\, a(\epsilon)))\, \subseteq\, \mB(\btas,\, \gamma(\epsilon))$.
In other words,
all elements of\linebreak
$\mB(\bmu(\btas),\, a(\epsilon))$ map to elements of $\mB(\btas,\, \gamma(\epsilon))$ and thus no element of $\fullspace \setminus \mB(\btas,\, \gamma(\epsilon))$ can map to an element of $\mB(\bmu(\btas),\, a(\epsilon))$.
In addition,
we have shown above that $\norm{\bmu(\bta(\bthetas)) - \bmu(\bta(\btheta))}_2 \geq C(\bthetas, \btheta)\, {|\mA| \choose 2}^{3/4}$.
Combining these results shows that there exists $\delta(\epsilon) > 0$ such that $a(\epsilon)$ can be written as $a(\epsilon) = \delta(\epsilon)\, {|\mA| \choose 2}^{3/4}$.
In summary,
no element of $\bTheta\, \setminus\, \mB(\bthetas,\, \epsilon)$ can map to an element of $\mB(\bmu(\btheta^\star),\, \delta(\epsilon)\, {|\mA| \choose 2}^{3/4})$.

Taken together,
these results show that there exists $\epsilon > 0$ such that
for all $\btheta \in \bTheta \setminus \mB(\bthetas, \epsilon)$,
\beno
\norm{\bmu(\bta(\bthetas)) - \bmu(\bta(\btheta))}_2
\gte \delta(\epsilon)\, \dis{|\mA| \choose 2}^{3/4},
\ee
provided $|\mA| \geq 4$.
Therefore,
identifiability condition \eqref{nonfullid} of Theorem \ref{nonfull} is satisfied with $\alpha = 3/4$ provided $|\mA| \geq 4$.

\hide{

\begin{lemma} 
\label{lemma:gw_grad_ident}
Consider the curved exponential family with within-neighborhood edge and geometrically weighted edgewise shared partner terms as defined in Section \ref{sec:curved}.
Let
\beno
\eta_1(\btheta)
\= \theta_1\s
\\
\eta_{i+1}(\btheta)
\= \theta_2\; \theta_3 \, \left[1 - \left(1 - \dfrac{1}{\theta_3}\right)^i\right],
&& i = 1, \dots, \ppp-1.
\ee
If $\bTheta = \mbR\, \times\, (0, \infty)\, \times\, (1, \infty)$,
then $\bta: \interior(\bTheta) \mapsto \interior(\fullspace)$ is one-to-one and continuously differentiable on $\interior(\bTheta)$ with gradient $\mJ$ of rank $\qqq$,
where the partial derivatives $\nabla_{\theta_i}\, \eta_j(\btheta)$ satisfy
\be
\label{partialderivatives}
\nabla_{\theta_1}\, \eta_1(\btheta)
\;=\; 1,\;\;\;\;
\nabla_{\theta_1}\, \eta_{2}(\btheta)
\;=\;0,\;\;\;\;
\nabla_{\theta_1}\, \eta_{i+1}(\btheta)
\;=\; 0,\s\s
\\
\nabla_{\theta_2}\, \eta_1(\btheta)
\;=\; 0,\;\;\;\;
\nabla_{\theta_2}\, \eta_2(\btheta)
\;=\; 1,\;\;\;\;
1
\;<\; \nabla_{\theta_2}\, \eta_{i+1}(\btheta)
\;<\; \theta_3,\s
\\
\nabla_{\theta_3}\, \eta_1(\btheta)
\;=\; 0,\;\;\;\;
\nabla_{\theta_3}\, \eta_2(\btheta)
\;=\; 0,\;\;\;\;
\dfrac{\theta_2}{\theta_3^2}
\;\leq\; \nabla_{\theta_3}\, \eta_{i+1}(\btheta)
\;<\; \theta_2,
\ee
where $i = 2, \dots, \ppp-1$.
If there exist $C_1 > 0$ and $C_2 > 1$ such that $\bTheta = \mR \times (0, C_1) \times (1, C_2)$,
then there exists $C > 0$ such that $|\nabla_{\theta_i}\, \eta_{j}(\btheta)| \leq C$ for all $\btheta \in \interior(\bTheta)$ ($i = 1, \dots, \qqq$, $j = 1, \dots, \ppp$).
\hide{
\beno
\min\left(1, \dfrac{C_1}{C_3^2}\right)
\,\leq\, \nabla_{\theta_i}\, \eta_{3}(\btheta) \leq \dots \leq \nabla_{\theta_i}\, \eta_{\ppp}(\btheta)
\,\leq\, \max(C_2, C_3),
& i = 2, 3
\ee
and 
}
\end{lemma}

\llproof \ref{lemma:gw_grad_ident}.
It is not too hard to see that  $\bta: \interior(\bTheta) \mapsto \interior(\fullspace)$ is one-to-one and continuously differentiable on $\interior(\bTheta)$ with gradient $\mJ$ of rank $\qqq$ and partial derivatives
\beno
\nabla_{\theta_1} \, \eta_1(\btheta)
\= 1,\;\;\;\;
\nabla_{\theta_1}\, \eta_{i+1}(\btheta)
\= 0,\s\s
\\
\nabla_{\theta_2} \, \eta_1(\btheta)
\= 0,\;\;\;\;
\nabla_{\theta_2} \, \eta_{i+1}(\btheta)
\= \theta_3 \left(1 - \vartheta(\theta_3)^i\right),\s
\\
\nabla_{\theta_3}\, \eta_1(\btheta)
\= 0,\;\;\;\;
\nabla_{\theta_3} \, \eta_{i+1}(\btheta)
\= \theta_2 \left(1 - \vartheta(\theta_3)^i - \dfrac{i}{\theta_3}\, \vartheta(\theta_3)^{i-1} \right),
\ee
where $\vartheta(\theta_3) = 1 - 1\, /\, \theta_3$ and $i = 1, \dots, \ppp-1$.
We bound the partial derivatives $\nabla_{\theta_2} \, \eta_{i+1}(\btheta)$ and $\nabla_{\theta_3} \, \eta_{i+1}(\btheta)$ ($i = 1, \dots, \ppp-1$) below.

\s

{\em Partial derivatives $\nabla_{\theta_2} \, \eta_{i+1}(\btheta)$ ($i = 1, \dots, \ppp-1$).}
The partial derivative $\nabla_{\theta_2} \, \eta_{2}(\btheta) = 1$ is constant on $\interior(\bTheta)$.
To bound the partial derivatives $\nabla_{\theta_2} \, \eta_{i+1}(\btheta)$ ($i = 2, \dots, \ppp-1$) from below,
observe that,
for all $\theta_3 \in (1, \infty)$,
\beno
\nabla_{\theta_2} \, \eta_{3}(\btheta) &<& \dots &<& \nabla_{\theta_2} \, \eta_{\ppp}(\btheta),
\ee
which implies that,
for all $\theta_3 \in (1, \infty)$,
\beno
\nabla_{\theta_2} \, \eta_{i+1}(\btheta)
\gte \nabla_{\theta_2} \, \eta_{3}(\btheta)
\= 2 - \dfrac{1}{\theta_3}
&>& 1,
& i = 2, \dots, \ppp-1.
\ee
To bound the partial derivatives $\nabla_{\theta_2} \, \eta_{i+1}(\btheta)$ ($i = 2, \dots, \ppp-1$) from above,
note that
\beno
\nabla_{\theta_2}\, \eta_{i+1}(\btheta)
\= \theta_3 \, \left(1 - \vartheta(\theta_3)^i\right)
&<& \theta_3,
& i = 2, \dots, \ppp-1,
\ee
because $0 < \vartheta(\theta_3) < 1$ for all $\theta_3 \in (1, \infty)$.
\hide{
In conclusion,
$\nabla_{\theta_2} \, \eta_{2}(\btheta) = 1$ is constant on $\interior(\bTheta)$ while
\beno
1
&<& \nabla_{\theta_2} \, \eta_{i+1}(\btheta)
&<& \theta_3,
&& i = 2, \dots, \ppp-1.
\ee
}

\s

{\em Partial derivatives $\nabla_{\theta_3} \, \eta_{i+1}(\btheta)$ ($i = 1, \dots, \ppp-1$).}
The partial derivative $\nabla_{\theta_3} \, \eta_{2}(\btheta) = 0$ is constant on $\interior(\bTheta)$.
To bound the partial derivatives $\nabla_{\theta_3} \, \eta_{i+1}(\btheta)$ ($i = 2, \dots, \ppp-1$) from below,
note that
\beno
\nonumber
\nabla_{\theta_3}\, \eta_{i+2}(\btheta) - \nabla_{\theta_3}\, \eta_{i+1}(\btheta)
\;=\; \theta_2\, \dfrac{i\; (\theta_3 - 1)^{i-1}}{\theta_3^{i+1}}
\;>\; 0,
& i = 1, \dots, \ppp-1
\ee
provided $\theta_2 \in (0, \infty)$ and $\theta_3 \in (1, \infty)$.
As a result,
\beno
\nabla_{\theta_3}\, \eta_{3}(\btheta) &<& \dots &<& \nabla_{\theta_3} \, \eta_{\ppp}(\btheta),
\ee
which implies that
\beno
\nabla_{\theta_3}\, \eta_{i+1}(\btheta)
\gte \nabla_{\theta_3}\, \eta_{3}(\btheta)
\= \dfrac{\theta_2}{\theta_3^2},
& i = 2, \dots, \ppp-1.
\ee
To bound the partial derivatives $\nabla_{\theta_3} \, \eta_{i+1}(\btheta)$ ($i = 2, \dots, \ppp-1$) from above,
note that $\vartheta(\theta_3) \geq 0$ for all $\theta_3 \in (1, \infty)$,
which implies that
\beno
\nabla_{\theta_3}\, \eta_{i+1}(\btheta)
&<& \theta_2,
& i = 2, \dots, \ppp-1.
\ee
\hide{
In conclusion,
$\nabla_{\theta_3} \, \eta_{2}(\btheta) = 0$ is constant on $\interior(\bTheta)$ while
\beno
\dfrac{\theta_2}{\theta_3^2}
\lte \nabla_{\theta_3}\, \eta_{i+1}(\btheta)
&<& \theta_2,
& i = 2, \dots, \ppp-1.
\ee
}

{\em Boundedness of partial derivatives.}
By inspecting \eqref{partialderivatives},
it is evident that there exists $C > 0$ such that $|\nabla_{\theta_i}\, \eta_{j}(\btheta)| \leq C$ for all $\btheta \in \interior(\bTheta) = \mR \times (0, C_1) \times (1, C_2)$ ($i = 1, 2, 3$, $j = 1, \dots, \ppp$). 

}
\s

\begin{lemma}
\label{edge.transitive1}
Consider an exponential-family random graph with within-neighborhood edge and transitive edge terms as defined in Appendix \ref{sec:canonical}.
Let $\bTheta = \mR \times \mR^+$ and assume that $\bthetas \in \interior(\bTheta)$.
Then,
for all $\epsilon > 0$ small enough so that $\mathscr{B}(\bthetas, \epsilon) \subseteq \interior(\bTheta)$,
there exists $\delta(\epsilon) > 0$ such that,
for all $\btheta \in \bTheta \setminus \mathscr{B}(\bthetas, \, \epsilon)$,
\beno
\norm{\bmu(\bta(\bthetas)) - \bmu(\bta(\btheta))}_2
\gte \delta(\epsilon)\, \dsum_{k=1}^K {|\mA_k| \choose 2},
\ee
provided $|\mA_k| \geq 3$ ($k = 1, \dots, K$).
Therefore,
identifiability condition \eqref{nonfullid} of Theorem \ref{nonfull} is satisfied with $\alpha = 1$ provided $|\mA_k| \geq 3$ ($k = 1, \dots, K$). 
\end{lemma}

\llproof \ref{edge.transitive1}.
To ease the presentation,
we consider $K = 1$ neighborhood and drop the subscript $k$ from all neighborhood-dependent quantities.
The extension to $K \geq 2$ neighborhoods is straightforward.
In the following,
we write $\bmu(\bta(\btheta)) = \bmu(\btheta) = \mbE_{\btheta}\, s(\bX)$ and $\psi(\bta(\btheta)) = \psi(\btheta)$,
because the natural parameter vector of the exponential family is given by $\bta(\btheta) = \btheta$.

To verify identifiability condition \eqref{nonfullid} of Theorem \ref{nonfull},
we need to verify that,
for all $\epsilon > 0$ small enough so that $\mathscr{B}(\bthetas, \epsilon) \subseteq \interior(\bTheta)$,
there exists $\delta(\epsilon) > 0$ such that,
for all $\btheta \in \bTheta \setminus \mathscr{B}(\bthetas, \epsilon)$,
\beno
\norm{\bmu(\bthetas) - \bmu(\btheta)}_2
&\geq& \delta(\epsilon)\, \dis\binom{|\mA|}{2}^\alpha
&\mbox{for some}
& 0 \leq \alpha \leq 1.
\ee
To do so,
observe that
\be 
\nonumber 
\norm{\bmu(\bthetas) - \bmu(\btheta)}_2
&=&\sqrt{(\mu_1(\bthetas) - \mu_1(\btheta))^2 +  (\mu_2(\bthetas) - \mu_2(\btheta))^2}.
\ee 
The deviation $|\mu_1(\bthetas) - \mu_1(\btheta)|$ can be dealt with by using Lemma \ref{mean.bound}, 
which shows that there exists $C_1(\btheta) \in [\Lambda(\btheta),\, \Omega(\btheta)]$ such that
\be 
\nonumber 
\mu_1(\btheta)
&=& \mbE_{\btheta} \dsum_{i\in\mA\, <\, j\in\mA} X_{i,j}
&=& C_1(\btheta)\, \dis\binom{|\mA|}{2},
\ee 
where $\Lambda(\btheta)$ and $\Omega(\btheta)$ are defined in the introduction of Appendix \ref{m.id}.
Here, 
$\Lambda(\btheta)$ is given by
$\Lambda(\btheta) = [1 + \exp(-\theta_1)]^{-1} > 0$,
while $\Omega(\btheta)$ satisfies $0 < \Lambda(\btheta) \leq \Omega(\btheta) < 1$ for all $\btheta \in \mR \times \mR^+$.
Therefore, 
the deviation $|\mu_1(\bthetas) - \mu_1(\btheta)|$ satisfies
\be 
\nonumber 
|\mu_1(\bthetas) - \mu_1(\btheta)|
&=& |C_1(\bthetas) - C_1(\btheta)| \, \dis\binom{|\mA|}{2}.
\ee
To deal with the deviation $|\mu_2(\bthetas) - \mu_2(\btheta)|$,
observe that,
by Lemma \ref{tran_edge_alpha_exp},
there exists $C_2(\btheta) > 0$ such that
\be 
\nonumber 	
\mu_2(\btheta)
&=& \mbE_{\btheta}\, \dsum_{i\in\mA\, <\, j\in\mA} X_{i,j}\, \max\limits_{h \in \mA,\, h \neq i,j} X_{i,h} \; X_{j,h}
&=& C_2(\btheta) \, \dis\binom{|\mA|}{2},
\ee
where $C_2(\btheta)$ satisfies 
\be
\nonumber
\lambda\, \Lambda(\btheta)^3 + (1 - \lambda)\, \Omega(\btheta)^3
\,\leq\, C_2(\btheta)
\,\leq\, \lambda\, \Lambda(\btheta) + (1 - \lambda)\, \Omega(\btheta),
& 0 < \lambda < 1.
\ee
As a result,
\be 
\nonumber 
|\mu_2(\bthetas) - \mu_2(\btheta)|
&=& |C_2(\bthetas) - C_2(\btheta)|\, \dis\binom{|\mA|}{2}.
\ee  

\hide{

We note that,
to grow with the number of within-neighborhood edge variables ${|\mA| \choose 2}$,
$C_i(\btheta)$ ($i = 1, 2$) must be bounded away from $0$.
To ensure that,
$\Lambda(\btheta)$ must be bounded away from $0$,
which is guaranteed as long as
\bi
\item $\btheta \in \mR \times \mR^-$ and the neighborhoods are bounded,
in which case there exist $A_1 > 0$ and $A_2 > 0$ such that $3 \leq A_1 \leq |\mA| \leq A_2$ and hence
\beno
\Lambda(\btheta) 
\gte [1 + \exp(-\theta_1 - \theta_2\, (2\, A_2 - 3))]^{-1};
\ee
\item $\btheta \in \mR \times \mR^+$ and the neighborhoods are of the same order of magnitude,
in which case
\beno
\Lambda(\btheta) \= [1 + \exp(-\theta_1)]^{-1};
\ee
\ei
note that $\Omega(\btheta)$ satisfies $\Lambda(\btheta) \leq \Omega(\btheta)$ for all $\btheta \in \mR^2$.

}

We can hence conclude that,
for all $\btheta \in \bTheta \setminus \mathscr{B}(\bthetas, \, \epsilon)$,
there exists a function $C(\btheta, \bthetas) > 0$ such that
\beno
\norm{\bmu(\bthetas) - \bmu(\btheta)}_2
\= C(\bthetas, \btheta)\, \dis\binom{|\mA|}{2}
&>& 0,
\ee
where strict positivity follows from the fact that the map $\bmu: \interior(\bTheta) \mapsto \rint(\mM)$ is one-to-one by classic exponential-family theory \citep[][Theorem 3.6, p.\ 74]{Br86},
which implies that $\norm{\bmu(\bthetas) - \bmu(\btheta)}_2 > 0$.

\hide{
Last,
but not least,
we show that,
for all $\epsilon > 0$ small enough so that $\mathscr{B}(\bthetas, \epsilon) \subseteq \interior(\bTheta)$, 
there exists $\delta(\epsilon) > 0$ such that, 
for all $\btheta \in \bTheta \setminus \mB(\bthetas, \epsilon)$,
we have $\norm{\bmu(\bthetas) - \bmu(\btheta)}_2 \geq \delta(\epsilon)\, \binom{|\mA|}{2} > 0$.
To do so,
note that $\bta(\btheta) = \btheta$ is the natural parameter vector of the exponential family.
Thus,
by classic exponential-family theory,
$\psi(\btheta)$ is a strictly convex function of the natural parameter vector $\btheta$ on the interior of the natural parameter space $\{\btheta \in \mR^{\qqq}: \psi(\btheta) < \infty\} = \mR^2$ \citep[][Theorem 1.13, p.\ 19]{Br86}.
Therefore,
the gradient $\nabla_{\btheta}\, \psi(\btheta) = \bmu(\btheta)$ is a strictly increasing function of $\btheta$ and the smallest value of $\norm{\bmu(\bthetas) - \bmu(\btheta)}_2$,
considered as a function of $\btheta$ for fixed $\bthetas \in \interior(\bTheta)$, 
is attained on the boundary of ball $\mB(\bthetas, \epsilon) \subseteq \interior(\bTheta)$. 
}
Last,
but not least,
we show that,
for all $\epsilon > 0$ small enough so that $\mathscr{B}(\bthetas, \epsilon) \subseteq \interior(\bTheta)$,
there exists $\delta(\epsilon) > 0$ such that,
for all $\btheta \in \bTheta \setminus \mB(\bthetas, \epsilon)$,
we have $\norm{\bmu(\bthetas) - \bmu(\btheta)}_2 \geq \delta(\epsilon)\, \binom{|\mA|}{2} > 0$.
To do so,
note that $\bta(\btheta) = \btheta$ is the natural parameter vector of the exponential family.
By classic exponential family theory,
the map $\bmu : \interior(\bTheta) \mapsto \rint(\mbM)$ is a homeomorphism,
i.e.,
$\bmu: \interior(\bTheta) \mapsto \rint(\mbM)$ is one-to-one and continuous,
and so is its inverse $\bmu^{-1}: \rint(\mbM) \mapsto \interior(\bTheta)$ \citep[][Theorem 3.6, p.\ 74]{Br86}.
By the continuity of $\bmu^{-1}: \rint(\mbM) \mapsto \interior(\bTheta)$,
we know that,
for each $\epsilon > 0$,
there exists $a(\epsilon) > 0$ such that $\bmu^{-1}(\mB(\bmu(\btheta^\star),\, a(\epsilon)))\, \subseteq\, \mB(\bthetas,\, \epsilon)$.
Thus,
all elements of $\mB(\bmu(\btheta^\star),\, a(\epsilon))$ map to elements of $\mB(\bthetas,\, \epsilon)$ and thus no element of $\bTheta \setminus \mB(\bthetas,\, \epsilon)$ can map to an element of $\mB(\bmu(\btheta^\star),\, a(\epsilon))$.
In addition,
we have shown above that $\norm{\bmu(\bthetas) - \bmu(\btheta)}_2 = C(\bthetas, \btheta)\, \binom{|\mA|}{2}$.
Combining these results shows that there exists $\delta(\epsilon) > 0$ such that $a(\epsilon)$ can be written as $a(\epsilon) = \delta(\epsilon)\, \binom{|\mA|}{2}$.
In summary,
no element of $\bTheta \setminus \mB(\bthetas,\, \epsilon)$ can map to an element of $\mB(\bmu(\btheta^\star),\, \delta(\epsilon)\, \binom{|\mA|}{2})$.

As a result,
for all $\epsilon > 0$ small enough so that $\mathscr{B}(\bthetas, \epsilon) \subseteq \interior(\bTheta)$,
there exists $\delta(\epsilon) > 0$ such that,
for all $\btheta \in \bTheta \setminus \mathscr{B}(\bthetas, \, \epsilon)$,
\be 
\nonumber
\norm{\bmu(\bthetas) - \bmu(\btheta)}_2
\gte \delta(\epsilon)\, \dis\binom{|\mA|}{2},
\ee 
provided $|\mA| \geq 3$.
Therefore,
identifiability condition \eqref{nonfullid} of Theorem \ref{nonfull} is satisfied with $\alpha = 1$ provided $|\mA| \geq 3$. 

\subsection{Auxiliary lemmas: Theorem \ref{theorem:parameters.correct}}
\label{appendix.m.id}

We prove an auxiliary lemma that is useful for verifying the conditions of Theorem \ref{theorem:parameters.correct}.

\s

\begin{lemma}
\label{edge.transitive.attribute}
Consider an exponential-family random graph with within-neighborhood edge, transitive edge, and same-attribute edge terms as defined in Section \ref{mainresults}.
Then identifiability condition [C.3] of Theorem \ref{theorem:parameters.correct} is satisfied with $\alpha = 1$ provided $|\mA_k| \geq 3$ ($k = 1, \dots, K$).
\end{lemma}

\llproof \ref{edge.transitive.attribute}.
To verify identifiability condition [C.3] of Theorem \ref{theorem:parameters.correct},
we need to verify that,
for all $\epsilon > 0$ small enough so that $\mathscr{B}(\bthetad, \epsilon) \subseteq \bThetad$,
there exist $\delta(\epsilon) > 0$ such that, 
for all $\btheta \in \bThetad \setminus \mathscr{B}(\bthetad, \epsilon)$,
\beno
\norm{\mbE_{\bthetas}\, b(\bX) - \mbE_{\btheta}\, b(\bX)}_2
&\geq& \delta(\epsilon)\, \dsum_{k=1}^K \dis{|\mA_k| \choose 2}^\alpha
\;\mbox{ for some }\; 0 \leq \alpha \leq 1,
\ee
where $\mbE_{\bthetas}\, b(\bX)$ is the expectation of $b(\bX)$ under the data-generating exponential-family distribution with natural parameter vector $\bthetas \in \bThetas$ and $\mbE_{\btheta}\, b(\bX)$ is the expectation of $b(\bX)$ under the misspecified exponential-family distribution with natural parameter vector $\btheta \in \bThetad$.

\hide{

\s

{\bf Expectation $\mbE_{\bthetas}\, b(\bX)$ under $\bthetas \in \interior(\bThetas)$.}
By using an argument along the lines of Lemma \ref{edge.transitive1},
we can show that there exist $C_i(\bthetas) > 0$ ($i = 1, 2$) such that
\be
\nonumber
\mbE_{\bthetas}\, b_i(\bX)
&=& C_i(\bthetas)\, \dis\binom{|\mA|}{2},
& i = 1, 2.
\ee
We note that,
to show that $C_i(\bthetas) > 0$ ($i = 1, 2$),
we used the fact that
\bi
\item Case $\bthetas \in \interior(\bThetas) = \mbR \times \mbR^+ \times \mR^-$:
$\Lambda(\bthetas) \geq [1 + \exp(-\thetas_1 - \thetas_3)]^{-1} > 0$;
\item Case $\bthetas \in \interior(\bThetas) = \mbR \times \mbR^+ \times \mR^+$:
$\Lambda(\bthetas) \geq [1 + \exp(-\thetas_1)]^{-1} > 0$;
\ei
while $\Omega(\bthetas)$ satisfies $\Lambda(\bthetas) \leq \Omega(\bthetas)$ for all $\bthetas \in \interior(\bThetas) = \mbR \times \mbR^+ \times \mR$.
In other words,
in both cases,
the full conditional probabilities are bounded away from $0$,
and so are $C_i(\bthetas) > 0$ ($i = 1, 2$) by using an argument along the lines of Lemma \ref{edge.transitive1}.

\s

{\bf Expectation $\mbE_{\btheta}\, b(\bX)$ under $\btheta \in \interior(\bThetad)$.}
The expectation $\mbE_{\btheta}\, b(\bX)$ is the expectation of $b(\bX)$ under the exponential-family distribution with natural parameter vector $\btheta = (\theta_1, \theta_2) \in \bThetad = \mR^2$ and sufficient statistic vector $b(\bx) = (s_1(\bx), s_2(\bx)) \in \mathbb{B}$.
Therefore,
by an argument along the lines of Lemma \ref{edge.transitive1},
it can be shown that there exist $C_i(\btheta) > 0$ ($i = 1, 2$) such that
\be
\nonumber
\mbE_{\btheta}\, b_i(\bX)
&=& C_i(\btheta)\, \dis\binom{|\mA|}{2},
& i = 1, 2.
\ee

{\bf Conclusion.}
By using an argument along the lines of Lemma \ref{edge.transitive1},
we conclude that,

}

Since the expectations $\mbE_{\bthetas}\, b(\bX)$ and $\mbE_{\btheta}\, b(\bX)$ are expectations of the number of within-neighborhood edges and transitive edges,
an argument along the lines of Lemma \ref{edge.transitive1} shows that,
for all $\epsilon > 0$ small enough so that $\mB(\bthetad, \epsilon) \subseteq \bThetad$,
there exist $\delta(\epsilon) > 0$ such that,
for all $\btheta \in \bThetad \setminus \mathscr{B}(\bthetas, \, \epsilon)$,
\be 
\nonumber
\norm{\mbE_{\bthetas}\, b(\bX) - \mbE_{\btheta}\, b(\bX)}_2 
\gte \delta(\epsilon)\, \dsum_{k=1}^K \dis\binom{|\mA_k|}{2},
\ee 
provided $|\mA_k| \geq 3$ ($k = 1, \dots, K$).
Therefore, 
identifiability condition [C.3] of Theorem \ref{theorem:parameters.correct} is satisfied with $\alpha = 1$ provided $|\mA_k| \geq 3$ ($k = 1, \dots, K$).

\hide{

\subsection{Other auxiliary lemmas}
\label{auxiliarylemmas}

We prove other auxiliary lemmas.

Throughout,
it is convenient to work with canonical exponential-family parameterizations.
Using canonical exponential-family parameterizations,
the expected loglikelihood function can be written as
\be 
\nonumber
\label{eq:g_lower_bound} 
g(\bta;\; \bmu(\bta^\star)) 
&=& \langle\bta,\;  \bmu(\bta^\star)\rangle - \psi(\bta),
\ee 
where $\bta^\star \equiv \bta(\bthetas) \in \etaspace \subseteq \interior(\fullspace)$, 
$\bta \equiv \bta(\btheta) \in \etaspace \subseteq \interior(\fullspace)$, 
and $\bmu(\bta^\star) = \mbE_{\bta^\star}\, s(\bX) \equiv \mbE_{\bta(\bthetas)}\, s(\bX) \equiv \mbE\, s(\bX)$.

\s\s

\hide{

\begin{lemma}
\label{lemma1}
Consider a full or non-full, curved exponential family with countable support $\mbX$ and local dependence.
Suppose that $\bta^\star \in \etaspace \subseteq \interior(\fullspace)$.
Then,
for all $\bta \in \etaspace \subseteq \interior(\fullspace)$,
there exists $\dot\bta = \lambda \; \bta^\star + (1 - \lambda) \; \bta \in \interior(\fullspace)$ $(0 < \lambda < 1)$ such that
\be 
\nonumber
\label{eq:g_lower_bound} 
g(\bta^\star;\; \bmu(\bta^\star)) - g(\bta;\; \bmu(\bta^\star)) 
&=& \langle\bta^\star - \bta,\; \bmu(\bta^\star) - \bmu(\dot\bta)\rangle.
\ee 
\end{lemma} 
\llproof \ref{lemma1}.
An application of the mean-value theorem to $g(\bta^\star;\; \bmu(\bta^\star))$ shows that,
for all $\bta \in \etaspace \subseteq \interior(\fullspace)$,
there exists $\dot\bta = \lambda \; \bta^\star + (1 - \lambda) \; \bta$ $(0 < \lambda < 1)$ such that
\be 
\nonumber 
g(\bta^\star;\; \bmu(\bta^\star)) - g(\bta;\; \bmu(\bta^\star)) 
&=& \langle\bta^\star - \bta,\, \nabla_{\bta}\; g(\bta;\; \bmu(\bta^\star))\, |_{\bta=\dot\bta}\rangle.
\ee
We note that $\dot\bta = \lambda\; \bta^\star + (1 - \lambda) \; \bta \in \interior(\fullspace)$ is in the interior $\interior(\fullspace)$ of the natural parameter space $\fullspace$,
because $\bta^\star \in \interior(\fullspace)$ and $\bta \in \interior(\fullspace)$ and $\fullspace$ is convex \citep[][Theorem 1.13, p.\ 19]{Br86}.
Since $\dot\bta \in \interior(\fullspace)$,
the gradient $\nabla_{\bta}\; g(\bta;\; \bmu(\bta^\star))\, |_{\bta=\dot\bta}$ is given by $\nabla_{\bta}\; g(\bta;\; \bmu(\bta^\star))\, |_{\bta=\dot\bta} = \bmu(\bta^\star) - \mbE_{\dot\bta}\; s(\bX) = \bmu(\bta^\star) - \bmu(\dot\bta)$ \citep[][Corollary 2.3, pp.\ 35--36]{Br86},
which implies that
\be \nonumber 
g(\bta^\star;\; \bmu(\bta^\star)) - g(\bta;\; \bmu(\bta^\star)) 
&=& \langle\bta^\star - \bta,\, \bmu(\bta^\star) - \bmu(\dot\bta)\rangle.
\ee

}

\begin{lemma} 
\label{reverse.cs} 
{\bf \citet[][Theorem 4, p.\ 9]{Dr05}.}
Let $\bv_1 \in \mR^d$ and $\bv_2 \in \mR^d$.
If there exist $A > 0$ and $B > 0$ such that
\be 
\label{condition.reverse.cs}
\langle A\, \bv_2 - \bv_1,\; \bv_1 - B\, \bv_2\rangle 
&\geq& 0,
\ee
then 
\be  
\nonumber
\dfrac{1}{2} \dfrac{A + B}{\sqrt{A\, B}}\; \left|\langle\bv_1,\, \bv_2\rangle\right|
&\geq& \norm{\bv_1}_2\; \norm{\bv_2}_2.
\ee
\end{lemma}

\llproof \ref{reverse.cs}. Lemma \ref{reverse.cs} is proved by \citet[][Theorem 4, p.\ 9]{Dr05}.

\s\s

\begin{lemma}
\label{corollary.reverse.cs}
Consider a full or non-full, curved exponential family with countable support $\mbX$ and local dependence.
Suppose that $\bta^\star \in \etaspace \subseteq \interior(\fullspace)$.
Consider any $\bta \neq \bta^\star$ such that $\bta \in \etaspace \subseteq \interior(\fullspace)$ and any $\dot\bta = \lambda \; \bta^\star + (1 - \lambda) \; \bta \in \interior(\fullspace)$ $(0 < \lambda < 1)$.
Then
\be
\nonumber 
\langle\bta^\star - \bta,\; \bmu(\bta^\star) - \bmu(\dot\bta)\rangle
&\geq& \dfrac{2 \; \norm{\bmu(\bta^\star) - \bmu(\dot\bta)}_2\; \norm{\bta^\star - \bta}_2^{3 / 2}}{\norm{\bta^\star - \bta}_2 + 1}.
\ee
\end{lemma}
\llproof \ref{corollary.reverse.cs}.
We want to bound $\langle\bta^\star - \bta,\; \bmu(\bta^\star) - \bmu(\dot\bta)\rangle$ from below by using the reverse Cauchy-Schwarz inequality in Lemma \ref{reverse.cs}.
To ease the presentation,
we write $\bv_1 = \bmu(\bta^\star) - \bmu(\dot\bta)$ and $\bv_2 = \bta^\star - \bta$.

We first prove that condition \eqref{condition.reverse.cs} of Lemma \ref{reverse.cs} is satisfied by showing that there exist $A > 0$ and $B > 0$ such that 
\be 
\nonumber 
\langle A\, \bv_2 - \bv_1,\; \bv_1 - B\, \bv_2\rangle
\hide{
&=& A\, \langle \bv_2\, \bv_1\rangle - A\, B\, \langle\bv_2,\; \bv_2\rangle - \langle\bv_1,\; \bv_1\rangle + B\, \langle\bv_1,\; \bv_2\rangle
\s
\\
}
&=& (A + B)\, \langle \bv_1,\, \bv_2\rangle - \norm{\bv_1}_2^2 - A\, B\, \norm{\bv_2}_2^2
&\geq& 0,
\ee
which is equivalent to showing that that there exist $A > 0$ and $B > 0$ such that
\be
\nonumber
(A + B)\, \langle \bv_1,\, \bv_2\rangle
&\geq& \norm{\bv_1}_2^2 + A\, B\, \norm{\bv_2}_2^2.
\ee
By using the Cauchy-Schwarz inequality,
we obtain
\be
\nonumber
(A + B)\, \norm{\bv_1}_2\, \norm{\bv_2}_2
&\geq& (A + B)\, \langle \bv_1,\, \bv_2\rangle
&\geq& \norm{\bv_1}_2^2 + A\, B\, \norm{\bv_2}_2^2,
\ee
which implies that we must find $A > 0$ and $B > 0$ such that the following inequality is satisfied:
\be
\label{inequality}
A + B
&\geq& \dfrac{\norm{\bv_1}_2^2 + A\, B\, \norm{\bv_2}_2^2}{\norm{\bv_1}_2\, \norm{\bv_2}_2}
&\geq& \dfrac{\norm{\bv_1}_2^2 + A\, B\, \norm{\bv_2}_2^2}{\norm{\bv_1}_2\, (1 + \norm{\bv_2}_2)}.
\ee
We note that the denominator $\norm{\bv_1}_2\, \norm{\bv_2}_2$ of the first ratio in inequality \eqref{inequality} is strictly positive,
because $\norm{\bv_1}_2 > 0$ and $\norm{\bv_2}_2 > 0$.
First,
$\norm{\bv_2}_2 = \norm{\bta^\star - \bta}_2 > 0$, 
because $\bta^\star \neq \bta$.
Second,
to show that $\norm{\bv_1}_2 = \norm{\bmu(\bta^\star) - \bmu(\dot\bta)}_2 > 0$, 
note that $\bta^\star \in \etaspace \subseteq \interior(\fullspace)$ and $\bta \in \etaspace \subseteq \interior(\fullspace)$,
which implies that $\dot\bta = \lambda \; \bta^\star + (1 - \lambda) \; \bta$ $(0 < \lambda < 1)$ is contained in the interior $\interior(\fullspace)$ of the natural parameter space $\fullspace$ since $\fullspace$ is convex \citep[][Theorem 1.13, p.\ 19]{Br86}.
In addition,
$\dot\bta = \lambda \; \bta^\star + (1 - \lambda) \; \bta$ $(0 < \lambda < 1)$ implies that $\bta^\star \neq \dot\bta$.
Since the map $\bmu: \interior(\fullspace) \mapsto \rint(\mM)$ is one-to-one \citep[][Theorem 3.6, p.\ 74]{Br86},
$\bta^\star \neq \dot\bta$ implies $\bmu(\bta^\star) \neq \bmu(\dot\bta)$ and hence $\norm{\bv_2}_2 = \norm{\bmu(\bta^\star) - \bmu(\dot\bta)}_2 > 0$.
To show that condition \eqref{condition.reverse.cs} of Lemma \ref{reverse.cs} is satisfied,
choose $A = \norm{\bv_1}_2\; /\;  \norm{\bv_2}_2 > 0$ and $B = \norm{\bv_1}_2 > 0$.
Then,
using inequality \eqref{inequality} along with the chosen $A > 0$ and $B > 0$,
we obtain
\be
\nonumber
A + B 
\hide{
&=& \dfrac{\norm{\bv_1}_2}{\norm{\bv_2}_2} + \norm{\bv_1}_2
}
\;\geq\; \dfrac{\norm{\bv_1}_2^2 + A\, B\, \norm{\bv_2}_2^2}{\norm{\bv_1}_2\, (1 + \norm{\bv_2}_2)}
\;=\; \dfrac{\norm{\bv_1}_2^2 + \norm{\bv_1}_2^2\; \norm{\bv_2}_2}{\norm{\bv_1}_2\, (1 + \norm{\bv_2}_2)}
\;=\; \norm{\bv_1}_2
\;>\; 0.
\ee
As a result,
condition \eqref{condition.reverse.cs} of Lemma \ref{reverse.cs} is satisfied and we can invoke Lemma \ref{reverse.cs} to bound $\langle\bta^\star - \bta,\; \bmu(\bta^\star) - \bmu(\dot\bta)\rangle$ from below.
To do so,
recall that $g(\bta;\; \bmu(\bta^\star))$ is the expected loglikelihood function,
which implies that
\be
\nonumber
KL(\bta^\star;\, \bta)
\,=\, g(\bta^\star; \, \bmu(\bta^\star)) - g(\bta; \, \bmu(\bta^\star))
\,=\, \langle\bta^\star - \bta,\; \bmu(\bta^\star) - \bmu(\dot\bta)\rangle
\,>\, 0,
\ee
where $KL(\bta^\star;\, \bta)$ denotes the Kullback-Leibler divergence from $\mbP_{\bta^\star}$ to $\mbP_{\bta}$ and where $KL(\bta^\star;\, \bta) > 0$ since $\bta^\star \neq \bta$.
Thus,
by Lemma \ref{reverse.cs} along with $\bv_1 = \bmu(\bta^\star) - \bmu(\dot\bta)$ and $\bv_2 = \bta^\star - \bta$,
\be  
\nonumber 
\langle\bta^\star - \bta,\; \bmu(\bta^\star) - \bmu(\dot\bta)\rangle
&=& \left|\langle\bta^\star - \bta,\; \bmu(\bta^\star) - \bmu(\dot\bta)\rangle\right|\s
\\
&\geq& 2\, \dfrac{\sqrt{A\, B}}{A + B}\; \norm{\bmu(\bta^\star) - \bmu(\dot\bta)}_2\;\, \norm{\bta^\star - \bta}_2.
\ee
Last,
but not least,
using $A = \norm{\bv_1}_2\; /\;  \norm{\bv_2}_2 > 0$ and $B = \norm{\bv_1}_2 > 0$ along with $\bv_1 = \bmu(\bta^\star) - \bmu(\dot\bta)$ and $\bv_2 = \bta^\star - \bta$ gives 
\be 
\nonumber  
\langle\bta^\star - \bta,\; \bmu(\bta^\star) - \bmu(\dot\bta)\rangle
&\geq& \dfrac{2 \; \norm{\bmu(\bta^\star) - \bmu(\dot\bta)}_2\; \norm{\bta^\star - \bta}_2^{3 / 2}}{\norm{\bta^\star - \bta}_2 + 1}.
\ee

}

\hide{

\begin{proposition} 
\label{identifiability}
Consider a full or non-full, curved exponential family with countable support $\mbX$ and local dependence.
Suppose that $\bta^\star \in \etaspace \subseteq \interior(\fullspace)$.
Then,
for all $\bta \in \etaspace \subseteq \interior(\fullspace)$,
there exists $\dot\bta = \lambda \; \bta^\star + (1 - \lambda) \; \bta \in \interior(\fullspace)$ $(0 < \lambda < 1)$ such that
\be
\nonumber
KL(\bta^\star;\, \bta)
&\geq& \dfrac{2\; \norm{\bmu(\bta^\star) - \bmu(\dot\bta)}_2\; \norm{\bta^\star - \bta}_2^{3 / 2}}{\norm{\bta^\star - \bta}_2 + 1},
\ee
where $KL(\bta^\star;\, \bta)$ denotes the Kullback-Leibler divergence from $\mbP_{\bta^\star}$ to $\mbP_{\bta}$.
\end{proposition}

\com 
In non-full exponential families with $\etaspace \subset \interior(\fullspace)$---e.g., 
curved exponential families---$\dot\bta$ may not be an element of $\etaspace$.
However,
as long as $\bta^\star \in \etaspace \subseteq \interior(\fullspace)$ and $\bta \in \etaspace \subseteq \interior(\fullspace)$,
$\dot\bta$ is an element of $\interior(\fullspace)$,
because the natural parameter space $\fullspace$ is convex \citep[][Theorem 1.13, p.\ 19]{Br86}.
Therefore,
even when $\dot\bta \not\in \etaspace$,
we are assured that $\dot\bta \in \interior(\fullspace)$ and hence $\psi(\dot\bta)\, <\, \infty$.

\pproof \ref{identifiability}.
Since $g(\bta; \, \bmu(\bta^\star))$ is the expected loglikelihood function,
we have
\be
\nonumber
KL(\bta^\star;\, \bta)
&=& g(\bta^\star; \, \bmu(\bta^\star)) - g(\bta; \, \bmu(\bta^\star)).
\ee
By Lemmas \ref{lemma1} and \ref{corollary.reverse.cs} along with $\bta^\star \in \etaspace \subseteq \interior(\fullspace)$ and $\bta \in \etaspace \subseteq \interior(\fullspace)$,
there exists $\dot\bta = \lambda\; \bta^\star + (1 - \lambda)\; \bta \in \interior(\fullspace)$ $(0 < \lambda < 1)$ such that
\be
\nonumber
g(\bta^\star; \, \bmu(\bta^\star)) - g(\bta; \, \bmu(\bta^\star))
&=& \langle \bta^\star - \bta, \; \bmu(\bta^\star) - \bmu(\dot\bta)\rangle\s
\\
&\geq& \dfrac{2 \; \norm{\bmu(\bta^\star) - \bmu(\dot\bta)}_2\; \norm{\bta^\star - \bta}_2^{3 / 2}}{\norm{\bta^\star - \bta}_2 + 1}.
\ee

}

\hide{
\begin{lemma}
\label{lemma:geometric_param_properties}
Consider a non-full, curved exponential family with natural parameter vector $\bta(\btheta)$,
where the coordinates $\eta_i(\btheta)$ of $\bta(\btheta)$ are given by
\be
\nonumber
\eta_i(\btheta)
& = & \theta_1 \left[\theta_2 - \theta_2 \, \left(1 - \dfrac{1}{\theta_2}\right)^i\right],
&& i = 1, \dots, I,
&& I \geq 2.
\ee
Let $\bthetas \in \bTheta$ be the data-generating parameter vector,
where $\bTheta = \mR^+ \times (C,\, \infty)$ and $C > 1/2$.
Then,
for all $\epsilon > 0$ and all $\btheta \in \bTheta$ such that $\norm{\bthetas - \btheta}_2 \geq \epsilon$,
there exists $\delta(\epsilon) > 0$ such that
\be
\nonumber
\norm{\bta(\btheta^{\star}) - \bta(\btheta)}_2
&\geq& \delta(\epsilon).
\ee
\end{lemma}
\llproof \ref{lemma:geometric_param_properties}.
Since $I \geq 2$,
the natural parameter vector $\bta(\btheta)$ has at least two coordinates,
denoted by $\eta_1(\btheta)$ and $\eta_2(\btheta)$.
The first coordinate $\eta_1(\btheta)$ is sensitive to changes of $\theta_1$ but is not sensitive to changes of $\theta_2$,
because
\be
\nonumber
\eta_1(\btheta)
&=& \theta_1\, \left[\theta_2 - \theta_2 \left(1 - \dfrac{1}{\theta_2} \right)\right]
&=& \theta_1.
\ee
The second coordinate $\eta_2(\theta)$ of $\bta(\theta)$ is sensitive to changes of both $\theta_1$ and $\theta_2$ as long as $\theta_1 \neq 0$,
because
\be
\nonumber
\eta_{2}(\btheta)
\,=\, \theta_1\, \left[\theta_2 - \theta_2 \left(1 - \dfrac{1}{\theta_2} \right)^2 \right]
\,=\, \theta_1\, \left[\theta_2 - \dfrac{(\theta_2 - 1)^2}{\theta_2} \right]
\,=\, \theta_1\, (2 - \theta_2^{-1}).
\ee
Therefore,
as long as $I \geq 2$,
\be 
\nonumber 
\norm{\bta(\bthetas) - \bta(\btheta)}_2 
&\geq& \left|\eta_2(\bthetas) - \eta_2(\btheta)\right|.
\ee 
By the mean-value theorem, 
there exists $\dot\btheta = \lambda\; \bthetas + (1-\lambda)\; \btheta$ $(0 < \lambda < 1)$ such that
\be 
\nonumber 
\left|\eta_2(\bthetas) - \eta_2(\btheta)\right| 
&=& \left|\langle \nabla_{\btheta}\, \eta_2(\btheta)\vert_{\btheta = \dot\btheta},\, \bthetas - \btheta\rangle\right|,
\ee 
where
\beno
\nabla_{\theta_1}\, \eta_2(\btheta)
&=& 2 - \theta_2^{-1}
&&\text{and}&&
\nabla_{\theta_2}\, \eta_2(\btheta)
&=& \theta_1\, \theta_2^{-2}.
\ee
We lower bound the inner product $\left|\langle \nabla_{\btheta}\, \eta_2(\btheta)\vert_{\btheta = \dot\btheta},\, \bthetas - \btheta\rangle\right|$ 
by using the same argument we used in the proof of Lemma \ref{corollary.reverse.cs}.
To do so,
let $\bm{v}_1 = \nabla_{\btheta}\, \eta_2(\btheta) \vert_{\btheta = \dot\btheta}$ and $\bm{v}_2 = \bthetas - \btheta$.
By the same argument used in the proof of Lemma \ref{corollary.reverse.cs},
the condition \eqref{condition.reverse.cs} of Lemma \ref{corollary.reverse.cs} is satisfied as long as there exist $A > 0$ and $B > 0$ such that 
\be 
\nonumber
A + B
&\geq& \dfrac{\norm{\bv_1}_2^2 + A\, B\, \norm{\bv_2}_2^2}{\norm{\bv_1}_2\, (1 + \norm{\bv_2}_2)}
&>& 0.
\ee
Choose $A = \norm{\bm{v}_1}_2 \, / \, \norm{\bm{v}_2}_2$ and $B = \norm{\bm{v}_1}_2$.
Then 
\be 
\nonumber 
A + B
&\geq& \dfrac{\norm{\bv_1}_2^2 + A\, B\, \norm{\bv_2}_2^2}{\norm{\bv_1}_2\, (1 + \norm{\bv_2}_2)}
&=& \norm{\bm{v}_1}_2
&>& 0,
\ee 
where $\norm{\bm{v}_1}_2 = \norm{\nabla_{\btheta} \, \eta_2(\btheta) \vert_{\btheta = \dot\btheta}}_2 > 0$,
because $\nabla_{\theta_1}\, \eta_2(\btheta) = 2\, \theta_2 - 1 > 0$ for all $\theta_2 \in (C, \, \infty)$ with $C > 1/2$.
Therefore,
the condition \eqref{condition.reverse.cs} of Lemma \ref{reverse.cs} is satisfied and Lemma \ref{reverse.cs} can be applied to obtain
\be 
\label{eq:eta_diff_lower_1}
\left|\left\langle \nabla_{\btheta}\, \eta_2(\btheta) \vert_{\btheta = \dot\btheta},\; \bthetas - \btheta\right\rangle \right| 
&\geq& 2 \, \dfrac{\norm{\bthetas - \btheta}_2 \, \norm{\nabla_{\btheta} \, \eta_2(\btheta) \vert_{\btheta = \dot\btheta}}_2^{3/2}}{\norm{\nabla_{\btheta} \, \eta_2(\btheta) \vert_{\btheta = \dot\btheta}}_2 + 1}\s
\\
&\geq& 2 \, \dfrac{\epsilon \, \norm{\nabla_{\btheta} \, \eta_2(\btheta) \vert_{\btheta = \dot\btheta}}_2^{3/2}}{\norm{\nabla_{\btheta} \, \eta_2(\btheta) \vert_{\btheta = \dot\btheta}}_2 + 1}.
\ee 
Lemma \ref{lemma:f} states that the lower bound in \eqref{eq:eta_diff_lower_1} is strictly increasing in $\norm{\nabla_{\btheta}\, \eta_2(\btheta) \vert_{\btheta = \dot\btheta}}_2$.
Thus, 
we can lower bound \eqref{eq:eta_diff_lower_1} by lower bounding $\norm{\nabla_{\btheta}\, \eta_2(\btheta) \vert_{\btheta = \dot\btheta}}_2$.
By definition,
\be 
\nonumber 
\norm{\nabla_{\btheta}\, \eta_2(\btheta) \vert_{\btheta = \dot\btheta}}_2
\hide{
&=& \sqrt{ \left( \nabla_{\theta_1}\, \eta_2(\btheta) \vert_{\btheta = \dot\btheta} \right)^2 + \left(  \nabla_{\theta_2}\, \eta_2(\btheta) \vert_{\btheta = \dot\btheta} \right)^2}\s
\\
}
&=& \sqrt{\left(2 - \dot\theta_2^{-1}\right)^2 + \left(\dot\theta_1\, \dot\theta_2^{-2}\right)^2}
&\geq& \left|2 - \dot\theta_2^{-1}\right|,
\ee
where $|2 - \dot\theta_2^{-1}| = 2 - \dot\theta_2^{-1} \geq 2 - C^{-1} > 0$ for all $\dot\theta_2 \in (C, \infty)$ with $C > 1/2$.
Therefore,
there exists $D > 0$ such that $\norm{\nabla_{\btheta}\, \eta_2(\btheta) \vert_{\btheta = \dot\btheta}}_2 > D$.
As a result,
by Lemma \ref{lemma:f},
for all $\epsilon > 0$ and all $\btheta \in \bTheta$ such that $\norm{\bthetas - \btheta}_2 \geq \epsilon$,
\be 
\nonumber 
\left|\left\langle \nabla_{\btheta}\, \eta_2(\btheta) \vert_{\btheta = \dot\btheta},\; \bthetas - \btheta\right\rangle \right| 
&\geq& 2\, \dfrac{\epsilon \, \norm{\nabla_{\btheta} \, \eta_2(\btheta) \vert_{\btheta = \dot\btheta}}_2^{3/2}}{\norm{\nabla_{\btheta} \, \eta_2(\btheta) \vert_{\btheta = \dot\btheta}}_2 + 1}
&\geq& 2\, \dfrac{\epsilon\; D^{3/2}}{D + 1}.
\ee 
Collecting terms shows that,
for all $\epsilon > 0$ and all $\btheta \in \bTheta$ such that $\norm{\bthetas - \btheta}_2 \geq \epsilon$,
there exists $\delta(\epsilon) > 0$ such that
\be
\nonumber
\norm{\bta(\bthetas) - \bta(\btheta)}_2
\,\geq\, \left|\eta_2(\bthetas) - \eta_2(\btheta)\right|
\,\geq\, \left|\left\langle \nabla_{\btheta}\, \eta_2(\btheta) \vert_{\btheta = \dot\btheta},\; \bthetas - \btheta\right\rangle \right|
\,\geq\, \delta(\epsilon),
\ee
where $\delta(\epsilon) = 2\, \epsilon\, D^{3/2} / (D + 1) > 0$.

\begin{lemma} 
\label{eta_max} 
Consider a non-full, curved exponential family with natural parameter vector $\bta(\btheta)$,
where the coordinates $\eta_i(\btheta)$ of $\bta(\btheta)$ are given by
\be
\nonumber
\eta_i(\btheta)
& = & \theta_1 \left[\theta_2 - \theta_2 \, \left(1 - \dfrac{1}{\theta_2}\right)^i\right],
&& i = 1, \dots, I,
&& I \geq 2.
\ee
Then,
for all $\btheta \in \mbR^{+} \times (1/2, \, \infty)$,
\be 
\nonumber 
\max\limits_{1 \leq i \leq I} \, \eta_i(\btheta)
&=& \max\left(\theta_1,\, \theta_1\, \theta_2\right).
\ee 
\end{lemma} 
\llproof \ref{eta_max}.
Observe that
\be 
\nonumber
\eta_{i}(\btheta)
\;=\; \theta_1 \left[\theta_2 - \theta_2 \left(1 - \dfrac{1}{\theta_2}\right)^{i}\right]
\;=\; \theta_1\, \theta_2\, \left[1 - \left(1 - \dfrac{1}{\theta_2}\right)^{i}\right],
\;\; i = 1, \dots, I,
\ee
so that $\max_{1 \leq i \leq I} \eta_i(\btheta)$ is determined by $(1 - 1/\theta_2)^i$.
We distinguish three cases.

\s

{\em Case $\theta_2 \in (1/2, 1)$:}
If $\theta_2 \in (1/2, 1)$,
then $(1 - 1/\theta_2)^i$ is negative for all odd integers $i \in \{1, \dots, I\}$ and positive for all even integers $i \in \{1, \dots, I\}$.
Hence $\max_{1 \leq i \leq I} \eta_i(\btheta)$ is determined by the smallest odd integer $i \in \{1, \dots, I\}$, 
implying
\be
\nonumber
\max\limits_{1 \leq i \leq I} \eta_i(\btheta)
&=& \eta_{1}(\btheta)
&=& \theta_1.
\ee
 
{\em Case $\theta_2 = 1$:}
If $\theta_2 = 1$,
then $\eta_i(\btheta) = \theta_1$ for all integers $i \in \{1, \dots, I\}$,
implying 
\be
\nonumber
\max\limits_{1 \leq i \leq I} \eta_i(\btheta) 
&=& \theta_1.
\ee

{\em Case $\theta_2 \in [1, \infty)$:}
If $\theta_2 \geq 1$,
then $0 \leq (1 - 1/\theta_2)^i < 1$ 
for all integers $i \in \{1, \dots, I\}$,
which implies
\be
\nonumber
\max\limits_{1 \leq i \leq I} \eta_i(\btheta) 
&=& \max\limits_{1 \leq i \leq I} \theta_1\, \theta_2\, \left[1 - \left(1 - \dfrac{1}{\theta_2}\right)^I\right]
&\leq& \theta_1\, \theta_2.
\ee

\noindent
As a result,
for all $\btheta \in \mR^+ \times (1/2, \infty)$,
\be
\nonumber
\max\limits_{1 \leq i \leq I} \, \eta_i(\btheta)
&=& \max\left(\theta_1,\, \theta_1\, \theta_2\right).
\ee

\begin{lemma} 
\label{edge.transitive} 
Consider exponential families with support $\mbX = \{0, 1\}^{\sum_{k=1}^K \binom{|\mA_k|}{2}}$ and local dependence generated by the number of within-neighborhood edges and transitive edges:
\be 
\nonumber 
s_1(\bx) 
&=& \dsum_{k=1}^K\; \dsum_{i\in\mA_k\, <\, j\in\mA_k} x_{i,j}\s
\\
s_2(\bx) 
&=& \dsum_{k=1}^K\; \dsum_{i\in\mA_k\, <\, j\in\mA_k} x_{i,j}\, \max\limits_{h \in \mA_k,\, h \neq i,j} x_{i,h} \; x_{j,h}.
\ee 
Let $\bta(\btheta) = \btheta$ and $\bTheta = \mR \times \mR^{+}$.
Then,
for all $\epsilon > 0$ and all $\btheta \in \bTheta$ such that $\norm{\bthetas - \btheta}_2 \geq \epsilon$,
\be
\nonumber
g(\bthetas; \, \bmu(\bta(\bthetas))) - g(\btheta; \, \bmu(\bta(\bthetas)))
&\geq & \delta(\epsilon)\, \dsum_{k=1}^{K}\, \dis\binom{|\mA_k|}{2}.
\ee
\end{lemma}

\llproof \ref{edge.transitive}.
To ease the presentation,
we consider $K = 1$ neighborhood and drop the subscript $k$ from all neighborhood-dependent quantities.
The extension to $K \geq 2$ neighborhoods is straightforward.
By Proposition \ref{identifiability},
for all $\epsilon > 0$ and all $\btheta \in \bTheta$ such that $\norm{\bthetas - \btheta}_2 \geq \epsilon$,
there exists $\dot\btheta = \lambda \; \bthetas + (1 - \lambda) \; \btheta \in \bTheta$ ($0 < \lambda < 1$) such that
\be 
\label{eq:edge_tran_objective}
g(\bthetas;\, \bmu(\bta(\bthetas))) - g(\btheta;\, \bmu(\bta(\bthetas)))
&\geq& \dfrac{2\; \norm{\bthetas - \btheta}_2^{3 \, / \, 2}\; \norm{\bmu(\bta(\bthetas)) - \bmu(\dot\btheta)}_2}{\norm{\bthetas - \btheta}_2 + 1}.
\ee 
We bound the right-hand side of \eqref{eq:edge_tran_objective} with respect to $\norm{\bthetas - \btheta}_2$ as follows.
Lemma \ref{lemma:f} implies that, 
for any $\norm{\bmu(\bta(\bthetas)) - \bmu(\dot\btheta)}_2$,
the right-hand side of \eqref{eq:edge_tran_objective} is minimized by the smallest possible value of $\norm{\bthetas - \btheta}_2$,
which implies that 
\be
\label{LB}
g(\bthetas;\, \bmu(\bta(\bthetas))) - g(\btheta;\, \bmu(\bta(\bthetas)))
&\geq& \dfrac{2\; \epsilon^{3 \, / \, 2}\; \norm{\bmu(\bta(\bthetas)) - \bmu(\dot\btheta)}_2}{\epsilon + 1}.
\ee
It remains to bound $\norm{\bmu(\bta(\bthetas)) - \bmu(\dot\btheta)}_2$ from below.
To do so,
first note that the natural parameter vectors $\bthetas \in \bTheta  \subset \mR^2$ and $\dot\btheta \in \bTheta \subset \mR^2$ are distinct and the map $\bmu: \mR^2 \mapsto \rint(\mathbb{M})$ is one-to-one \citep[][Theorem 3.6, p.\ 74]{Br86}.
Therefore, 
the mean-value parameter vectors $\bmu(\bta(\bthetas))$ and $\bmu(\dot\btheta)$ are distinct and $\norm{\bmu(\bta(\bthetas)) - \bmu(\dot\btheta)}_2 > 0$.
To bound $\norm{\bmu(\bta(\bthetas)) - \bmu(\dot\btheta)}_2$ from below,
note that
\be 
\nonumber 
\norm{\bmu(\bta(\bthetas)) - \bmu(\dot\btheta)}_2 
&=&\sqrt{(\mu_1(\bta(\bthetas)) - \mu_1(\dot\btheta))^2 +  (\mu_2(\bta(\bthetas)) - \mu_2(\dot\btheta))^2}
&>& 0.
\ee 
We bound both deviations $|\mu_1(\bta(\bthetas)) - \mu_1(\dot\btheta)|$ and $|\mu_2(\bta(\bthetas)) - \mu_2(\dot\btheta)|$;
note that at least one of them must be non-zero.
The first deviation $|\mu_1(\bta(\bthetas)) - \mu_1(\dot\btheta)|$ can be bounded from below by using Lemma \ref{mean.bound}, 
which shows that there exists $C_1(\bthetas) \in (0, 1)$ such that $0 < \Lambda(\bthetas) \leq C_1(\bthetas) \leq \Omega(\bthetas) < 1$ and
\be 
\nonumber 
\mu_1(\bta(\bthetas))
&=& \mbE_{\bthetas} \dsum_{i\in\mA\, <\, j\in\mA} X_{i,j}
&=& C_1(\bthetas)\, \dis\binom{|\mA|}{2},
\ee 
where $\Lambda(\bthetas)$ and $\Omega(\bthetas)$ are given by
\be 
\nonumber 
\Lambda(\bthetas) 
&=& \left[1 + \exp\left( - \theta_1^\star\right) \right]^{-1}
&>& 0
\ee 
and 
\be \nonumber 
\Omega(\bthetas)
&=& \left[1 + \exp\left( - \theta_1^\star - \theta_2^\star\, (2 \, |\mA| - 3)\right) \right]^{-1}
&<& 1,
\ee
respectively.
Along the same lines, 
there exists $C_1(\dot\btheta) \in (0, 1)$ such that $0 < \omega(\dot\btheta) \leq C_1(\dot\btheta) \leq \alpha(\dot\btheta) < 1$ and
\be 
\nonumber 
\mu_1(\dot\btheta)
&=& C_1(\dot\btheta)\, \dis\binom{|\mA|}{2},
\ee 
which implies
\be 
\nonumber 
|\mu_1(\bta(\bthetas)) - \mu_1(\dot\btheta)|
&=& |C_1(\bthetas) - C_1(\dot\btheta)| \, \dis\binom{|\mA|}{2}.
\ee
In addition,
by Lemma \ref{tran_edge_alpha_exp},
there exists $C_2(\bthetas) > 0$ such that
\be 
\nonumber 	
\mu_2(\bta(\bthetas))
&=& \mbE_{\bthetas}\, \dsum_{i\in\mA\, <\, j\in\mA} X_{i,j}\, \max\limits_{h \in \mA,\, h \neq i,j} X_{i,h} \; X_{j,h}
&=& C_2(\bthetas) \, \dis\binom{|\mA|}{2},
\ee
where $C_2(\bthetas)$ satisfies 
\be
\nonumber
\lambda\, \Lambda(\bthetas)^3 + (1 - \lambda)\, \Omega(\bthetas)^3
&\leq& C_2(\bthetas)
&\leq& \lambda\, \Lambda(\bthetas) + (1 - \lambda)\, \Omega(\bthetas)
\ee
for some $\lambda \in (0, 1)$.
The lower bound satisfies
\be
\nonumber
\lambda\, \Lambda(\bthetas)^3 + (1 - \lambda)\, \Omega(\bthetas)^3 
\geq \lambda\, \left[1 + \exp\left( - \theta_1^\star\right)\right]^{-3} + (1 - \lambda)\, \left[1 + \exp\left( - \theta_1^\star - 3\, \theta_2^\star\right)\right]^{-3}
\ee
by using $\theta_2^\star > 0$ and $|A| \geq 3$.
The upper bound satisfies
\be
\nonumber
\lambda\, \Lambda(\bthetas) + (1 - \lambda)\, \Omega(\bthetas) 
&\leq& \lambda\, \left[1 + \exp\left( - \theta_1^\star\right)\right]^{-1} + (1 - \lambda).
\ee
Thus,
$C_2(\bthetas)$ is contained in a bounded interval whose endpoints do not depend on $|\mA|$. 
Along the same lines, 
it can be shown that there exists $C_2(\dot\btheta) > 0$ such that
\be 
\nonumber       
\mu_2(\dot\btheta)
&=& C_2(\dot\btheta) \, \dis\binom{|\mA|}{2},
\ee
where $C_2(\dot\btheta)$ is contained in a bounded interval whose endpoints do not depend on $|\mA|$.
As a result,
\be 
\nonumber 
|\mu_2(\bta(\bthetas)) - \mu_2(\dot\btheta)|
&=& |C_2(\bthetas) - C_2(\dot\btheta)|\, \dis\binom{|\mA|}{2}.
\ee  
Since the exponential family is minimal,
$\bthetas \neq \dot\btheta$ implies that either $|\mu_1(\bta(\bthetas)) - \mu_1(\dot\btheta)| > 0$ or $|\mu_2(\bta(\bthetas)) - \mu_2(\dot\btheta)| > 0$ or both,
hence there exists $\gamma(\epsilon) > 0$ such that 
\be 
\nonumber
\norm{\bmu(\bta(\bthetas)) - \mu(\dot\btheta)}_2 
= \sqrt{(C_1(\bthetas) - C_1(\dot\btheta))^2 + (C_2(\bthetas) - C_2(\dot\btheta))^2}\, \dis\binom{|\mA|}{2}
\geq \gamma(\epsilon) \dis\binom{|\mA|}{2},
\ee 
where $\gamma(\epsilon)$ 
does not depend on $|A|$.
Using \eqref{LB},
we conclude that,
for all $\epsilon > 0$ and all $\btheta \in \bTheta$ such that $\norm{\bthetas - \btheta}_2 \geq \epsilon$,
there exists $\delta(\epsilon) > 0$ such that
\be 
\nonumber 
g(\bthetas;\, \bmu(\bta(\bthetas))) - g(\btheta; \, \bmu(\bta(\bthetas)))
&\geq& \dfrac{2\; \epsilon^{3 \, / \, 2}\; \gamma(\epsilon)}{\epsilon + 1}\, \dis\binom{|\mA|}{2} 
&=& \delta(\epsilon) \dis\binom{|\mA|}{2},
\ee 
where $\delta(\epsilon) = 2\; \epsilon^{3 \, /  \, 2}\; \gamma(\epsilon)\,/\, (\epsilon + 1) > 0$.

\hide{

\begin{lemma}
\label{lemma:f}
The function $f: \mR^+ \mapsto \mR^+$ defined by
\be
\nonumber
f(z)
&=& \dfrac{2\; C\; z^{3 / 2}}{z + 1},
&& C > 0, && z \in \mR^+
\ee
is strictly increasing for all $z > 0$.
Thus,
for all $\epsilon > 0$,
\be
\nonumber
\min\limits_{z\, \geq\, \epsilon} f(z)
&=& f(\epsilon).
\ee
\end{lemma}
\llproof \ref{lemma:f}.
The function $f: \mR^+ \mapsto \mR^+$ is continuously differentiable on $\mR^+$ and its derivative is given by
\be
\nonumber
\nabla_{z}\, f(z)
\hide{
&=& \nabla_{z}\, (2\; C\; z^{3 / 2})\, (z + 1)^{-1}
&=& 2\, C \, \dfrac{3}{2} \, z^{1/2} \, (z + 1)^{-1} \, - 2\, C \, z^{3 / 2} \, (z + 1)^{-2}
}
\hide{
&=& C \, \dfrac{3 \, z^{1 / 2} \, (z + 1) - 2 \, z^{3 / 2}}{(z+1)^2}
&=& C \, \dfrac{3 \, z^{3 / 2} + 3 \, z^{1/ 2} - 2 \, z^{3 / 2}}{(z+1)^2}\s
\\
&=& C \, \dfrac{z^{3 / 2} + 3\,z^{1 / 2}}{(z + 1)^2}
}
&=& \dfrac{C\; z^{1/2}\; (z + 3)}{(z + 1)^2}
\;\;>\;\; 0,
\;\;\;\; z \in \mR^+.
\ee
Since $\nabla_{z} \, f(z) > 0$ for all $z > 0$,
$f(z)$ is strictly increasing for all $z > 0$ and hence,
for all $\epsilon > 0$,
\be
\nonumber
\min\limits_{z\, \geq\, \epsilon} f(z)
&=& f(\epsilon).
\ee

}

\s

\begin{lemma}
\label{c.curved}
Consider the curved exponential family with support\linebreak 
$\mbX = \{0, 1\}^{\sum_{k=1}^{K} \binom{|\mA_k|}{2}}$ generated by the number of within-neighborhood geometrically weighted edgewise shared partners.
Let $\Thetad = [A, B] \subseteq \Theta$ be a subset of the parameter space $\Theta$,
where $1\, /\, 2 < A < B$.
Then conditions [C.1]---[C.4] are satisfied.
\end{lemma}

\llproof \ref{c.curved}.
To show that conditions [C.1]---[C.4] hold,
note that the coordinates $\eta_{k,i}(\theta)$ of the neighborhood-dependent natural parameter vectors $\bta_k(\theta)$ can be written as 
\be
\label{reduced.form}
\eta_{k,i}(\theta)
= \theta\, \left[1 - \left(1 - \dfrac{1}{\theta}\right)^i\right],\; i = 1, \dots, I_k,\; I_k \geq 2,\; k = 1, \dots, K.
\ee
Observe that exponential families with natural parameters of the form \eqref{reduced.form} can be reduced to an exponential family with natural parameter vector $\bta(\theta)$ of dimension
$\ppp = \max_{1 \leq k \leq K} I_k$.
We denote the coordinates of the natural parameter vector $\bta(\theta)$ by $\eta_i(\theta)$,
which are given by
\beno
\eta_i(\theta)
\= \theta - \theta\, \alpha(\theta)^i,
&& i = 1, \dots, \max\limits_{1 \leq k \leq K} I_k,
\ee
where
\beno
\alpha(\theta)
\= 1 - \dfrac{1}{\theta}.
\ee
The coordinates $\eta_{k,i}(\theta)$ of the neighborhood-dependent natural parameter vectors $\bta_k(\theta)$ are related to the coordinates $\eta_i(\theta)$ of the natural parameter vector $\bta(\theta)$ of the exponential family as follows:
\beno
\eta_{k,i}(\theta)
\= \eta_i(\theta),
&& i = 1, \dots, I_k, 
&& k = 1, \dots, K.
\ee
A helpful observation is that the coordinates $\eta_i(\theta)$ of $\bta(\theta)$ are continuously differentiable on $(A, B)$ with derivatives
\beno
\nabla_\theta\, \eta_i(\theta)
\hide{
\= \nabla_\theta \left[\theta - \theta\, \alpha(\theta)^i\right]
}
\= 1 - \alpha(\theta)^i - \dfrac{i}{\theta}\, \alpha(\theta)^{i-1},
&& \theta \in (A, B).
\ee
We check conditions [C.1]---[C.4] one by one.

\s

\underline{Condition [C.1].}
To show that the map $\bta: \Thetad \mapsto \etaspace$ is one-to-one on $\Thetad$ for all $K \geq 1$,
we show that at least one coordinate of $\bta(\theta+\delta)$ must deviate from $\bta(\theta)$ for all $\theta \in (A, B)$,
all $\delta > 0$,
and all $K \geq 1$.
To do so,
note that $\bta(\theta)$ has at least two coordinates,
denoted by $\eta_1(\theta)$ and $\eta_2(\theta)$,
because $I_k \geq 2$ ($k = 1, \dots, K$).
The first coordinate $\eta_1(\theta)$ of $\bta(\theta)$ is constant on $(A, B)$:
\beno
\eta_1(\theta)
\= 1,
&& \theta \in (A, B).
\ee
The second coordinate $\eta_2(\theta)$ of $\bta(\theta)$ is continuously differentiable on $(A, B)$ with derivative
\beno
\nabla_\theta\, \eta_2(\theta)
\= 1 - \alpha(\theta)^2 - \dfrac{2}{\theta}\, \alpha(\theta)
\= \dfrac{1}{\theta^2}
&>& 0,
&& \theta \in (A, B).
\ee
By the mean-value theorem,
\beno
\eta_2(\theta + \delta) - \eta_2(\theta)
\gte \dfrac{\delta}{(\theta + \delta)^2}
&>& 0,
&& \theta \in (A, B),
&& \delta > 0.
\ee
Thus,
$\eta_2(\theta)$ is strictly increasing on $(A, B)$ and at least one coordinate of $\bta(\theta+\delta)$ must deviate from $\bta(\theta)$ for all $\theta \in (A, B)$,
all $\delta > 0$, 
and all $K \geq 1$. 
As a result,
the map $\bta: \Thetad \mapsto \etaspace$ is one-to-one and continuous on $\Thetad$.
Thus condition [C.1] is satisfied.

\s

\underline{Condition [C.2].}
Condition [C.2] follows from the continuity of $\bta: \Thetad \mapsto \etaspace$ and the upper semicontinuity of exponential-family loglikelihood functions \citep[][Lemma 5.3, p.\ 146]{Br86} along with the fact that $\bta: \Thetad \mapsto \etaspace$ is one-to-one and exponential-family loglikelihood functions are strictly concave on $\etaspace \subseteq \interior(\fullspace)$ \citep[][Lemma 5.3, p.\ 146]{Br86}.

\s

\underline{Condition [C.3].}
Choose any $\theta \in \Thetad$ and $\theta^\prime \in \Thetad$.
By the triangle inequality,
we obtain,
for all $\bmu \in \mM$,
\be
\label{ccurved1}
\left|\langle\bta(\theta^\prime) - \bta(\theta),\, \bmu\rangle\right|
\hide{
\= \left|\dsum_{k=1}^K \langle\bta_k(\theta^\prime) - \bta_k(\theta),\, \bmu_k\rangle\right|
}
\= \left|\dsum_{k=1}^K \dsum_{i=1}^{I_k} \left[\eta_{k,i}(\theta^\prime) - \eta_{k,i}(\theta)\right] \mu_{k,i}\right|\s
\\
\hide{
\= \left|\dsum_{i = 1}^{I_k} \left[\eta_i(\theta^\prime) - \eta_i(\theta)\right] \dsum_{k=1}^K \mu_{k,i}\right|
}
\lte \dsum_{k=1}^K \dsum_{i=1}^{I_k} |\eta_{k,i}(\theta^\prime) - \eta_{k,i}(\theta)| |\mu_{k,i}|.
\ee
We show that there exists $C > 2$ such that,
for all $\theta \in \Thetad$ and all $i \in \{1, 2, \dots\}$,
\beno
\left|\nabla_\theta\, \eta_i(\theta)\right|
\lte \max(3, C),
&& i \in \{1, 2, \dots\},
\ee
which,
by the mean-value theorem,
implies that
\beno
|\eta_i(\theta^\prime) - \eta_i(\theta)|
\hide{
\= |\theta^\prime - \theta|\, |\nabla_\theta\, \eta_i(\theta)|_{\theta=\dot\theta}|
}
\lte |\theta^\prime - \theta|\, \max(3, C),
&& i \in \{1, 2, \dots\}.
\ee
To show that $|\nabla_\theta\, \eta_i(\theta)|$ is bounded for all $\theta \in \Thetad$ and all $i \in \{1, 2, \dots\}$,
note that
\beno
\left|\nabla_\theta\, \eta_i(\theta)\right|
\hide{
\= \left|\nabla_\theta\, \left(\theta - \theta\, \alpha(\theta)^i\right)\right|
}
\= \left|1 - \alpha(\theta)^i - \dfrac{i}{\theta}\, \alpha(\theta)^{i-1}\right|
\lte 2 + 2\, i\, |\alpha(\theta)|^{i-1},
\ee
where we used the fact that $\theta > 1/2$ and $|\alpha(\theta)| < 1$ for all $\theta \in \Thetad$.
Observe that
\be
\label{bound}
\left|\nabla_\theta\, \eta_i(\theta)\right|
\lte 2 + 2\, i\, |\alpha(\theta)|^{i-1}
\lte 3
\ee
is satisfied as long as
\be
\label{bound.exponential}
|\alpha(\theta)|
\lte \exp\left(-\dfrac{\log 2\, i}{i-1}\right),
&& i \in \{2, 3, \dots\}.
\ee
We note that \eqref{bound} does not hold for small $i$ and while it does hold for sufficiently large $i$,
the upper bound $3$ in \eqref{bound} is not tight, 
but it is convenient.
To show that \eqref{bound} is satisfied by all $\theta\in\Thetad$ for all sufficiently large $i$,
we lower bound the exponential term on the right-hand side of \eqref{bound.exponential} as follows.
Choose $\gamma > 0$ so that $0 < \gamma \leq \min((2\, A - 1) / 2,\, 1/B)$,
where $A > 1/2$ and $B > 1/2$ were selected above so that $1 / 2 < A \leq \theta \leq B < \infty$ for all $\theta \in \Thetad$;
note that the choice of $\gamma > 0$ implies that $1 / (2 - \gamma) \leq A < B \leq 1 / \gamma$.
\hide{
Example:
> A <- .51
> B <- 51
> (2*A - 1) / A
[1] 0.03921569
> 1 / B
[1] 0.01960784
}
Then there exists $I_0 \geq 2$ such that,
for all $i > I_0$,
\beno
1-\gamma
\lte \exp\left(-\dfrac{\log 2\, i}{i-1}\right),
&& i \in \{I_0+1, I_0+2, \dots\},
&& I_0 \geq 2.
\ee
In other words,
we have $|\nabla_\theta\, \eta_i(\theta)| \leq 3$ as long as $i > I_0$ and $|\alpha(\theta)| \leq 1-\gamma$,
i.e.,
as long as $\theta \in \Thetad$ satisfies $1 / (2 - \gamma) \leq \theta \leq 1 / \gamma$. 
\hide{
First, observe that
|beta(theta)| 
= 1 - 1/theta if theta >= 1
= 1/theta - 1 if 1/2 < theta < 1
Case 1/2 < theta < 1:
beta(theta) = 1 - 1/theta < 0
thus   is negative and hence |beta(theta)| = 1/theta - 1, implying
1/theta - 1 < 1 - gamma
1 - theta < theta - gamma theta
1 < 2 theta - gamma theta
1 < theta (2 - gamma)
1 / (2 - gamma) < theta
or
theta > 1 / (2 - gamma) > 1 / 2
e.g.
> 1/(2-.1)
[1] 0.5263158
Case theta >= 1:
1 - 1/theta < 1 - gamma
theta - 1 < theta - gamma theta
-1 < - gamma theta
1 > gamma theta
theta < 1 / gamma
}
Since $\gamma > 0$ was chosen above so that $1 / (2 - \gamma) \leq A < B \leq 1 / \gamma$,
all $\theta \in \Thetad$ satisfy $1 / (2 - \gamma) \leq A \leq \theta \leq B \leq 1 / \gamma$
and hence $|\nabla_\theta\, \eta_i(\theta)| \leq 3$ holds for all $\theta \in \Thetad$ provided $i > I_0$.
In addition,
there exists $C > 2$ such that,
for all $\theta \in \Thetad$ and all $i \in \{1, 2, \dots, I_0\}$,
\beno
\left|\nabla_\theta\, \eta_i(\theta)\right|
\lte 2 + 2\, i\, |\alpha(\theta)|^{i-1}
\lte 2 + 2\, I_0
&=& C,
&& i \in \{1, 2, \dots, I_0\},
\ee
where we used $|\alpha(\theta)| < 1$ for all $\theta \in \Thetad$.
Therefore,
for all $\theta \in \Thetad$,
\beno
\left|\nabla_\theta\, \eta_i(\theta)\right|
\lte 2 + 2\, i\, |\alpha(\theta)|^{i-1}
\lte \max(3, C),
&& i \in \{1, 2, \dots\},
\ee
implying
\be
\label{eta.diff}
|\eta_i(\theta^\prime) - \eta_i(\theta)|
\hide{
\= |\theta^\prime - \theta|\, |\nabla_\theta\, \eta_i(\theta)|_{\theta=\dot\theta}|
}
\lte |\theta^\prime - \theta|\, \max(3, C),
&& i \in \{1, 2, \dots\}.
\ee
Using \eqref{ccurved1} and \eqref{eta.diff} shows that there exists $C_1 > 0$ such that
\be
\label{ccurved0}
&& \left|\langle\bta(\theta^\prime) - \bta(\theta),\, \bmu\rangle\right|
\;\leq\; \dsum_{k=1}^K \dsum_{i=1}^{I_k} |\eta_{k,i}(\theta^\prime) - \eta_{k,i}(\theta)| |\mu_{k,i}|\s
\\
\lte |\theta^\prime - \theta|\, \max(3, C) \dsum_{k=1}^K \dsum_{i=1}^{I_k} |\mu_{k,i}|
\;\leq\; C_1\, \norm{\theta^\prime - \theta}_2 \dsum_{k=1}^K {|\mA_k| \choose 2}^\alpha,
\ee
where $\sum_{i=1}^{I_k} |\mu_{k,i}| = \sum_{i=1}^{I_k} \mbE\, s_{k,i}(\bX)$ is the expected number of transitive edges in neighborhood $\mA_k$,
which is bounded above by ${|\mA_k| \choose 2}$ ($k = 1, \dots, K$). 
Thus condition [C.3] is satisfied.

\s

\underline{Condition [C.4].}
Since $|\alpha(\theta)| < 1$ for all $\theta > 1/2$,
including all $\theta \in \Thetad$,
we have
\beno
\label{bounded.eta}
|\eta_i(\theta)|
\hide{
\= |\theta - \theta\, \alpha(\theta)^i|
}
\lte |\theta| + |\theta|\, |\alpha(\theta)|^i
\lte 2\, B,
&& \theta\in\Thetad.
\ee
which implies that there exist $C_2 > 0$ such that
\beno
\left|\langle\bta(\theta),\, s(\bx_1) - s(\bx_2)\rangle\right|
\;=\; \left|\dsum_{k=1}^K \dsum_{i=1}^{I_k} \eta_{k,i}(\theta)\, \left[s_{k,i}(\bx_{k,1}) - s_{k,i}(\bx_{k,2})\right]\right|\s
\\
\leq\; 2\, B\, \left|\dsum_{k=1}^K \dsum_{i=1}^{I_k} s_{k,i}(\bx_{k,1}) - \dsum_{k=1}^K \dsum_{i=1}^{I_k} s_{k,i}(\bx_{k,2})\right|
\;\leq\;
\hide{
\; C_2\, \norm{\mA}_\infty\, \dsum_{k = 1}^K d(\bx_{k,1}, \bx_{k,2})
\;=\;
}
C_2\; \norm{\mA}_\infty\, d(\bx_1, \bx_2),
\ee
where the last line follows from the fact that $\sum_{i = 1}^{I_k} s_{k,i}(\bx_{k})$ is the number of transitive edges in neighborhood $\mA_k$ and changing a single edge cannot change the number of transitive edges by more than $2\, (|\mA_k| - 2) + 1\, \leq\, 2\, \norm{\mA}_\infty$ ($k = 1, \dots, K$).
Thus condition [C.4] is satisfied.

}

\bibliographystyle{asa1}

\bibliography{base}

\label{last.page}

\end{appendix}

\end{document}

%% file: aospreamble.tex
\usepackage[mathscr]{eucal}
\usepackage{afterpage}
\usepackage{amsmath}
\usepackage{tikz}

\usepackage{centernot}
\usepackage{dlfltxbcodetips,bbm,mathtools,subfigure,amsfonts,amsmath,amssymb,fancybox,graphicx,bm,latexsym}

\usepackage{tikz}

\newcommand{\thetas}{\theta^\star}

\newcommand{\bthetad}{\dot\btheta}

\usetikzlibrary{fit,positioning}

\newcommand{\lip}{{\tiny\mbox{Lip}}}

\newcommand{\uu}{a_1}
\newcommand{\vv}{a_2}

\newcommand{\tv}{{\tiny\mbox{TV}}}

\usepackage[tikz]{bclogo}

\newcommand{\rint}{\mathop{\mbox{rint}}}

\newcommand{\interior}{\mbox{int}}
\newcommand{\bt}{\begin{bclogo}[couleur={rgb:orange,0;yellow,0;white,1},arrondi=0.1,logo=\bcplume,ombre=true]}
\newcommand{\et}{\end{bclogo}\s}
\newcommand{\btt}{\begin{box}}
\newcommand{\ett}{\end{box}}

\newcommand{\btheorem}{\begin{bclogo}[couleur={rgb:orange,0;yellow,0;white,1},arrondi=0.1,logo=\bcplume,ombre=true]{Theorem}}
\newcommand{\ettheorem}{\end{bclogo}}

\newcommand{\bst}{\begin{bclogo}[couleur={rgb:orange,1;yellow,1;white,0.5},arrondi=0.1,logo=\bcpanchant]}
\newcommand{\est}{\end{bclogo}}

\DeclareMathOperator*{\argmin}{arg\,min}
\DeclareMathOperator*{\argmax}{arg\,max}

\newcommand{\benum}{\begin{enumerate}}
\newcommand{\eenum}{\end{enumerate}}

\newcommand{\tends}{\to}

\newcommand{\bq}{\begin{quote}\em}
\newcommand{\eq}{\end{quote}}
\newcommand{\bbq}{\begin{quote}\bf\em}
\newcommand{\ebq}{\end{quote}}

\renewcommand{\=}{&=&}

\newcommand{\lte}{&\leq&}
\newcommand{\gte}{&\geq&}
\newcommand{\mR}{\mathbb{R}}
\newcommand{\mK}{\mathscr{K}}
\newcommand{\mbR}{\mathbb{R}}

\newcommand{\mbX}{\mathbb{X}}

\newcommand{\one}{\mathbbm{1}}

\newcommand{\mbE}{\mathbb{E}}
\newcommand{\mbC}{\mathscr{C}}
\newcommand{\mbD}{\mathscr{D}}
\newcommand{\mbM}{\mathbb{M}}

\newcommand{\mbY}{\mathbb{Y}}

\newcommand{\mbP}{\mathbb{P}}

\newcommand{\hide}[1]{}
\newcommand{\ghost}[1]{}
\newcommand{\ba}{\begin{array}{llllllllll}}
\newcommand{\ea}{\end{array}}
\newcommand{\bea}{\begin{equation}\begin{array}{llllllllll}}
\newcommand{\eea}{\end{array}\end{equation}}
\newcommand{\be}{\begin{equation}\begin{array}{lllllllllllllllll}}
\newcommand{\beno}{\begin{equation}\begin{array}{lllllllllllll}\nonumber}
\newcommand{\ee}{\end{array}\end{equation}}
\newcommand{\bel}{\begin{equation}\begin{array}{lllllllllllll}\nonumber}
\newcommand{\eel}{\Box\end{array}\end{equation}}
\newcommand{\bi}{\begin{itemize}}
\newcommand{\ei}{\end{itemize}}
\newcommand{\ben}{\begin{enumerate}}
\newcommand{\een}{\end{enumerate}}
\newcommand{\alert}{\textcolor{black}}

\newcommand{\dis}{\displaystyle}
\newcommand{\dsum}{\displaystyle\sum\limits}

\newcommand{\dprod}{\displaystyle\prod\limits}

\newcommand{\mM}{\mbox{$\mathbb{M}$}}

\newcommand{\mS}{\mathbb{S}}

\newcommand{\mA}{\mathscr{A}}

\newcommand{\mE}{\mathscr{E}}
\newcommand{\s}{\vspace{0.175cm}}
\newcommand{\bx}{\bm{x}}

\newcommand{\bA}{\bm{A}}

\newcommand{\bv}{\bm{v}}

\newcommand{\bX}{\bm{X}}
\newcommand{\mX}{\mathbb{X}}
\newcommand{\mT}{\mathscr{T}}

\newcommand{\mY}{\mathscr{Y}}

\newcommand{\mB}{\mathscr{B}}

\newcommand{\mL}{\mathscr{L}}
\newcommand{\bY}{\bm{Y}}

\newcommand{\by}{\bm{y}}

\newcommand{\bbeta}{\mbox{\boldmath$\beta$}}

\newcommand{\bmu}{\mbox{\boldmath$\mu$}}

\newcommand{\bTheta}{\bm\Theta}
\newcommand{\bta}{\boldsymbol{\eta}}
\newcommand{\btas}{\bta^\star}

\newcommand{\btheta}{\boldsymbol{\theta}}

\newcounter{comment}
\newenvironment{comment}[1][]{\refstepcounter{comment}\par\noindent%
\textbf{Comment~\thecomment #1}:\\ \rmfamily}{}
\newcommand{\bc}{\begin{comment}\em}
\newcommand{\ec}{\end{comment}}

\newcounter{ex}

\newcounter{counterexample}

\newcounter{definition}

\newcounter{theorem}
\newenvironment{theorem}[1][]{\refstepcounter{theorem}\par\smallskip\indent%
\textbf{Theorem~\thetheorem #1}.\em \rmfamily}{}

\newcounter{proposition}
\newenvironment{proposition}[1][]{\refstepcounter{proposition}\par\smallskip\indent%
\textbf{Proposition~\theproposition #1}.\em \rmfamily}{}

\newcounter{ttproof}
\newcommand{\ttproof}{%
\addtocounter{ttproof}{1}%
{\medskip}\textbf{Proof of Theorem}{}
}

\newcounter{ttproofs}
\newcounter{pproof}
\newcommand{\pproof}{%
\addtocounter{pproof}{1}%
{\medskip}\textbf{Proof of Proposition}{}
}

\newcounter{corollary}
\newenvironment{corollary}[1][]{\refstepcounter{corollary}\par\smallskip\indent%
\textbf{Corollary~\thecorollary #1}.\em \rmfamily}{}

\newcounter{ccproof}
\newcommand{\ccproof}{%
\addtocounter{ccproof}{1}%
{\medskip}\textbf{Proof of Corollary}{}
}

\newcounter{ccproofs}
\newcounter{llproof}
\newcommand{\llproof}{%
\addtocounter{llproof}{1}%
{\smallskip}\textbf{Proof of Lemma}{}
}

\newcounter{lemma}
\newenvironment{lemma}[1][]{\refstepcounter{lemma}\par\medskip\noindent%
\textbf{Lemma~\thelemma #1}.\em \rmfamily}{\medskip}

\newcounter{example}

\newcounter{com}
\newenvironment{com}[1][]{\refstepcounter{com}\par\s\indent%
\textit{Remark~\thecom #1}. \rmfamily}{\medskip}

\newcounter{assumption}